\documentclass[a4paper,11pt]{amsart}
\usepackage[T1]{fontenc} 
\usepackage{amsmath,amsthm}
\usepackage[francais,english]{babel}
\usepackage[dvips]{graphicx}
\usepackage{todonotes}
\usepackage{mathdots} 
\usepackage{color}
\usepackage{tikz}
\usepackage{tikz-cd}
\usetikzlibrary{arrows,positioning,shapes,matrix,backgrounds}
\usetikzlibrary{fit,calc,decorations.pathreplacing}

\usepgflibrary{arrows}
\usepackage{soul}
\usepackage{amsfonts,amssymb}
\usepackage{geometry}
\usepackage{hyperref}
\usepackage{stmaryrd}
\usepackage{fancyhdr}
\usepackage{url}
\usepackage{todonotes}
\usepackage{dsfont}
\usepackage[utf8]{inputenc} 
\usepackage{enumerate}
\usepackage[shortlabels]{enumitem}
\usepackage[all]{xy}
\usepackage{mathtools}
\usepackage{tikz}

\theoremstyle{theorem}

\newtheorem*{thm*}{Theorem}

\newtheorem{prop}{Proposition}[section]

\newtheorem{lemma}[prop]{Lemma}

\newtheorem{thm}[prop]{Theorem}
\newtheorem{definition}[prop]{Definition}

\newtheorem{corollary}[prop]{Corollary}

\newtheoremstyle{pourlesremarques}{\topsep}{\topsep}{\normalfont}{}{\bfseries}{.}{ }{}
\theoremstyle{pourlesremarques}
\newtheorem{rem}[prop]{Remark}
\newtheorem*{rem*}{Remark}
\newtheoremstyle{pourlesexemples}{\topsep}{\topsep}{\normalfont}{}{\bfseries}{.}{ }{}
\theoremstyle{pourlesexemples}

\def\presuper#1#2%
  {\mathop{}%
   \mathopen{\vphantom{#2}}^{#1}%
   \kern-\scriptspace%
   #2}
   
\newcommand{\Wh}{\mathcal{W}}

\def\Rep{\operatorname{Rep}}
\def\Irr{\operatorname{Irr}}

\newcommand{\OOO}{\mathcal{O}}

\newcommand{\Fam}{\mathcal{F}}

\newcommand{\Ind}{\operatorname{Ind}}
\newcommand{\ind}{\operatorname{ind}}

\newcommand{\Kbar}{\overline{\mathcal{K}}}

\newcommand{\Hom}{\operatorname{Hom}}
\renewcommand{\subset}{\subseteq}

\newcommand{\swrz}{\mathcal{S}}

\newcommand{\W}{\mathrm{W}}

\newcommand{\GL}{\operatorname{GL}}

\newcommand{\Ker}{\mathrm{Ker}}
\def\Spec{\operatorname{Spec}}

\newcommand{\C}{\mathbb{C}}

\newcommand{\fri}{\mathfrak{i}}
\newcommand{\M}{\mathcal{M}}
\renewcommand{\S}{\mathcal{S}}
\newcommand{\cL}{\mathcal{L}}

\newcommand{\ZZ}{\mathbb{Z}}

\newcommand{\Z}{\mathfrak{Z}}

\def\End{\operatorname{End}}

\def\GL{\operatorname{GL}}

\def\\Hom{\operatorname{\Hom}}
\def\Irr{\operatorname{Irr}}

\def\Rep{\operatorname{Rep}}
\def\Zl{\overline{\mathbb{Z}_{\ell}}}

\def\leq{\leqslant}
\def\geq{\geqslant}

\def\wflb{W(\overline{\mathbb{F}_{\ell}})}

\title{Towards a theta correspondence in families for type II dual pairs}
\author{Gil Moss and Justin Trias}
\date{}

\begin{document}

\maketitle

\begin{abstract}
Let $R$ be a commutative $\ZZ[1/p]$-algebra, let $m \leq n$ be positive integers, and let $G_n=\GL_n(F)$ and $G_m=\GL_m(F)$ where $F$ is a $p$-adic field. The Weil representation is the smooth $R[G_n\times G_m]$-module $C_c^{\infty}(\text{Mat}_{n\times m}(F),R)$ with the action induced by matrix multiplication. When $R=\C$ or is any algebraically closed field of banal characteristic compared to $G_n$ and $G_m$, the local theta correspondence holds by the work of Howe and M\'inguez. At the level of supercuspidal support, we interpret the theta correspondence as a morphism of varieties $\theta_R$, which we describe as an explicit closed immersion. For arbitrary $R$, we construct a canonical ring homomorphism $\theta^\#_{R} : \Z_{R}(G_n)\to \Z_{R}(G_m)$ that controls the action of the center $\Z_{R}(G_n)$ of the category of smooth $R[G_n]$-modules on the Weil representation. We use the rank filtration of the Weil representation to first obtain $\theta_{\ZZ[1/p]}^\#$, then obtain $\theta^\#_R$ for arbitrary $R$ by proving $\Z_R(G_n)$ is compatible with scalar extension. In particular, the map $\Spec(\Z_R(G_m))\to \Spec(\Z_R(G_n))$ induced by $\theta_R^\#$ recovers $\theta_R$ in the $R=\C$ case and in the banal case. We use gamma factors to prove $\theta_R^\#$ is surjective for any $R$. Finally, we describe $\theta^\#_R$ in terms of the moduli space of Langlands parameters and use this description to give an alternative proof of surjectivity in the tamely ramified case.
\end{abstract}

\section{Introduction}

For a non-archimedean local field with residue field of characteristic $p\neq 2$, the theta correspondence involves two groups forming a dual pair in a symplectic group, and provides a bijection between certain subsets of irreducible representations of (central extensions of) the two groups. It is an important tool in the theory of automorphic forms, one famous application being the construction of counterexamples to the generalized Ramanujan--Petersson conjecture (\cite{howe-PS}), others being cases of the local Langlands conjectures (\cite{gan_takeda_gsp4}) and Gan-Gross-Prasad conjectures (\cite{gan_ichino_2016}). Some of the deepest arithmetic properties of automorphic forms are the congruences they satisfy, which in turn come from congruences in the local representation theory. Instead of working in the traditional setting of complex representations, we consider the local theta correspondence for $\ell$-modular representations of a $p$-adic field, where $\ell\neq p$, or, more generally, representations on $R$-modules where $R$ is a $\ZZ[1/p]$-algebra. The latter constitutes a ``family'' in the sense that its fiber at each $x\in \Spec(R)$ is a traditional representation on a $\kappa(x)$-vector space, where $\kappa(x)$ is the residue field of $R$ at $x$. The tools and perspectives needed in this framework can lead to new insights even when specialized to the complex setting.

Dual pairs divide into two main kinds: type I and type II. Type I involves isometry groups such as symplectic, orthogonal and unitary, and type II involves general linear groups over skew fields. Type II dual pairs over $p$-adic fields are a natural place to begin investigating the theta correspondence in families because $\ell$-adic families with $\ell\neq p$ have been well-studied, especially in the context of describing how the deformation theory of $\ell$-adic Galois representations is reflected in the local Langlands correspondence for $\GL_n$ (\cite{emerton_helm_families,helm_curtis, HMconverse}). Recently there has been growing interest in working independently of $\ell$ by using global coefficient rings over $\mathbb{Z}[1/p]$; we take that approach here.

Let $F$ be a non-archimedean local field of residual characateristic $p$ and residual cardinal $q$. We allow $p=2$, as opposed to the type I case. Let $R$ be a commutative $\ZZ[1/p]$-algebra, let $m \leq n$ be two positive integers, and set $G_n=\GL_n(F)$ and $G_m = \GL_m(F)$. The group $G_n\times G_m$ acts by left and right translation on the Weil representation $\omega_{n,m}^R =C_c^{\infty}(\mathcal{M}_{n,m}(F),R)$, which is the space of smooth compactly supported $R$-valued functions on the set of $n$ by $m$ matrices. When $R=\mathbb{C}$ and $\pi$ is in $\Irr_{\C}(G_m)$, the type II theta correspondence says there is a finite length $G_n$-representation $\Theta(\pi)$ such that over $\Theta(\pi) \otimes \pi$ is the $(G_m,\pi)$-isotypic quotient of $\omega_{n,m}$ and $\Theta(\pi)$ has a unique irreducible quotient $\theta_\C(\pi)$. This theorem is due to Howe in unpublished work; a proof can be found in the appendix of M\'inguez' thesis \cite{minguez_thesis}. In \cite{minguez_thesis}, M\'inguez works over an algebraically closed field $R$ of characteristic $\ell\neq p$ and establishes an injective map of irreducible $R$-representations 
$$\theta_R : \Irr_R(G_m) \to \Irr_R(G_n), \ m \leq n,$$
but only in the cases where $\ell$ is banal with respect to $G_n$ and $G_m$, i.e., when $\ell$ does not divide the pro-orders of these two groups (the case $\ell=0$ is always banal by convention). M\'inguez makes the map explicit in terms of the Langlands quotient classification of irreducibles (\cite{minguez_howe}). However, when $\ell$ divides $q^n-1$, this already fails: restriction to $\{ 0 \}$ and the Haar measure $\mu_{F^n}$ give a surjection
$$\omega_{n,1} \twoheadrightarrow (1_n \otimes 1) \oplus (|\cdot|_n  \otimes 1),$$ where $1_n$ denotes the trivial representation of $G_n$ and $|\cdot|_n = |\cdot|_F\circ \det$. Thus there is no map $\theta_R$ in the traditional sense and a new perspective is needed to formulate a theta correspondence mod-$\ell$ or in families. On closer inspection, $1_n$ and $|\cdot|_n$ have the same supercuspidal support when $\ell|(q^n-1)$, which suggests considering a map on supercuspidal supports, as we now explain.

Using the explicit description in \cite{minguez_howe} in terms of the Langlands classification, one finds $\theta_R$ is indeed compatible with supercuspidal supports in the banal or complex settings. More precisely, there exists an injective map -- which we still call $\theta_R$ -- between the sets of supercuspidal supports $\Omega_R(G_m) \to \Omega_R(G_n)$ such that when $\ell$ is banal the following diagram commutes:
$$\xymatrix{
\Irr_R(G_m) \ar[r]^{\theta_R} \ar[d]^{\textup{scs}}  & \Irr_R(G_n) \ar[d]^{\textup{scs}}  \\
\Omega_R(G_m) \ar[r]^{\theta_R} & \Omega_R(G_n).
}$$
We can describe $\theta_R$ explicitly as follows. The Bernstein decomposition is a disjoint union $\Omega_R(G_m) = \bigsqcup_{\mathfrak{s}\in\mathcal{B}(G_m)}\Omega_R^{\mathfrak{s}}(G_m)$ where $\mathcal{B}(G_m)$ is the set of inertial supercuspidal supports. If $(M,\rho)_{\text{scs}}$ is a supercuspidal support consisting of the $G_m$-conjugacy class of a Levi $M$ and a cuspidal representation $\rho$ of $M$, we have (in the banal case) an injection of sets
\begin{align*}
\theta_R:\Omega_R^{\mathfrak{s}}(G_m)&\to \Omega_R^{\theta(\mathfrak{s})}(G_n)\\
(M,\rho)_{\text{scs}}&\mapsto (M\times T_{n-m}, \rho^{\vee}\chi_M\otimes_R \chi_{T_{n-m}})_{\textup{scs}}\ ,
\end{align*}
where
$$\chi_M = |\cdot|^{-\frac{n-m}{2}} \textup{ and } \chi_{T_{n-m}} = |\cdot|_1^{(m+1-n)+\frac{(n-1)}{2}} \otimes_R \cdots \otimes_R |\cdot|_1^{\frac{(n-1)}{2}}.$$

\subsection{Theta map as a morphism of varieties} Let $R$ be an algebraically closed field with characteristic $\ell$ different from $p$. The set $\Omega_R(G_m)$ enjoys a richer structure of a disjoint union of affine algebraic varieties and it is natural to ask whether $\theta_R$ preserves this geometric structure. As a result of Schur's lemma, the center $\Z_R(G_m)$ of the category $\Rep_R(G_m)$ acts on $\pi \in \Irr_R(G_m)$ by a character $\eta_\pi : \Z_R(G_m) \to R$, and 
\begin{align*}
\Irr_R(G_m) &\to \Hom_{R-\textup{alg}}(\Z_R(G_m),R)\\
\pi &\mapsto \eta_{\pi}
\end{align*}
is a surjective map whose finite fibers are precisely the partition of $\Irr_R(G_m)$ by supercuspidal supports. Via this identification, supercuspidal supports are $R$-points (in the sense of Zariski) of $X_{n,R} = \textup{Spec}(\Z_R(G_m))$. In banal characteristic, the Bernstein decomposition coincides with a decomposition
$$\Z_R(G_m) = \prod_{\mathfrak{s} \in \mathcal{B}_R(G_m)} \Z_R^\mathfrak{s}(G_m)$$ into integral domains that are finite type $R$-algebras (see Lemma~\ref{bernstein_center_components_lemma}), so we are asking for a homomorphism of rings $$\theta_R^\# : \Z_R^{\theta(\mathfrak{s})}(G_n)\to \Z_R^{\mathfrak{s}}(G_m).$$ In Section~\ref{section:theta_and_scs} we prove $\theta_R$ is a closed immersion of varieties in the banal setting, hence $\theta_R^\#$ is surjective. When $R=\C$ this question has been addressed, but only for irreducible unramified representations: in \cite{rallis}, Rallis produced a map of spherical Hecke algebras giving the type~II theta correspondence on Satake parameters.

When $\ell$ is non-banal, an explicit description of the Bernstein components is much less straightforward. Worse, $\Z_R^{\mathfrak{s}}(G_m)$ can fail to be reduced, so even if we produce a candidate for a map on points $\theta_R:\Omega_R(G_n) \to \Omega_R(G_m)$, it won't uniquely determine a morphism of schemes, and might not be the ``right one.'' However, the center $\Z_{\ZZ[1/p]}(G_n)$ of the category of smooth $\ZZ[1/p][G]$-modules is reduced. Our strategy, therefore, is to produce a canonical map $$\theta^\#_{\ZZ[1/p]}:\Z_{\ZZ[1/p]}(G_n)\to \Z_{\ZZ[1/p]}(G_m)$$ by studying the subquotients of the rank filtration of the Weil representation over $\ZZ[1/p]$. 

\begin{thm}\label{main_thm_intro}
There exists a unique homomorphism 
$$\theta^\#_{\ZZ[1/p]}:\Z_{\ZZ[1/p]}(G_n)\to \Z_{\ZZ[1/p]}(G_m),$$  such that the kernel of the natural map
$$\Z_{\ZZ[1/p]}(G_n)\otimes_{\ZZ[1/p]} \Z_{\ZZ[1/p]}(G_m)\to \End_{\ZZ[1/p][G_n\times G_m]}(\omega_{n,m}^{\ZZ[1/p]})$$ is the ideal generated by $\{z\otimes 1 - 1 \otimes\theta_{\ZZ[1/p]}^\#(z):\ z\in \Z_{R}(G_n)\}$.
\end{thm}

For any $\ZZ[1/p]$-algebra $R$, let $\Rep_R(G_n) = \prod_{r\geq 0}\Rep_R(G_n)_r$ denote the decomposition of $\Rep_R(G_n)$ according to depth in the sense of Moy--Prasad (see \cite[Appendix A]{dat_finitude}) and let $\Z_R(G_n)=\prod_r\Z_R(G_n)_r$ denote the corresponding factorization of the center. We show in Section~\ref{DEPTH SECTION} that $\theta_{\ZZ[1/p]}^\#$ preserves depth in the sense that $\theta_{\ZZ[1/p]}^\#\left(\Z_{\ZZ[1/p]}(G_n)_r\right)\subset \Z_{\ZZ[1/p]}(G_m)_r$. In particular the depth-$r$ summand of $\omega_{n,m}^{\ZZ[1/p]}$ in $\Rep_{\ZZ[1/p]}(G_n)$ is contained in $\Rep_{\ZZ[1/p]}(G_m)_r$.

In order to define $\theta_R^\#$ for an arbitrary $\ZZ[1/p]$-algebra $R$ we need some compatibility with extension of scalars. The Weil representation is easily seen to be compatible with arbitrary scalar extensions, but this compatibility for the Bernstein center is not obvious when the extension is not flat (for example, $\ZZ[1/p]\to \overline{\mathbb{F}}_{\ell}$ is not flat). Thus an essential input for our strategy is the following theorem. 
\begin{thm}\label{compatibility_scalars}
For any $\ZZ[1/p]$-algebra $R$, the natural map $\Z_{\ZZ[1/p]}(G_m)_r\otimes_{\ZZ[1/p]}R\to \Z_R(G_m)_r$ is an isomorphism.
\end{thm}
\noindent Note that Theorem~\ref{compatibility_scalars} would fail without restricting to a finite number of factors of $\Z_{\ZZ[1/p]}(G_m)$ as tensor product does not commute with infinite products. This property has been widely expected among experts. We give a proof in the appendix, following a suggestion of D. Helm. 

Using Theorem~\ref{compatibility_scalars} we define, for any $\ZZ[1/p]$-algebra $R$ and any depth $r$, a homomorphism
$$\theta^\#_{R,r}:\Z_{R}(G_n)_r\to \Z_{R}(G_m)_r,$$ and let $\theta_R^\#=\prod_r\theta_{R,r}^\#$. Because of the way it is constructed, we obtain:
\begin{thm}\label{thm_intro_theta_R}
The map $\theta_R^\#$ satisfies the following two conditions. 
\begin{enumerate}
\item The kernel of the natural map
$$\Z_{R}(G_n)\otimes_{R} \Z_{R}(G_m)\to \End_{R[G_n\times G_m]}(\omega_{n,m}^{R})$$ is the ideal generated by $\{z\otimes 1 - 1 \otimes\theta_{R}^\#(z):\ z\in \Z_{R}(G_n)\}$.
\item For any field $R$ of banal characteristic, the map between algebraic varieties induced by $\theta_R^\#$ is the algebraic morphism $\theta_R$ previously defined.
\end{enumerate}
\end{thm}

For example, if $R=\mathbb{C}$ and $\mathfrak{s}$ is the inertial class of the principal block, Rallis's map is exactly $\theta_{\C}^\#$  restricted to $\Z^{\mathfrak{s}}_{\C}(G_n)$. While our correspondence can be compared to that of M\'{i}nguez or Rallis, we do not use their results as an input in our construction.

When $n=m$, $\theta_R^\#$ is a natural duality involution $\Z_R(G_n)\to \Z_R(G_n)$.

\subsection{Finiteness, inductive relations, and surjectivity of $\theta_R^\#$}

The subquotients in the rank filtration of $\omega_{n,m}^{\ZZ[1/p]}$ can be realized as parabolic inductions from Levi subgroups. Thus $\theta_{\ZZ[1/p]}^\#$ factors through a so-called Harish-Chandra homomorphism $\text{HC}: \Z_{\ZZ[1/p]}(G_n)\to \Z_{\ZZ[1/p]}(M^n_m)$, where $M^n_m$ is the Levi subgroup $\left(\begin{smallmatrix} G_m & \\ & G_{n-m}\end{smallmatrix}\right)$, where $\text{HC}$ is defined by the property that the action of $\Z_{\ZZ[1/p]}(G_n)$ on objects parabolically induced from $M^n_m$ factors through $\text{HC}$. 

It has recently been established that Harish-Chandra morphisms over $\ZZ[1/p]$ are finite \cite[Th 4.1]{dhkm_finiteness}, from which we can deduce:
\begin{thm}\label{thm:finiteness_intro}
The homomorphisms $\theta_R^\#$ of Theorem~\ref{main_thm_intro} are finite.
\end{thm}

By realizing $\theta_{\ZZ[1/p]}^\#$ in terms of a Harish-Chandra morphism and using induction in stages, we deduce an interesting recurrence relation. Denoting by $\theta^\#_{R,n,m}$ the map $\theta^\#_R : \Z_{R}(G_n)\to \Z_{R}(G_m)$, we establish $$\theta_{R,n,m}^\# = \theta_{R,k,m}^\# \circ \theta_{R,k,k}^\# \circ \theta_{R,n,k}^\#$$ for all $m \leq k \leq n$. 

When $R$ is an algebraically closed field of banal characteristic, our explicit description of $\theta_R$ in terms of supercuspidal support allows us to deduce $\theta_R$ is a closed immersion, i.e., $\theta_R^\#$ is surjective. Actually, this phenomenon has a bigger scope. In Section~\ref{sec:surjectivity} we prove:
\begin{thm}\label{thm:surjectivity_intro}
For any $\ZZ[1/p]$-algebra $R$, the homomorphism $\theta_R^\#$ is surjective. 
\end{thm}
Extending surjectivity beyond the banal setting requires significant new techniques; we use Rankin-Selberg gamma factors and converse theorems in families, which have only recently been developed in \cite{matringe_moss, Moss16.1, HMconverse}. By the recurrence relation above, we reduce to the case $m=n-1$. Then, generalizing the ideas of Watanabe in \cite{watanabe}, we realize $\theta_R$ in terms of the action of the Bernstein center on the induced module $\ind_{P^n_{n-1}}^{G_n}(\Gamma_n\otimes 1)$, where $\Gamma_n$ is the Gelfand--Graev representation for a fixed additive character $\psi$ on $F$. We then establish a multiplicativity property for gamma factors of induced modules, and apply a ``gamma factor descent'' technique used in \cite{HMconverse}. Note that while Theorem~\ref{thm:surjectivity_intro} implies Theorem~\ref{thm:finiteness_intro}, we use Theorem~\ref{thm:finiteness_intro} as an input for proving Theorem~\ref{thm:surjectivity_intro}. 

In Section~\ref{sec:galois} we give a second proof of surjectivity in depth zero by using the local Langlands correspondence in families and interpreting $\theta_R$ in terms of the geometry of the space of semisimple Langlands parameters.

We can now formulate a modular theta correspondence on the level of supercuspidal supports for any any algebraically closed field $R$ of characteristic $\ell\neq p$. More precisely, let $\eta:\Z_R(G_m)\to R$ be a character 
(\textit{i.e.}, a supercuspidal support), let $(\omega_{n,m})_{\eta}$ be the largest quotient of $\omega_{n,m}$ on which $\Z_R(G_m)$ acts via $\eta$, and let $\theta_R(\eta) = \eta\circ \theta_R^\#$. In Section~\ref{sec:applications} we use the surjectivity and finiteness statements above to deduce the following:
\begin{corollary}\label{cor:modular_theta}
The largest quotient of $\omega_{n,m}$ on which $\Z_R(G_n)$ acts by $\theta_R(\eta)$ is $(\omega_{n,m})_{\eta}$.
\end{corollary}
Given $\pi\in \Irr_R(G_m)$, let $\pi\otimes \Theta(\pi)$ be the canonical decomposition of the $(G_m,\pi)$-isotypic quotient of $\omega_{n,m}$. In Section~\ref{sec:applications}, we prove that $\Theta(\pi)$ has finite length. It follows from Corollary~\ref{cor:modular_theta} that $\Z_R(G_n)$ acts on every irreducible constituent of $\Theta(\pi)$ via the character $\theta_R(\eta_{\pi})$; in particular all the constituents of $\Theta(\pi)$ have the same supercuspidal support.

\subsection{Functoriality}
\label{sec:functoriality_intro}

Let $\widehat{G}_n$ denote the algebraic group $\GL_n$ defined over $\ZZ[1/p]$. Let $W_F$ be the Weil group of $F$ and let $\Phi^{WD}_\mathbb{C}(n)$ denote the set of $\widehat{G}_n(\mathbb{C})$-conjugacy classes of Frobenius-semisimple complex Weil-Deligne representations $(\rho,N)$ where $\rho:W_F\to \widehat{G}_n(\C)$ is a smooth homomorphism and $N$ is a nilpotent operator satisfying the usual relation. The classical local Langlands correspondence for $G_n$ is a canonical bijection 
$$\cL_{n,\C}:\Irr_{\mathbb{C}}(G_n) \to \Phi^{WD}_{\mathbb{C}}(n).$$
Let $\nu$ denote the character of $W_F$ corresponding to $|\cdot|_F$ under local class field theory and define $$\hat{\phi} := \phi^{\vee}\cdot\nu^{-\frac{n-m}{2}}\oplus \nu^{m+1-n+\frac{n-1}{2}}\oplus \cdots \oplus \nu^{\frac{(n-1)}{2}}.$$ In \cite{minguez_howe}, M\'inguez showed the following diagram commutes
$$\begin{tikzcd}
\Irr_{\C}(G_n)\arrow[r,"\cL_{n,\C}"] & \Phi^{WD}_{\C}(n)\\
\Irr_{\C}(G_m) \arrow[u,"\theta_\C"]\arrow[r,"\cL_{m,\C}"] & \Phi^{WD}_{\C}(m)\arrow[u,"\phi \mapsto \hat{\phi}"']\ .
\end{tikzcd}
$$
It is natural to ask whether our integral morphism $\theta_{\ZZ[1/p]}$ is reflected on the $W_F$ side of the local Langlands correspondence. As we will now explain, the proper context for such a question is the local Langlands correspondence \emph{in families}, which interpolates $\cL_{n,\C}$ to a morphism on the integral Bernstein variety.

The semisimplified local Langlands correspondence $\cL_{n,\C}^{ss}$ is the composition of $\cL_{n,\C}$ with
\begin{align*}
\Phi^{WD}_{\C}(n)&\to \Phi_{\C}(n)\\
(\rho,N)&\mapsto \rho\ ,
\end{align*}
where $\Phi_{\C}(n)$ denotes the conjugacy classes of smooth homomorphisms $\rho:W_F\to \widehat{G}_n(\C)$ whose images consist of semisimple elements. Since $\cL_{n,\C}^{ss}$ is constant on the fibers of the supercuspidal support map $\Irr_{\mathbb{C}}(n)\xrightarrow{\text{scs}} \Omega_{\mathbb{C}}(G_n)$, it defines a map 
$$\Omega_{\C}(G_n)\to \Phi_{\C}(n),$$ which we will also denote $\cL_{n,\C}^{ss}$. 

As described above, $\Omega_\C(G_n)$ possesses more structure: it is the $\C$-points of the integral Bernstein variety $\Spec(\Z_{\ZZ[1/p]}(G_n))$, a disjoint union of finite type affine schemes over $\ZZ[1/p]$. In fact, $\Phi_{\C}(n)$ also has a geometric structure and $\cL_{n,\C}^{ss}$ can be upgraded to a morphism of $\mathbb{Z}[\sqrt{q}^{-1}]$-schemes. More precisely, let $(P_F^e)_{e\in \mathbb{N}}$ be a descending filtration of the wild inertia subgroup $P_F$ by open subgroups that are normal in $W_F$ with $P_F^0=P_F$. Let $s$ be a topological generator of the tame quotient $I_F/P_F$, let $\text{Fr}$ be a lift of Frobenius in $W_F/P_F$ and let $W_F^0$ be the preimage in $W_F$ of the discrete subgroup $\langle \text{Fr},s\rangle\subset W_F/P_F$. We define $X_n^e$ as the scheme over $\ZZ[1/p]$ representing homomorphisms from the finitely presented group $W_F^0/P_F^e$ to $\widehat{G}_n$.  The subset $\Phi^e_{\C}(n)$ consisting of semisimple parameters trivial on $P_F^e$ is equivalent to the $\C$-points of the affine GIT quotient scheme $X_n^e\sslash \widehat{G}_n = \Spec(\OOO[X_n^e]^{\widehat{G}_n})$. This coarse moduli space was introduced by Helm in \cite{helm_curtis} over $\wflb$ and expanded to $\ZZ[1/p]$ and other groups beyond $\GL_n$ in \cite{dhkm_moduli,fargues_scholze, zhu}. The local Langlands correspondence in families states that for each $e$ there exists a direct factor $\Z_{\ZZ[\sqrt{q}^{-1}]}^e(G_n)$ of $\Z_{\ZZ[\sqrt{q}^{-1}]}(G_n)$ and an isomorphism of rings:
$$\cL_{n,e}^{ss,\#}:\OOO[X_n^e]^{\widehat{G}_n}\to \Z^e_{\ZZ[\sqrt{q}^{-1}]}(G_n),$$ that recovers $\cL_n^{ss}$ on $\Phi^e_{\C}(n)$ upon extending scalars to $\C$. It was established in \cite{helm_curtis, HMconverse} over $\wflb$ (taking $\cL_{n,\C}$ as an input) and the morphism descends to $\ZZ[\sqrt{q}^{-1}]$ with formal methods (forthcoming work \cite{dhkm_banal}).

Consider the map on GIT quotients ${^L\theta_e}:X_m^e\sslash \widehat{G}_m \to X_n^e\sslash \widehat{G}_n$ that makes the following square commute
$$\begin{tikzcd}
\OOO[X^e_n]_{\ZZ[\sqrt{q}^{-1}]}^{\widehat{G}_n}\arrow[d,"{^L\theta}_e^\#"']\arrow[r,"\cL^{ss,\#}_{n,e}","\sim"'] & \Z^e_{\ZZ[\sqrt{q}^{-1}]}(G_n)\arrow[d,"\theta_{\ZZ[\sqrt{q}^{-1}],e}^\#"]\\
\OOO[X^e_m]_{\ZZ[\sqrt{q}^{-1}]}^{\widehat{G}_m}\arrow[r,"\cL^{ss,\#}_{m,e}","\sim"'] & \Z^e_{\ZZ[\sqrt{q}^{-1}]}(G_m)
\end{tikzcd}\ .
$$
We give a straightforward description of ${^L\theta}_e$ in Proposition~\ref{prop:description_of_Ltheta}. For simplicity we will restrict ourselves in this subsection to the depth zero case where $e=0$, i.e., where parameters are tamely ramified (trivial on $P_F$). In this case, the tame parameter space $X_m^0$ is the closed affine subscheme of $\widehat{G}_m \times \widehat{G}_m$ representing pairs $(\Fam, \sigma)$ such that $\Fam\sigma \Fam^{-1} = \sigma^q$. The map sending $(\Fam,\sigma)\in X_m^0$ to the pair
\begin{align*}
\left(
\begin{pmatrix}
q^{-\frac{n-m}{2}}I_m\cdot {^t\Fam^{-1}} \\ 
& q^{m+1 - n + \frac{n-1}{2}}\\
&&\ddots\\
&&&q^{\frac{(n-1)}{2}}
\end{pmatrix},
\begin{pmatrix}
{^t\sigma^{-1}}\\
&1\\
&&\ddots\\
&&&1
\end{pmatrix}\right)\in X_n^0
\end{align*}
is equivariant for the conjugation action on source and target, and ${^L\theta}_0$ is the induced map on GIT quotients $X_m^0\sslash \widehat{G}_m\to X_n^0\sslash \widehat{G}_n$. 

Thus $\theta_{\ZZ[\sqrt{q}^{-1}]}$ is \emph{almost} an instance of Langlands' functoriality principle in that, up to accounting for some unramified twisting, ${^L\theta_e}$ is given by a homomorphism of $L$-groups \begin{align*}
\widehat{G}_m&\to \widehat{G}_n\\
g&\mapsto \begin{pmatrix}{^tg^{-1}} & 0 \\ 0 & I_{n-m} \end{pmatrix}\ . 
\end{align*}

\noindent By invoking the local Langlands isomorphism in families and our Theorem~\ref{thm:surjectivity_intro}, we deduce
\begin{corollary}\label{cor:surjectivity_galois_intro}
The map ${^L\theta_e}$ is a closed immersion, i.e., ${^L\theta_e^\#}$ is surjective.
\end{corollary}

Given the description of ${^L\theta_e}$, it is natural to seek a direct proof that ${^L\theta_e}$ is a closed immersion, one which is geometric and does not rely on gamma factors or the local Langlands correspondence in families. This turns out to be deep because the GIT quotient is a subtle construction in the integral setting. In Section~\ref{sec:galois}, we give a proof in the depth zero case relying heavily on \cite[Th VIII.0.2]{fargues_scholze}, which is a difficult result relating $\OOO[X_n^e]^{\widehat{G}_n}$ to an algebra of excursion operators. Our approach could probably be extended to positive depth using the reduction-to-depth-zero machinery in \cite{dhkm_moduli}, but we have not pursued this.



\subsection{Further directions} The questions we have investigated in this paper can be posed in other contexts. For example, some techniques in this paper might generalize to the setting of type II dual pairs over division algebras. While the geometric lemma and the Harish-Chandra morphisms generalize, the local Langlands correspondence in families and the theory of gamma factors would require new developments. As another example, the compatibility of the theta correspondence with supercuspidal support for dual pairs $(G,G')$ of type I has been established in \cite{kudla_86} and one can ask about the algebraicity and integrality of such a map. Even though the algebraicity seems rather straightforward from Kudla's formulas, integrality seems to be difficult. The authors plan to investigate these questions in future work.

\subsection{Acknowledgements} The authors are grateful for helpful comments and suggestions from Jean-Fran\c{c}ois Dat, David Helm, Rob Kurinczuk, Alberto M\'{i}nguez, Vincent S\'echerre and Shaun Stevens. The first author was partially supported by NSF grants DMS-2001272 and DMS-2302591. The second author would also like to thank Jack Sempliner for useful conversations on algebraic geometry.

\section{The filtration of the Weil representation and its endomorphisms}

In this section, our coefficient ring $R$ is any commutative $\mathbb{Z}[1 / p]$-algebra. In particular we do not assume the existence of a square root of $p$ in $R$. We denote by $1$ the trivial representation, which is the free $R$-module of rank one $R$ with trivial group action. Recall that the linear action of $G_n \times G_m$ on the Weil representation is given by matrix multiplication:
$$(\omega_{n,m}^R(g_n,g_m) \cdot f ) (x) = f(g_n^{-1} x g_m),$$
for $x\in \mathcal{M}_{n,m}(F)$ and $f\in C_c^{\infty}(\mathcal{M}_{n,m}(F),R)$.
As the context here should be clear, we drop the reference to $R$ in  spaces of functions and in the Weil representation to lighten notations.

We will study the rank filtration of the Weil representation and show that the endomorphism rings of its subquotients are identified to the center of the category of smaller and smaller general linear groups. The main tool we use is the so-called geometric lemma. To preserve the flow of the exposition, we have relegated the details of the geometric lemma to Appendix \ref{geometric_lemma_appendix}. We briefly recall some notations from the appendix. Choose as a minimal parabolic subgroup of $G_n$, also called a Borel subgroup in this situation, the subgroup of upper triangular matrices $B_n$ with Levi decomposition $T_n N_n$ where $T_n$ is the subgroup of diagonal matices in $G_n$ and $N_n$ the subgroup of unipotent matrices in $B_n$. For $0 \leq k \leq n$, set:
$$M_k^n = \left\{ \left. \left[ \begin{array}{cc}
a_k & 0 \\
0 & b_{n-k}
\end{array} \right] \in G_n \ \right\vert \ a_k \in G_k \textup{ and } b_{n-k} \in G_{n-k} \right\}.$$
It is a standard Levi of $G_n$, which is contained in a unique standard parabolic subgroup denoted by $P_k^n = M_k^n N_k^n$. Let $Q_k^n = M_k^n \bar{N}_k^n$ be the opposite parabolic to $P_k^n$ with respect to $B_n$. We use similar notations for $G_m$. From now on assume that $n \geq m$.

\subsection{Filtration by the rank}\label{subsection:filtration_by_rank} Let $\mathcal{O}_k$ be the set of rank $k$ matrices in $\mathcal{M}_{n,m}(F)$ and write:
$$\mathcal{M}_{n,m}(F) = \coprod_{0 \leq k \leq m} \mathcal{O}_k\ .$$
Each $\mathcal{O}_k$ is a single $(G_n \times G_m)$-orbit that is also a locally closed subset of $\mathcal{M}_{n,m}(F)$. Denote by $U_k = \coprod_{l \geq k} \mathcal{O}_l$ the set of matrices of rank at least $k$. The set $U_k$ is a $(G_n \times G_m)$-stable open subset of $\mathcal{M}_{n,m}(F)$ and $\mathcal{O}_k$ is closed in $U_k$, yielding a stratification of the space $\mathcal{M}_{n,m}(F)$. Take representatives for $(\mathcal{O}_k)_{0 \leq k \leq m}$ by setting:
$$x_k = \left[ \begin{array}{cc}
\textup{Id}_k & 0 \\
0 & 0 \end{array} \right] \in \mathcal{M}_{n,m}(F).$$
Denote by $\textup{St}_k$ the stabiliser of $x_k$, which is the normal subgroup of $P_k^n \times Q_k^m$ defined by:
$$\textup{St}_k = \left. \left\{ \left(
\left[ \begin{array}{cc}
\alpha & * \\
0 & * \end{array} \right] ,
\left[ \begin{array}{cc}
\alpha & 0 \\
* & *
\end{array} \right] \right) \in G_n \times G_m \ \right\vert \ \alpha \in G_k \right\}.$$
Write $C_c^\infty(G_k) \otimes_R 1$ for the representation of $P_k^n \times Q_k^m$, where the $G_k$ factor of $P_k^n$ acts by left multiplication on $C_c^{\infty}(G_k)$ and that of $Q_k^m$ by right multiplication, \textit{i.e.}, where the action of $P_k^n \times Q_k^m$ is given, for $f \in C_c^\infty(G_k)$, by:
$$\left(
\left[ \begin{array}{cc}
\alpha & * \\
0 & * \end{array} \right] ,
\left[ \begin{array}{cc}
\alpha' & 0 \\
* & *
\end{array} \right] \right) \cdot f : x \in G_k \mapsto f(\alpha^{-1} x \alpha') \in R.$$

\begin{prop} \label{filtration_weil_rep_proposition} Set $\omega_{n,m}^{(k)} = C_c^\infty(U_k)$.
\begin{enumerate}[label=\textup{\alph*)}]
\item \label{filtration_weil_rep_proposition_pt_1} The rank induces a filtration in $\textup{Rep}_R(G_n \times G_m)$ of the Weil representation:
$$0 \subset \omega_{n,m}^{(k)} \subset \dots \subset \omega_{n,m}^{(1)} \subset \omega_{n,m}^{(0)} = \omega_{n,m},$$
where each subquotient is canonically isomorphic to some $C_c^\infty(\mathcal{O}_k)$.
\item \label{filtration_weil_rep_proposition_pt_2} Define:
$$W_{n,m}^k = \ind_{P_k^n \times Q_k^m}^{G_n \times G_m} ( C_c^\infty(G_k) \otimes_R 1)\ .$$
Then the orbit map $g \in ( G_n \times G_m ) \mapsto g^{-1} \cdot x_k \in \mathcal{O}_k$ factors through an homeomorphism $\textup{St}_k \backslash ( G_n \times G_m ) \simeq \mathcal{O}_k$, which induces canonical isomorphisms:
$$C_c^\infty(\mathcal{O}_k) \simeq C_c^\infty(\textup{St}_k \backslash (G_n \times G_m)) \simeq W_{n,m}^k.$$ \end{enumerate} \end{prop}

\begin{proof}\ref{filtration_weil_rep_proposition_pt_1}
Because $\mathcal{O}_k$ is closed in $U_k$, we have for all $k$ that:
$$0 \to C_c^\infty(U_{k+1}) \overset{i_{k+1}}{\to} C_c^\infty(U_k) \overset{p_k}{\to} C_c^\infty(\mathcal{O}_k) \to 0$$
where $i_{k+1}$ is the obvious inclusion from $U_{k+1} \subset U_k$ and $p_k$ is the support restriction to $\mathcal{O}_k$. Collecting these many exact sequences for $0 \leq k \leq m$ yields the filtration.

\noindent \ref{filtration_weil_rep_proposition_pt_2} First, by \cite[II.3.3 Cor]{ren}, the map $g \in G_n \times G_m \mapsto g^{-1} \cdot x \in \mathcal{O}_k$ induces a homeomorphism $\textup{St}_k \backslash (G_n \times G_m) \simeq \mathcal{O}_k$. So $C_c^\infty(\mathcal{O}_k) \simeq \textup{ind}_{\textup{St}_k}^{G_n \times G_m} (1)$ where $1$ is the trivial representation. As $\textup{St}_k$ is a normal subgroup of $P_k^n \times Q_k^m$, we get that:
$$\textup{ind}_{\textup{St}_k}^{P_k^n \times Q_k^m} (1) \simeq C_c^\infty(G_k) \otimes_R 1.$$
Furthermore the action of $N_k^m \times \bar{N}_k^m$ is trivial on $C_c^\infty(G_k) \otimes_R 1$ because it is contained in $\textup{St}_k$. So by transitivity of induction $\textup{ind}_{\textup{St}_k}^{G_n \times G_m} (1) \simeq \ind_{P^k_n \times Q_k^m}^{G_n \times G_m}(C_c^\infty(G_k) \otimes_R 1)$ in $\textup{Rep}_R(G_n \times G_m)$. \end{proof}

\subsection{Example when $n=2$ and $m=1$} In this situation we have a filtration:
$$0 \subset \omega_{2,2}^{(1)} \subset \omega_{2,2}^{(0)} = \omega_{2,2}.$$
Denoting by $B_2$ the standard Borel subgroup of $G_2$ of upper triangular matrices, we have: 
$$W_{2,2}^1 = \textup{ind}_{B_2}^{G_2} (C_c^\infty(G_1) \otimes_R 1) \textup{ and } W_{2,2}^0 = 1.$$
Note that $\End_{G_2 \times G_2}(1) \simeq \Z_R(G_0) = R$. In addition $\Hom_{G_2 \times G_2}(W_{2,2}^0,W_{2,2}^1)=0$ as, for compact support reasons, there is no function in $ \textup{ind}_{B_2}^{G_2} (C_c^\infty(G_1) \otimes_R 1)$ with support $G_2$. Studying the endomorphism ring of the remaining subquotient, Frobenius reciprocity reads:
$$\End_{G_2 \times G_1}(W_{2,2}^1) \simeq \Hom_{T_2 \times G_1} (  \mathfrak{r}_{G_2}^{T_2} (W_{2,2}^1) , C_c^\infty(G_1) \otimes_R 1)$$
and as a consequence of the geometric lemma:
$$0 \to \delta_{B_2} \cdot (1 \otimes_R C_c^\infty(G_1)) \to \mathfrak{r}_{G_2}^{T_2} (W_{2,2}^1) \to C_c^\infty(G_2) \otimes_R 1 \to 0.$$
Any morphism deduced from Frobenius reciprocity must restrict to zero on $\delta_{B_2} \cdot (1 \otimes_R C_c^\infty(G_1))$. Indeed, this is a consequence of $\Hom_{T_2 \times G_1} (\delta_{B_2} \cdot (1 \otimes_R C_c^\infty(G_1)) , C_c^\infty(G_1) \otimes_R 1) = 0$ because, for compact support reasons again, we must have $\Hom_{G_1}(\chi , C_c^\infty(G_1)) = 0$. Therefore all these maps factor through the quotient in the exact sequence:
\begin{eqnarray*} \Hom_{T_2 \times G_1} (  \mathfrak{r}_{G_2}^{T_2} (W_{2,2}^1) , C_c^\infty(G_1) \otimes_R 1) & \simeq & \Hom_{T_2 \times G_1} (  C_c^\infty(G_1) \otimes_R 1 , C_c^\infty(G_1) \otimes_R 1) \\
 & \simeq & \End_{G_1 \times G_1}(C_c^\infty(G_1)) \simeq \mathfrak{z}_R(G_1).\end{eqnarray*}
Actually these rather simple ideas (nullity of some homomorphism space for compact support reasons and using the geometric lemma) transfer well to the general case, at the cost of introducing less digestible notation.

\subsection{Representations $W_{n,m}^k$} In order to study the properties of the Weil representation $\omega_{n,m}$, one can start considering the subquotients $W_{n,m}^k$ for $0 \leq k \leq m$. It happens that the endomorphism ring of $W_{n,m}^k$ is isomorphic to the Bernstein center of $G_k$, as already noted for $n=2$ and $m=1$ in the previous paragraph.

\begin{prop} \label{endomorphisms_of_W_k_nm_prop} By setting $G = G_n \times G_m$, $H= P_k^n \times Q_k^m$ and $V_k = C_c^\infty(G_k) \otimes_R 1$, the induced representation $W_{n,m}^k$ satisfies the hypothesis $\textup{Hom}_H(\textup{ker}(ev_1),V_k) = 0$ of Corollary \ref{all_morphism_are_coming_from_H_cor1}. In particular, we obtain an isomorphism of $R$-algebras:
$$\textup{End}_{G_n \times G_m}(W_{n,m}^k) \overset{\curvearrowleft}{\simeq} \textup{End}_{G_k \times G_k}(C_c^\infty(G_k)).$$ \end{prop}

\begin{proof} The idea of the proof has already been presented in the case $n=2$ and $m=1$. The proof below can be more easily navigated keeping this example in mind, as the core idea remains the same for general $n$ and $m$.

We would like to apply Corollary \ref{all_morphism_are_coming_from_H_cor1}, so we need to show that $\textup{Hom}_H(\textup{ker}(ev_1),V_k) = 0$. First of all, because the action of the radical unipotent of $H$ is trivial on $V_k$, we deduce that:
$$\textup{Hom}_H(\textup{ker}(ev_1),V_k) \simeq \textup{Hom}_{M_k^n \times M_k^m}(\mathfrak{r}_{G_n \times G_m}^{P_k^n \times Q_k^m} (\textup{ker}(ev_1)),V_k).$$
Note that $\textup{ker}(ev_1) \subset \mathfrak{i}_{P_k^n \times Q_k^m}^{G_n \times G_m} (V_k)$ is the subset of functions on $G_n \times G_m$ supported on the complement of $P_k^n \times Q_k^m$. The geometric lemma as stated in Appendix \ref{geometric_lemma_appendix} gives a filtration of $\mathfrak{r}_{G_n \times G_m}^{P_k^n \times Q_k^m}(\textup{ker}(ev_1))$ in $\textup{Rep}_R(M_k^n \times M_k^m)$. Its subquotients are: 
$$I_{w,w'} \simeq \mathfrak{i}_{P_{(k-i,i,i)}^n \times Q_{(k-j,j,j)}^m}^{M_k^n \times M_k^m} \left( \delta_{w_{k,i}^n} \delta_{w_{k,j}^m} \otimes_R \big( (w_{k,i}^n,w_{k,j}^m) \circ \mathfrak{r}_{M_k^n \times M_k^m}^{P_{(k-i,i,i)}^n \times Q_{(k-j,j,j)}^m} ( V_k ) \big) \right)$$
for $(w,w') = (w_{k,i}^n,w_{k,j}^m) \neq (\textup{Id}_n,\textup{Id}_m)$. In order to prove the condition $\textup{Hom}_H(\textup{ker}(ev_1),V_k)=0$, it is sufficient to show that $\textup{Hom}_{M_k^n \times M_k^m}(I_{w,w'},V_k) = 0$ for all $(w,w') \neq (\textup{Id}_n,\textup{Id}_m)$.

Suppose that $w \neq \textup{Id}_n$. Second adjunction is valid in this context \cite[Cor 1.3]{dhkm_finiteness}, so the $R$-module $\textup{Hom}_{M_k^n \times M_k^m}(I_{w,w'},V_k)$ is isomorphic to:
$$\textup{Hom}_{M_{(k-i,i,i)}^n \times M_k^m} \left( \mathfrak{i}_{Q_{(k-j,j,j)}^m}^{M_k^m} \left( \delta_{(w_{k,i}^n,w_{k,j}^m)} \otimes_R \big( (w_{k,i}^n,w_{k,j}^m) \circ \cdots \big) \right) , \bar{\mathfrak{r}}_{M_k^n}^{P_{(k-i,i,i)}^n} (V_k)\right).$$
Because $w \neq \textup{Id}_n$, we have $i \neq 0$. Consider the following non-trivial torus:
$$T_i^n = \left[ \begin{array}{cc} T_i & \\ & \textup{Id}_{n-k} \end{array} \right], \textup{ where } T_i = \left\{ \left. \left[ \begin{array}{cccc}
\textup{Id}_{k-i} &  \\
 & \lambda \textup{Id}_i \end{array} \right] \in G_k \ \right\vert \ \lambda \in F^\times \right\}.$$
Then we claim that $T_i^n$ acts as a character on the left-hand side of the last Hom-space above whereas it can not act as a character on the right-hand side. Indeed we have:
$$\bar{\mathfrak{r}}_{M_k^n}^{P_{(k-i,i,i)}^n} (V_k) = \mathfrak{r}_{M_k^n}^{Q_{(k-i,i,i)}^n} (V_k) \simeq C_c^\infty(\bar{N}_i^k \backslash G_k) \otimes_R 1.$$
If $T_i^n$ acts as a character on $v \in C_c^\infty(\bar{N}_i^k \backslash G_k) \otimes_R 1$, then $T_i$ acts a character on some element $v' \in C_c^\infty(\bar{N}_i^k \backslash G_k)$. But if $v' \neq 0$ then its support must contain $\bar{N}_i^k T_i$, which is not compact in $\bar{N}_i^k \backslash G_k$, so $v'$ must be zero \textit{i.e.} $v=0$. This means that the Hom-space above must be zero too as $T_i^n \subset M_{(k-i,i,i)}^n$.

Therefore $\textup{Hom}_{M_k^n \times M_k^m}(I_{w,w'},V_k) = 0$ for all $(w,w')$ with $w \neq \textup{Id}_n$. Alternatively we can conclude this is also zero for all $(w,w') \neq (\textup{Id}_n,\textup{Id}_m)$ when $w' \neq \textup{Id}_m$, just switching the roles of $w$ and $w'$ in the proof above. Therefore Corollary \ref{all_morphism_are_coming_from_H_cor1} applies as $\textup{ker}(\textup{ev}_1)$ has a filtration whose subquotients are $(I_{w,w'})$ for $(w,w') \neq (\textup{Id}_n,\textup{Id}_m)$. \end{proof}

With similar arguments to the preceding proof, we can prove:

\begin{prop} \label{hom_space_Wk_Wk'_prop} For all $k' > k$, we have $\textup{Hom}_{G_n \times G_m} ( W_{n,m}^k , W_{n,m}^{k'} ) = 0$. \end{prop}

\begin{proof} As this proof is just a variation of the previous one, we go through the arguments in a more direct way. Set $V_l = C_c^\infty(G_l) \otimes_R 1 \in \textup{Rep}_R(P_l^n \times Q_l^m)$. By Frobenius reciprocity:
$$\textup{Hom}_{G_n \times G_m} ( W_{n,m}^k , W_{n,m}^{k'} ) \simeq \textup{Hom}_{M_{k'}^n \times G_m} ( \mathfrak{r}_{G_n}^{P_{k'}^n} \circ \mathfrak{i}_{P_k^n}^{G_n} ( \mathfrak{i}_{Q_k^m}^{G_m} ( V_k ) ) , \mathfrak{i}_{Q_{k'}^n}^{G_m} ( V_{k'} )).$$
Here the version of the geometric lemma we use is again explained in Appendix \ref{geometric_lemma_appendix}, where the index set is $W(k,k',n)$ and its elements are $w_i = w_{k,k',i}^n \in W(k,k',n)$ for $i \in [ \! [ 0 , \textup{min}(k,n-k') ] \! ]$. To ease the notation set $V_l^m = \mathfrak{i}_{Q_l^m}^{G_m} ( V_l ) \in \textup{Rep}_R(P_l^n \times G_m)$. The subquotients read:
$$I_{w_i} = \mathfrak{i}_{M_{(k-i,k'-k+i,i)}^n}^{M_{k'}^n} \left( \delta_{w_i} \otimes_R (w_i \circ \mathfrak{r}_{M_k^n}^{M_{(k-i,i,k'-k+i)}^n} ( V_k^m ) \right)$$
and we want to prove that $\textup{Hom}_{M_{k'}^n \times G_m}(I_{w_i} , V_{k'}^m) = 0$ for all $w_i \in W(k,k',n)$.

Applying the second adjunction \cite[Cor 1.3]{dhkm_finiteness}, we want to prove that:
$$\textup{Hom}_{M_{(k-i,k'-k+i,i)}^n \times G_m} (\delta_{w_i} \otimes_R ( w_i \circ \mathfrak{r}_{M_k^n}^{M_{(k-i,i,k'-k+i)}^n}(V_k^m)),\bar{\mathfrak{r}}_{M_{k'}^n}^{M_{(k-i,k'-k+i,i)}^n} (V_{k'}^m) ) = 0.$$
Similarly to the previous proof, the non-trivial torus:
$$\left\{ \left.  \left[ \begin{array}{ccc} \textup{Id}_{k} & & \\ & \lambda \textup{Id}_{k'-k} & \\ & & \textup{Id}_{n-k'} \end{array} \right] \in G_n \ \right\vert \ \lambda \in F^\times \right\}$$
acts as a character on the left-hand side of the Hom-space, but it cannot act as a character on the right-hand side $\bar{\mathfrak{r}}_{M_{k'}^n}^{M_{(k-i,k'-k+i,i)}^n} (V_{k'}^m) \simeq C_c^\infty( \bar{N}_{k-i}^{k'} \backslash G_{k'}) \otimes_R 1$ for compact support reasons. Therefore $\textup{Hom}_{M_{k'}^n \times G_m}(I_{w_i} , V_{k'}^m) = 0$ for all $w_i$. \end{proof}

\section{Action of the Bernstein center on the Weil representation}

In this section we define a ring morphism $\theta_{\ZZ[1/p]}^\#$ coming from the compatibility of the action of $\Z_{\ZZ[1/p]}(G_n)\otimes \Z_{\ZZ[1/p]}(G_m)$ with the subquotients $W^k_{n,m}$ of the filtration. Then we show that the actions of $G_n$ and $G_m$ on $W^k_{n,m}$ are compatible with respect to depth. This allows us to define $\theta_R^\#$ for an arbitrary $\ZZ[1/p]$-algebra $R$ using the results of Appendix \ref{extension_of_scalars_appendix_section}.

\subsection{Action of the center on the filtration} As a result of Proposition \ref{endomorphisms_of_W_k_nm_prop}, we can consider the natural ring morphisms: 
$$\varphi_{n,m}^k : \Z_R(G_n) \otimes_R \Z_R(G_m) \to \textup{End}_{R[G_n\times G_m]}(W_{n,m}^k)$$
whose images lie naturally in $\Z_R(G_k)$. Also we define the natural action:
$$\varphi_{n,m} : \Z_R(G_n) \otimes_R \Z_R(G_m) \to \textup{End}_{R[G_n\times G_m]}(\omega_{n,m}^R).$$
By abuse of notation, we represent $\varphi_{n,m}(z)$ for $z \in \Z_R(G_n) \otimes \Z_R(G_m)$ as a matrix with respect to the rank filtration of the Weil representation. The center preserves subrepresentations, so we must have:
$$\varphi_{n,m}(z) =
\left[ \begin{array}{ccc}
\varphi_{n,m}^m(z) & * & *  \\
0 & \ddots & * \\
0 & 0 & \varphi_{n,m}^0(z) \end{array} \right].$$
In addition, the $*$ maps above are $0$ because of Proposition \ref{hom_space_Wk_Wk'_prop} and this implies the following relation on kernels:
$$\Ker(\varphi_{n,m}) = \bigcap_{k=0}^m \Ker(\varphi_{n,m}^k).$$
We will prove these kernels are well-ordered for inclusion.

\subsection{Action of the center: the case of $\ZZ[1/p]$}

\begin{prop}
\label{prop:inclusion_of_kernels} Let $R = \mathbb{Z}[1/p]$. We have inclusions of kernels:
$$\Ker(\varphi_{n,m}^m) \subset \dots \subset \Ker(\varphi_{n,m}^1) \subset \Ker(\varphi_{n,m}^0).$$ \end{prop}

\begin{proof} Let $z \in \Z_{\mathbb{Z}[1/p]}(G_n) \otimes   \Z_{\mathbb{Z}[1/p]}(G_m)$ such that $\varphi_{n,m}^{k+1}(z) = 0$, then the goal is $\varphi_{n,m}^k(z) = 0$. Because $\varphi_{n,m}^k(z)$ belongs to $\Z_{\mathbb{Z}[1/p]}(G_k)$, it is enough to check it over the complex numbers via the canonical embedding: 
$$\Z_{\mathbb{Z}[1/p]}(G_k) \hookrightarrow \Z_{\mathbb{C}}(G_k)$$
and the compatibility of the Weil representation with scalar extension. 

Consider the subquotients $W_{n,m}^k$ of the filtration of Proposition \ref{filtration_weil_rep_proposition} when $R = \mathbb{C}$. Applying Lemma \ref{reduced_Jacobson_rings_lemma} to each direct factor of the center, we have $\varphi_{n,m}^k(z) = 0$ if and only if for any Zariski open dense subset $U$ we have $\eta(\varphi_{n,m}^k(z)) = 0$ for all $\eta : \Z_{\mathbb{C}}(G_k) \to \mathbb{C}$ in $U$. Then combining it with Corollary \ref{reduced_jacobson_big_modules_cor} and Proposition \ref{regular_representation_generic_quotient_semisimple_prop}, we are left to check $\eta(\varphi_{n,m}^k(z)) = 0$ where this scalar is the action of $z$ on: 
$$(W_{n,m}^k)_\eta = \textup{Ind}_{P_k^n}^{G_n}(\pi_\eta \otimes 1_{n-k}) \otimes_\mathbb{C} \textup{Ind}_{Q_k^m}^{G_m}(\pi_\eta^\vee \otimes 1_{m-k}).$$
In order to introduce a term coming from $\textup{Ind}_{P_{k+1}}^{G_n}$ and use our hypothesis $\varphi_{n,m}^{k+1}(z)=0$, we use induction in stages. First, embed: 
$$1_{n-k} \subset \textup{Ind}_{P_1^{n-k}}^{G_{n-k}}(1_1 \otimes 1_{n-(k+1)}).$$
Similarly for $m$. Remark that the induced representation above has always length $2$ by the theory of segments: it contains $1_{n-k}$ as well as an irreducible quotient $\sigma$. Now by transitivity of induction we can embed $(W_{n,m}^k)_\eta$ in:
\begin{equation} \label{equation_big_tensor_product_W_n_m_k} \textup{Ind}_{P_{k+1}^n}^{G_n}(\textup{Ind}_{P_k^{k+1}}^{G_{k+1}}(\pi_\eta \otimes 1_1) \otimes 1_{n-(k+1)}) \otimes \textup{Ind}_{Q_{k+1}^m}^{G_m}(\textup{Ind}_{Q_k^{k+1}}^{G_{k+1}}(\pi_\eta^\vee \otimes 1_1) \otimes 1_{m-(k+1)}).\end{equation}
Thanks to the theory of segments, we can always assume for all $\eta$ in our Zariski dense open set $U$ that $\textup{Ind}_{P_k^{k+1}}^{G_{k+1}}(\pi_\eta \otimes 1_1)$ and $\textup{Ind}_{Q_k^{k+1}}^{G_{k+1}}(\pi_\eta^\vee \otimes 1_1)$ are irreducible.
Therefore it defines an irreducible quotient of the regular representation: 
$$C_c^\infty(G_{k+1}) \twoheadrightarrow \textup{Ind}_{P_k^{k+1}}^{G_{k+1}}(\pi_\eta \otimes 1_1) \otimes \textup{Ind}_{P_k^{k+1}}^{G_{k+1}}(\pi_\eta \otimes 1_1)^\vee$$
which induces a quotient of $W_{n,m}^{k+1} = \textup{Ind}_{P_{k+1}^n \times Q_{k+1}^m}^{G_n \times G_m}(C_c^\infty(G_{k+1}) \otimes 1)$, namely:
\begin{equation} \label{equation_big_tensor_product_W_n_m_k+1} \textup{Ind}_{P_{k+1}^n}^{G_n}(\textup{Ind}_{P_k^{k+1}}^{G_{k+1}}(\pi_\eta \otimes 1_1) \otimes 1_{n-(k+1)}) \otimes \textup{Ind}_{Q_{k+1}^m}^{G_m}(\textup{Ind}_{P_k^{k+1}}^{G_{k+1}}(\pi_\eta \otimes 1_1)^\vee \otimes 1_{m-(k+1)}).\end{equation}
Note that $z$ acts as $\varphi_{n,m}^{k+1}(z)=0$ on $W_{n,m}^{k+1}$, so it acts trivially on the latter quotient (\ref{equation_big_tensor_product_W_n_m_k+1}) of $W_{n,m}^{k+1}$. We will now relate this quotient to our big tensor product (\ref{equation_big_tensor_product_W_n_m_k}) above.

It remains to check that the right-hand side of the tensor products (\ref{equation_big_tensor_product_W_n_m_k}) and (\ref{equation_big_tensor_product_W_n_m_k+1}) define the same supercuspidal support. It is easier to switch to normalized induction to read supercuspidal support as with the non-normalized version, one needs to introduce twists by modulus characters. Choosing a square root of $q$ in $\mathbb{C}$, we only need to compare the supercuspidal supports of: 
$$\textup{Ind}_{Q_k^{k+1}}^{G_{k+1}}(\pi_\eta^\vee \otimes 1_1) = i_{Q_k^{k+1}}^{G_{k+1}}(\delta_{Q_k^{k+1}}^{-\frac{1}{2}} (\pi_\eta^\vee \otimes 1)) \textup{ and } \textup{Ind}_{P_k^{k+1}}^{G_{k+1}}(\pi_\eta \otimes 1_1)^\vee = i_{P_k^{k+1}}^{G_{k+1}}(\delta_{P_k^{k+1}}^{\frac{1}{2}} (\pi_\eta^\vee \otimes 1_1)).$$
But $\delta_{Q_k^{k+1}} = \delta_{P_k^{k+1}}^{-1}$ as $Q_k^{k+1}$ is the opposite parabolic of $P_k^{k+1}$. So the supercuspidal supports are the same. To sum up what we have obtained: the scalar $\eta(\varphi_{n,m}^k(z))$ corresponds to the action of $z$ on $(W_{n,m}^k)_\eta$ and it is also equal to the scalar $\eta'(\varphi_{n,m}^{k+1}(z))$ corresponding to the action of $z$ on $(W_{n,m}^{k+1})_{\eta'}$ where $\eta' : \Z_\mathbb{C}(G_{k+1}) \to \mathbb{C}$ is the supercuspidal support associated to the irreducible representation $\textup{Ind}_{Q_k^{k+1}}^{G_{k+1}}(\pi_\eta^\vee \otimes 1)$. \end{proof}

We deduce from the previous proposition the existence of our morphism:
\begin{prop} \label{defining_theta_n_m_prop} There exists a unique $\mathbb{Z}[1/p]$-algebra morphism: 
$$\theta_{\ZZ[1/p]}^\# : \Z_{\mathbb{Z}[1/p]}(G_n) \to \Z_{\mathbb{Z}[1/p]}(G_m)$$
such that for all $z \in \Z_{\mathbb{Z}[1/p]}(G_n)$ we have $\varphi_{n,m}^m(z \otimes 1) = \varphi_{n,m}^m( 1 \otimes \theta_{\ZZ[1/p]}^\#(z))$. 
\end{prop}

\begin{proof}
Looking at the definition of $\varphi_{n,m}^m$, the homomorphism is uniquely determined by the action of $\Z_{\mathbb{Z}[1/p]}(G_n)$ on $W_{n,m}^m$ and the canonical identification $\textup{End}_{G_n \times G_m} (W_{n,m}^m) \simeq \Z_{\mathbb{Z}[1/p]}(G_m)$ from Proposition \ref{endomorphisms_of_W_k_nm_prop}. \end{proof}

\subsection{Preservation of depth} \label{DEPTH SECTION}

Let $G$ be any connected reductive $p$-adic group, let $\mathcal{B}(G, F)$ be the Bruhat--Tits building of $G$. For $x\in \mathcal{B}(G,F)$ Moy and Prasad defined in \cite{moy_prasad_unrefined} a decreasing filtration of the parahoric subgroup $G_x=G_{x,0}$ by open pro-$p$ subgroups $G_{x,r}$, $r\in \mathbb{R}_{\geq 0}$, and the jumps in the filtration form a discrete sequence $r_1,r_2,\dots$. There exists a finite set $\Sigma$ of ``optimal points'' $x$ for which the $r_i$ are rational. We fix such a $\Sigma$ for $G$, and thus also for its Levi subgroups, and restrict our attention to $x\in \Sigma$. Moy and Prasad defined subgroups $G_{x,r^+} = \bigcup_{s>r}G_{x,s}$ and a set of characters $\chi_{x,r_i}:G_{x,r_i}/G_{x,r_i^+}\to \ZZ[1/p,\zeta_p]^{\times}$, called unrefined minimal types of depth $r_i$. Any $V\in \Rep_{\C}(G)$ contains at least one $\chi_{x,r_i}$ for some $i$ and $x$, and $V$ is said to have depth $r_i$ if all the unrefined minimal types it contains have depth $r_i$. An irreducible $V$ has depth $r_i$ if and only if $r_i$ is minimal with the property that $\pi^{G_{x,r_i^+}}\neq 0$ for some $x$. 

One can construct (c.f. \cite[II.5.7a]{vig_book},\cite[Appendix A]{dat_finitude}) finitely generated projective objects $$Q_i=P(r_i) = \bigoplus_x\bigoplus_{\chi_{x,r_i}}\ind_{G_{x,r_i}}^{G}\chi_{x,r_i}.$$ Considered as objects in $\Rep_{\ZZ[1/p]}(G)$, the $P(r_i)$ induce a decomposition of the category:
$$\Rep_{\ZZ[1/p]}(G) = \prod_i \Rep_{\ZZ[1/p]}(G)_{r_i}$$
where $\Rep_{\ZZ[1/p]}(G)_{r_i}$ is the full subcategory of objects $V$ satisfying $\Hom_{{\ZZ[1/p]}[G]}(Q_{j}^n,V)=0$ for $j\neq i$. Any object $V\in \Rep_{\ZZ[1/p]}(G)$ decomposes as a direct sum $V = \bigoplus_iV_i$ where $V_i = \sum_{\phi\in \Hom(Q_i,V)} \text{Im}(\phi)$. 

Let $R$ be a $\ZZ[1/p]$-algebra. Since compact induction commutes with base change, $P(r_i)\otimes_{\ZZ[1/p]}R$ is a projective generator of the subcategory $\Rep_R(G)_{r_i}$ of objects satisfying the analogous property in $\Rep_R(G)$. When $R$ is any algebraically closed field containing $p^{-1}$, this is precisely the subcategory of objects having depth $r_i$ (\cite[II.5.7 Rem]{vig_book}).

Let $e_i^G$ denote the idempotent in $\Z_{\ZZ[1/p]}(G)$ defined by the depth-$r_i$ projector $V\mapsto V_i$.

\begin{lemma}\label{lem:depth_parabolic} Let $R$ be a $\ZZ[1/p]$-algebra. For all parabolics $P = M N$ in  $G$, the parabolic induction functor $\textup{ind}_P^G : \Rep_R(M) \to \Rep_R(G)$ preserves the depth. 
\end{lemma}

\begin{proof} 
This lemma was proved in the case $R = \C$ in \cite{moy_prasad_parabolic} and in \cite[II.5.12]{vig_book} when $R$ is a general algebraically closed field of characteristic different from $p$. It suffices to prove our statement over $R=\ZZ[1/p]$ because we have observed above that the depth-$r_i$ subcategory of $\Rep_{\ZZ[1/p]}(G)$ is preserved under any scalar extension $\ZZ[1/p]\to R$.

Consider $\ind_P^G(Q_i)$, where $Q_i$ is the progenerator of $\Rep_{\ZZ[1/p]}(M)_{r_i}$ considered above. Since parabolic induction preserves depth in the case where $R=\C$, $e_i^G$ acts on $\ind_P^G(Q_i\otimes \C) = \ind_P^G(Q_i)\otimes \C$ as the identity endomorphism. Since $\ind_P^G(Q_i)$ is torsion free, it is a submodule of $\ind_P^G(Q_i)\otimes \C$ so $e_i^G$ acts on $\ind_P^G(Q_i)$ as the identity. Similarly, $e_i^G$ acts by zero on $\ind_P^G(Q_j)$ for $j\neq i$. This shows that $\ind_P^G(Q_i)$ is in the depth-$r_i$ subcategory of $\Rep_{\ZZ[1/p]}(G)$.

Reciprocally if $V$ is in $\Rep_{\ZZ[1/p]}(M)_{r_i}$, there is a surjection from a direct sum of copies of $Q_i$ to $V$ in $\Rep_{\ZZ[1/p]}(M)$, in which case we have a surjection from a direct sum of copies of $\ind_P^G(Q_i)$ to $\ind_P^G(V)$, which shows that $\ind_P^G(V)$ is in $\Rep_{\ZZ[1/p]}(G)_{r_i}$. \end{proof}

Given groups $G$, $G'$, and an object $V$ in $\Rep_R(G\times G')$, we will denote by $e_i^{G}V$ the depth-$i$ summand of $V$ for the $G$-action and by $Ve_i^{G'}$ the depth-$i$ summand for the $G'$-action. Note that $e_i^G(Ve_{i'}^{G'})=(e_i^GV)e_{i'}^{G'}$ since the idempotents commute in $\End_{G\times G'}(V)$.

\begin{lemma}\label{lem:reg_left_right_action}
$e_i^{G_k}C_c^{\infty}(G_k) = C_c^{\infty}(G_k)e_i^{G_k}$.
\end{lemma}
\begin{proof}
Let $\rho_l$ denote the left translation action on $C_c^{\infty}(G_k)$ and $\rho_r$ the right translation action, respectively. We have an isomorphism of $R[G_k]$-modules
\begin{align*}
(C_c^{\infty}(G_k),\rho_l)&\overset{\sim}{\to}(C_c^{\infty}(G_k),\rho_r)\\
f&\mapsto f^{\vee}\ .
\end{align*}
Thus $e_i^{G_k}$ acts as the identity (respectively, zero) through the left action if and only if it acts the same way through the right action.
\end{proof}

\begin{lemma}\label{lem:depth_levi_factors}
Let $R$ be any $\ZZ[1/p]$-algebra and let $k\leq m\leq n$ be positive integers. Let $M^n_k\times M^m_k$ act on $C_c^{\infty}(G_k)\otimes 1$ as in Subsection~\ref{subsection:filtration_by_rank}. Then $e_i^{M^n_k}(C_c^{\infty}(G_k)\otimes 1)=(C_c^{\infty}(G_k)\otimes 1)e_i^{M^m_k}$.
\end{lemma}
\begin{proof} Since everything is torsion-free, it suffices to prove the lemma in the setting where $R=\C$. We claim the $M^n_k$-representation $(e_i^{G_k}C_c^{\infty}(G_k))\otimes 1$ has depth $r_i$. Let $M = M^n_k$ for ease of notation. Given $x\in \mathcal{B}(M,F)$, the Moy-Prasad filtration is given by $M_{x,r} = M\cap (G_n)_{x,r}$ after considering $x$ as a point in $\mathcal{B}(G_n,F)$ (\cite[Prop 13.2.5]{kaletha_prasad}) and $M_{x,r} = (G_k)_{y,r}\times (G_{n-k})_{z,r}$ for $y\in \mathcal{B}(G_k,F)$, $z\in \mathcal{B}(G_{n-k},F)$. Suppose $(e_i^{G_k}C_c^{\infty}(G_k))\otimes 1$ contains an unrefined minimal type $\chi$ of depth $r$. If $r=0$, then $\chi$ is the trivial character on a group $M_{x,0^+}$ and therefore $e_i^{G_k}C_c^{\infty}(G_k)$ has non-zero fixed vectors under $(G_k)_{y,0^+}$ \textit{i.e.} $e_i^{G_k}C_c^{\infty}(G_k)$ contains a depth 0 type. So $r_i=0$. If $r>0$ we have: $$\chi:M_{x,r}/M_{x,r^+}\to \C^{\times}$$
and we can factor $\chi$ as $\chi_1 \otimes \chi_2$. Since $\chi$ is minimal and $\chi_2$ must be the trivial character of $(G_{n-k})_{z,r}$ we have that $\chi_1$ must be an unrefined minimal type of $(G_k)_{y,r} / (G_k)_{y,r^+}$ because $r>0$. But all unrefined minimal types in $e_i^{G_k}C_c^{\infty}(G_k)$ have same depth $r_i$, which must be equal to $r$. As a result of the claim, we have $e_i^{M_k^n} (C_c^{\infty}(G_k))\otimes 1) = (e_i^{G_k}C_c^{\infty}(G_k))\otimes 1 = (C_c^{\infty}(G_k)e_i^{G_k})\otimes 1 =  (C_c^{\infty}(G_k))\otimes 1)e_i^{M_k^m}$ by Lemma~\ref{lem:reg_left_right_action}.
\end{proof}

\begin{lemma}
We have $e_i^{G_n}W^k_{n,m} = W^k_{n,m}e_i^{G_m}$.
\end{lemma}
\label{lem:subquotient_left_right_depth}
\begin{proof}
We begin by noting that
\begin{align*}
W^k_{n,m} &= \ind_{P^n_k}^{G_n}\left(\ind_{Q^m_k}^{G_m}(C_c^{\infty}(G_k)\otimes 1) \right)\\
&=\ind_{Q^m_k}^{G_m}\left(\ind_{P^n_k}^{G_n}(C_c^{\infty}(G_k)\otimes 1) \right),\ 
\end{align*}
Now Lemma~\ref{lem:depth_parabolic} (several times) combined with Lemma~\ref{lem:depth_levi_factors} gives
\begin{align*}
e_i^{G_n}W^k_{n,m} &= \ind_{P^n_k}^{G_n}\left(e_i^{M^n_k}\ind_{Q^m_k}^{G_m}(C_c^{\infty}(G_k)\otimes 1) \right)\\
&=\ind_{P^n_k}^{G_n}\left(\ind_{Q^m_k}^{G_m}\left(e_i^{M^n_k}(C_c^{\infty}(G_k)\otimes 1)\right) \right)\\
&=\ind_{Q^m_k}^{G_m}\left(\ind_{P^n_k}^{G_n}\left((C_c^{\infty}(G_k)\otimes 1)e_i^{M^m_k}\right) \right) = W^k_{n,m}e_i^{G_m} 
\end{align*}

\end{proof}

\begin{corollary}
Let $R=\ZZ[1/p]$. We have $e_i^{G_n}\omega_{n,m}^{\ZZ[1/p]} = \omega_{n,m}^{\ZZ[1/p]}e_i^{G_m}$. In particular, $$\theta^\#_{\ZZ[1/p]}(\Z_{\ZZ[1/p]}(G_n)_{r_i})\subset \Z_{\ZZ[1/p]}(G_m)_{r_i}.$$
\end{corollary}
\begin{proof}
In our previous notation, Lemma~\ref{lem:subquotient_left_right_depth} shows $\varphi_{n,m}^k(e_i^{G_n}\otimes 1 - 1\otimes e_i^{G_m})=0$, so the result follows from Propositions~\ref{hom_space_Wk_Wk'_prop} and~\ref{prop:inclusion_of_kernels}.
\end{proof}
Therefore we can set $\theta^\#_{\ZZ[1/p],r} : \Z_{\ZZ[1/p]}(G_n)_r \to \Z_{\ZZ[1/p]}(G_m)_r$ and we have $\theta_{\ZZ[1/p]}^\# = \prod_r \theta^\#_{\ZZ[1/p],r}$ according to the depth decomposition.

\subsection{Action of the center: the case of arbitrary $R$} \begin{corollary} Let $R$ be any $\mathbb{Z}[1/p]$-algebra. We have an inclusion of kernels:
$$\Ker(\varphi_{n,m}^m) \subset \dots \subset \Ker(\varphi_{n,m}^1) \subset \Ker(\varphi_{n,m}^0).$$ \end{corollary}

\begin{proof} The $R$-algebra $\Z_R(G_n)_r \otimes \Z_R(G_m)_r$ is generated by the image of $\Z_{\mathbb{Z}[1/p]}(G_n)_r \otimes \Z_{\mathbb{Z}[1/p]}(G_m)_r$ thanks to Corollary \ref{scalar_extension_center_cor}. As the maps $\varphi_{n,m}^k$ are compatible with scalar extension, and so is the Weil representation, the statement follows from the fact that it holds over $\mathbb{Z}[1/p]$. \end{proof}

For any $\ZZ[1/p]$-algebra $R$, the previous corollary allows us to define $\theta_{R,r}^\# = \theta_{\ZZ[1/p],r}^\# \otimes R$ and $\theta_R^\# = \prod_r\theta_{R,r}^\#$, and we obtain Theorem~\ref{thm_intro_theta_R} of the introduction. We re-state it as follows:
\begin{prop} \label{defining_theta_n_m_prop} There exists a unique $R$-algebra morphism: 
$$\theta_R^\# : \Z_R(G_n) \to \Z_R(G_m)$$
such that for all $z \in \Z_R(G_n)$ we have $\varphi_{n,m}^m(z \otimes 1) = \varphi_{n,m}^m( 1 \otimes \theta_R^\#(z))$. Furthermore this construction is compatible with scalar extension $R\to R'$. \end{prop}

\begin{proof}This holds thanks to the compatibility with scalar extension of the Weil representation together with Corollary \ref{scalar_extension_center_cor}. \end{proof}

\section{Theta correspondence and supercuspidal support}\label{section:theta_and_scs}

In this section, let $R$ be an algebraically closed field of characeristic $\ell \neq p$. Let $\omega_{n,m}$ be the Weil representation with coefficients in $R$ for the group $G_n \times G_m$ and assume that $m \leq n$.

\subsection{Banal theta correspondence} \label{banal_theta_corresp_sec} Recall that if the characteristic $\ell$ of $R$ does not divide the pro-orders of $G_m$ and $G_n$, we say that $\ell$ is banal with respect to $G_n$ and $G_m$. The following theorem constitutes the heart of the theta correspondence and has been proved by Roger Howe for complex coefficients $R=\mathbb{C}$ and by Alberto M\'inguez for any algebraically closed field $R$ of banal characteristic with respect to $G_n$ and $G_m$:

\begin{thm}[\cite{minguez_thesis},\cite{minguez_howe}] \label{banal_theta_correspondence_thm} Let $\pi \in \textup{Irr}_R(G_m)$. Then: 
\begin{enumerate} 
\item there exists a unique $\theta(\pi) \in \textup{Irr}_R(G_n)$ such that $\omega_{n,m} \twoheadrightarrow \theta(\pi) \otimes_R \pi$;
\item the map $\pi \mapsto \theta(\pi)$ thus defined is injective;
\item the quotient is multiplicity one \textit{i.e.} $\textup{dim}_{G_n \times G_m}(\omega_{n,m},\theta(\pi) \otimes_R \pi)= 1$. \end{enumerate} \end{thm}
\noindent Write $\textup{Irr}_R^\theta(G_n)$ to denote the image of the map $\theta$. Then the theorem asserts a bijection:
$$\textup{Irr}_R(G_m) \overset{\theta}{\simeq} \textup{Irr}_R^\theta (G_n).$$
The map $\theta$ of the theorem is reputed to be compatible with the supercuspidal support, which means there exists a map $\theta_{\textup{scs}} : \Omega_{\textup{scs}}(G_m) \to \Omega_{\textup{scs}}(G_n)$ such that the following diagram commutes:
$$\xymatrix{
		\textup{Irr}_R(G_m) \ar[r]^{\theta} \ar[d] & \textup{Irr}_R(G_n)  \ar[d] \\
		\Omega_{\textup{scs}}(G_m) \ar[r]^{\theta_{\textup{scs}}} & \Omega_{\textup{scs}}(G)}.$$
Similarly to the map $\theta$, we can denote by $\Omega_{\textup{scs}}^\theta(G)$ the image of $\theta_{\textup{scs}}$, which alternatively is the image of $\textup{Irr}_R^\theta(G_n)$ through the supercuspidal support. It is not \textit{a priori} clear whether the map $\theta_{\textup{scs}}$ thus defined is injective, however, we can deduce its injectivity from the following more precise description:
$$\begin{array}{ccc}
\tilde{\Omega}_R(G_m) & \to & \tilde{\Omega}_R(G_n) \\
(M,\rho) & \mapsto & (M \times T_{n-m},\rho^\vee \chi_M \otimes_R \chi_{T_{n-m}}) \end{array}$$
where $\tilde{\Omega}_R(G_m)$ is the set of supercuspidal pairs $(M,\rho)$ where $M$ is a Levi and $\rho$ is a supercuspidal representation of this Levi and the characters of $M$ and $T_{n-m}$ respectively are: 
$$\chi_M = |\cdot|^{-\frac{n-m}{2}} \textup{ and } \chi_{T_{n-m}} = |\cdot|_1^{(m+1-n)+\frac{(n-1)}{2}} \otimes_R \cdots \otimes_R |\cdot|_1^{\frac{(n-1)}{2}}.$$ 
This map is well-defined on the equivalence classes on each side, also called the association classes of supercuspidal pairs, and therefore defines a map:
$$\begin{array}{ccc}
\Omega_R(G_m) & \to & \Omega_R(G_n) \\
 ( M,\rho )_{\textup{scs}} & \mapsto & (M \times T_{n-m},\rho^\vee \chi_M \otimes_R \chi_{T_{n-m}})_{\textup{scs}}
 \end{array}$$
where we denote equivalence classes by \textup{scs} subscripts. This map is precisely $\theta_R$ and one can easily check that the explicit latter map is injective, so $\theta_R$ induces a bijection $\Omega_R(G_m) \simeq \Omega_R^\theta(G_n)$.

\subsection{Map of varieties} Another way to think about supercuspidal supports is to identify them with the points of the Bernstein center. For all irreducible $\pi \in \textup{Irr}_R(G)$, the Bernstein center acts as a character $\eta_\pi : z \in \Z_R(G) \mapsto z_\pi \in R$ thanks to Schur's lemma. It is a result of Vignéras that the equivalence relation on $\textup{Irr}_R(G)$ defined by ``having the same character'' agrees with supercuspidal support \textit{i.e.} $\eta_{\pi} = \eta_{\pi'}$ if and only if $\textup{scs}(\pi) = \textup{scs}(\pi')$. For $\mathfrak{m}$ a supercuspidal support, we denote by $\eta_\mathfrak{m}$ the associated character. Therefore we have a bijection:
$$\begin{array}{ccc} 
\Omega_R(G) & \to & \textup{Hom}_{R-\textup{alg}}(\Z_R(G),R) \\
\mathfrak{m} & \mapsto & \eta_{\mathfrak{m}} \end{array}.$$
As $\Omega_R(G)$ is naturally identified with $\textup{Spec}(\Z_R(G))(R)$, it can be endowed with a structure of affine scheme. When $\ell$ is banal with respect to $G$, one can describe the irreducible components of $\Z_R(G)$, which correspond to the set of primitive idempotents in $\Z_R(G)$. Each one of these irreducible components is also connected and finite type over $R$. When $\ell$ is non banal, these irreducible components are still known to be finite type over $R$ but they may fail to be reduced. The situation is far worse as there is no explicit description of these components purely in terms of representation theory.

With this point of view, one can ask about the algebraicity of the map $\theta_R$ defined above:
\begin{prop} Let $\ell$ be banal with respect to $G_n$ and $G_m$. The map $\theta_R$ induces a morphism of algebraic varieties $\Omega_R(G_m) \to \Omega_R(G_n)$. \end{prop}

\begin{proof} We decompose the map $\theta_R$ as the composition of:
$$( M,\rho )_{\textup{scs}} \mapsto (M,\rho^\vee)_{\textup{scs}} \textup{ and } (M,\rho)_{\textup{scs}} \mapsto (M \times T_{n-m},\rho \chi_M \otimes_R \chi_{T_{n-m}})_{\textup{scs}}.$$
The first one is algebraic as it corresponds to the contragredient involution on the center, which is an automorphism of $\mathcal{Z}_R(G_m)$. 

Regarding the second one, let $\mathfrak{s} \in \mathcal{B}_R(G_m)$ be an inertial support and let $(M,\sigma)$ be a supercuspidal pair such that $(M,\sigma)_\textup{scs} \in \mathfrak{s}$. Let $X_R(M)$ be the variety of unramified characters for the Levi $M$. The map:
$$\psi \in X_R(M) \mapsto (M,\sigma \psi) \in \Omega_R^\mathfrak{s}(G_m)$$
identifies $\Omega_R^\mathfrak{s}(G_m)$ with the quotient $X_R(M) / H_{(M,\sigma)}$ where $H_{(M,\sigma)}$ is the finite group corresponding to all characters $\psi \in X_R(M)$ such that $(M,\sigma \psi)_\textup{scs} = (M,\sigma)_\textup{scs}$. We have this relation if we can find $w \in N_{G_m}(M) / M$ such that $(\sigma \psi)^w \simeq \sigma \psi$. Similarly let $\mathfrak{s}' \in \mathcal{B}_R(G_n)$ such that $(M \times T_{n-m},\sigma \chi_M \otimes_R \chi_{T_{n-m}})_\textup{scs} \in \mathfrak{s}'$ and identify the variety $\Omega_R^{\mathfrak{s}'}(G_n)$ with the quotient $X_R(M \times T_{n-m})/H_{(M \times T_{n-m},\sigma \chi_M \otimes_R \chi_{T_{n-m}})}$. We can define the algebraic map:
$$(M,\sigma \psi) \mapsto (M \times T_{n-m},\sigma \psi \chi_M \otimes_R \chi_{T_{n-m}}).$$
In order for the map:
$$(M,\rho)_{\textup{scs}} \mapsto (M \times T_{n-m},\rho \chi_M \otimes_R \chi_{T_{n-m}})_{\textup{scs}}$$ 
to be algebraic it is sufficient to check whether the algebraic map:
$$(M,\sigma \psi) \mapsto (M \times T_{n-m},\sigma \psi \chi_M \otimes_R \chi_{T_{n-m}})$$
induces a map $\Omega_R^\mathfrak{s}(G_m) \to \Omega_R^{\mathfrak{s}'}(G_n)$ on quotients of $X_R(M)$ and $X_R(M \times T_{n-m})$. But for all $\psi \in H_{(M,\sigma)}$ we claim that:
$$(M \times T_{n-m} , \sigma \psi \chi_M \otimes_R \chi_{T_{n-m}})_\textup{scs} = (M \times T_{n-m} , \sigma \chi_M \otimes_R \chi_{T_{n-m}})_\textup{scs}.$$
Indeed our map is equivariant for:
$$w \in N_{G_m}(M) / M \to (w,\textup{Id}_{T_{n-m}}) \in N_{G_n}(M \times T_{n-m}) / (M \times T_{n-m})$$
in the sense that: 
$$(\sigma \psi \chi_M \otimes_R \chi_{T_{n-m}})^{(w,\textup{Id}_{T_{n-m}})} \simeq (\sigma \psi \chi_M)^w \otimes_R \chi_{T_{n-m}} \simeq (\sigma \psi)^w \chi_M \otimes_R \chi_{T_{n-m}}$$
where we have used $\chi_M^w = \chi_M$ for all $w \in N_{G_m}(M) / M$ because $\chi_M$ factors through the determinant over $M$. Therefore we obtain a map $\Omega_R^\mathfrak{s}(G_m) \to \Omega_R^{\mathfrak{s}'}(G_n)$. \end{proof}

By the banality assumption, the ring $\Z_R(G_m)$ is reduced, so $\theta_R$ is uniquely determined by its behavior on points. Between principal blocks this algebraic map:
$$(T_m,\psi) \mapsto (T_n, \psi^{-1} |\cdot|_m^{- \frac{n-m}{2}} \otimes_R |\cdot|_1^{(m+1-n)+\frac{(n-1)}{2}} \otimes_R \cdots \otimes_R |\cdot|_1^{\frac{(n-1)}{2}})$$
has already been studied by Rallis. Writing $T_n = T_m\times T_{n-m}$, we can describe the corresponding morphism $\theta_R^\#$ on coordinate rings as follows: first, define
$$(\theta_R^{(T_m,1_m)})^\# : R[X_1^{\pm 1},\cdots , X_n^{\pm 1}] \mapsto R[X_1^{\pm 1}, \cdots, X_m^{\pm 1}]$$
by sending $X_i$ to $q^{- \frac{n-m}{2}} X_i^{-1}$ if $1 \leq i \leq m$ and to $q^{(i-n)+\frac{n-1}{2}}$ if $m+1 \leq i \leq n$. This map is compatible with the permutation action on variables from $\mathfrak{S}_n$ and $\mathfrak{S}_m$ in the sense that there is a commutative diagram:
$$\xymatrix@R+2pc@C+2pc{
R[X_1^{\pm 1},\cdots , X_n^{\pm 1}] \ar[r]^{(\theta_R^{(T_m,1_m)})^\#}  &  R[X_1^{\pm 1},\cdots , X_m^{\pm 1}]  \\
R[X_1^{\pm 1},\cdots , X_n^{\pm 1}]^{\mathfrak{S}_n} \ar[r]^{(\theta_R^{(T_m,1_m)_\textup{scs}})^\#} \ar@{^{(}->}[u]  & R[X_1^{\pm 1},\cdots , X_m^{\pm 1}]^{\mathfrak{S}_m} \ar@{^{(}->}[u] 
}.$$
It can be easily checked that the $R$-algebra morphism obtained on invariants, which we refer to as Rallis' map, is surjective. This extends beyond the principal block thanks to the explicit description of Bernstein blocks in the banal setting:

\begin{prop} \label{surjectivity_of_theta_banal_prop} Under the banal assumption, the map $\theta_R^\# : \Z_R(G_n) \to \Z_R(G_m)$ is surjective. So $\theta_R$ is a closed immersion. \end{prop}

\begin{proof} We can explicitly write what the map $\theta_R$ is on Bernstein components. First of all, note that the morphism $a \in \mathbb{G}_{m,R} \to \lambda a^{-1} \in \mathbb{G}_{m,R}$ corresponds to the morphism of $R$-algebras $X \in R[X^{\pm 1}] \mapsto \lambda X^{-1} \in R[X^{\pm 1}]$. Recall the notations from the previous proof where $(M,\sigma)$ is a supercuspidal pair and $\mathfrak{s} \in \mathcal{B}_R(G_m)$ such that $(M,\sigma)_{\textup{scs}} \in \mathfrak{s}$, as well as the image $\mathfrak{s}' \in \mathcal{B}_R(G_n)$. We choose the supercuspidal support $(M \times T_{n-m},\sigma \otimes_R 1_{n-m})_\textup{scs} \in \mathfrak{s}'$ as our base point in $\Omega_R^{\mathfrak{s}'}(G_n)$. We can also assume that $(M,\sigma)$ is factored according to its principal unramified part in the sense that $(M,\sigma) = (M_0 \times T_k,\sigma_0 \otimes 1_k)$ where $\sigma_0$ does not contain any unramified character of a torus.

Consider the $R$-algebra morphism:
$$\bar{m} \otimes_R \bar{t}_{n-m} \in R[M/M^0] \otimes_R R[T_{n-m}/T_{n-m}^0] \mapsto \chi_{T_{n-m}}(\bar{t}_{n-m})\chi_M(\bar{m}) \bar{m}^{-1} \in R[M/M^0]$$
which corresponds to $(M,\sigma \psi) \mapsto (M,\sigma^\vee \psi^{-1} \chi_M \otimes_R \chi_{T_{n-m}})$. 
Thanks to the identification $T_{k+n-m} = T_k \times T_{n-m}$  through $t_{n-m+k} = (t_k , t_{n-m})$ where $k$ was our principal unramified index, we can rewrite it as $(M_0 \times T_k,\sigma_0 \psi_{M_0} \otimes_R \psi_k) \mapsto (M_0 \times T_{k+n-m} , \sigma^\vee \psi_0^{-1} \chi_{M_0} \otimes_R \chi_{T_{k+n-m}} \psi_k^{-1})$ for:
$$\chi_{M_0} = \chi_M|_{M_0} \textup{ and } \chi_{T_{k+n-m}} = \chi_{M}|_{T_k} \chi_{T_{n-m}} = |\cdot|_{T_k}^{- \frac{n-m}{2}} \otimes_R |\cdot|_1^{(m+1-n)+\frac{(n-1)}{2}} \otimes_R \cdots \otimes_R |\cdot|_1^{\frac{(n-1)}{2}}.$$
This corresponds on $R[M_0/M_0^0] \otimes_R R[T_{k+n-m} / T_{k+n-m}^0] \to R[M_0/M_0^0] \otimes_R R[T_k / T_k^0]$ to:
$$(\theta_R^{(M,\sigma)})^\# : \bar{m}_0 \otimes_R \bar{t}_{n-m+k} \mapsto \chi_{M_0}(\bar{m}_0) \bar{m}_0^{-1} \otimes_R \chi_{T_{n-m+k}}(\bar{t}_{n-m+k}) \bar{t}_k^{-1}.$$
The latter map is equivariant for the action of the groups:
$$H_{(M \times T_{n-m},\sigma^\vee \otimes_R 1_{n-m})} = H_{(M_0,\sigma_0)} \times H_{(T_{k+n-m},1_{k+n-m})} = H_{(M_0,\sigma_0)} \times \mathfrak{S}_{k+n-m}$$
and its subgroup:
$$H_{(M,\sigma)} =  H_{(M_0,\sigma_0)} \times \mathfrak{S}_k$$
obtained from the embedding $\mathfrak{S}_k = N_{G_k}(T_k) /T_k \subset \mathfrak{S}_{k+n-m} = N_{G_{k+n-m}}(T_{k+n-m})/T_{k+n-m}$. The map $(\theta_R^{(M,\sigma)})^\#$ is compatible to the action of these groups and induces a $R$-algebra morphism $(\theta_R^{(M,\sigma)_\textup{scs}})^\#$ for invariant subrings:
$$\xymatrix{
R[M_0/M_0^0] \otimes_R R[T_{k+n-m} / T_{k+n-m}^0] \ar[r]  &  R[M_0/M_0^0] \otimes_R R[T_k / T_k^0] \\
R[M_0/M_0^0]^{H_{(M_0,\sigma_0)}} \otimes_R R[T_{k+n-m} / T_{k+n-m}^0]^{\mathfrak{S}_{k+n-m}} \ar[r] \ar@{^{(}->}[u]  & R[M_0/M_0^0]^{H_{(M_0,\sigma_0)}} \otimes_R R[T_k / T_k^0]^{\mathfrak{S}_k} \ar@{^{(}->}[u] 
}.$$
Here $R[M_0/M_0^0]^{H_{(M_0,\sigma_0)}} \to R[M_0/M_0^0]^{H_{(M_0,\sigma_0)}}$ is an isomorphism induced by $\bar{m}_0 \mapsto \chi_{M_0} (\bar{m}_0) \bar{m}_0^{-1}$ whereas $R[T_{k+n-m} / T_{k+n-m}^0]^{\mathfrak{S}_{k+n-m}} \to R[T_k / T_k^0]^{\mathfrak{S}_k}$ is Rallis' map. \end{proof}

\begin{rem}
In the non-banal setting, point (1) of Theorem \ref{banal_theta_correspondence_thm} already fails, as we discussed in the introduction with the counterexample when $\ell|(q^n-1)$. The map 
$$\begin{array}{cccc}
\theta_R : & \Omega_R(G_m) & \to & \Omega_R(G_n) \\
& ( M,\rho )_{\textup{scs}} & \mapsto & (M \times T_{n-m},\rho^\vee \chi_M \otimes_R \chi_{T_{n-m}})_{\textup{scs}}
\end{array}$$
is still well-defined and presents a good candidate for a theta correspondence on the level of sets. However, the scheme $\Omega_R(G_n) = \textup{Spec}(\Z_R(G_n))$ is no longer reduced in the non-banal setting, so a morphism of schemes is not uniquely determined by its values on points. The strategy of Proposition~\ref{surjectivity_of_theta_banal_prop} will not work to prove the morphism $\theta_R^\# = \theta_{\ZZ[1/p]}^\#\otimes_{\ZZ[1/p]}R$ from Proposition \ref{defining_theta_n_m_prop} is surjective, because even though a surjective ring morphism $\Z_R(G_n)\to \Z_R(G_m)$ realizing the previous map on points might exist, it would no longer be the unique ring morphism realizing this ``good candidate'' on points. \end{rem}

\section{Finiteness and inductive relations}

\subsection{Finiteness of $\theta_R^\#$.}

In  Proposition \ref{defining_theta_n_m_prop}, we obtained the morphism of $\mathbb{Z}[1/p]$-algebras $\theta_{\ZZ[1/p]}^\#$ by considering the natural action of $\Z_{\mathbb{Z}[1/p]}(G_n)$ on $W_{n,m}^m = \textup{ind}_{P_m^n}^{G_n} ( C_c^\infty(G_m) \otimes 1)$. This moprhism of $\mathbb{Z}[1/p]$-algebras can also be interpreted, using \cite[Sec 4]{dhkm_finiteness}, in terms of Harish-Chandra morphisms:

\begin{lemma} \label{identity_with_HC_and_W_n_m_lem} Let $\sigma_m = C_c^\infty(G_m) \otimes 1 \in \textup{Rep}_{\mathbb{Z}[1/p]}(M_m^n \times G_m)$ as in Proposition \ref{filtration_weil_rep_proposition}. Then for all $z \in \mathcal{Z}_{\mathbb{Z}[1/p]}(G_n)$ we have:
$$z_{W_{n,m}^m} = \textup{ind}_{P_m^n}^{G_n}(HC(z)_{\sigma_m})$$
where $HC : \Z_{\mathbb{Z}[1/p]}(G_n) \to \Z_{\mathbb{Z}[1/p]}(M_m^n)$ is the Harish-Chandra morphism. \end{lemma}

\begin{proof} We have $W_{n,m}^m = \textup{ind}_{P_m^n}^{G_n} ( C_c^\infty(G_m) \otimes 1)$ and we can apply \cite[Th 4.1]{dhkm_finiteness}. \end{proof}

\begin{prop}\label{prop:theta_HC} There exists a surjective map $\alpha_{\sigma_m} : \Z_{\mathbb{Z}[1/p]}(M_m^n) \twoheadrightarrow \Z_{\mathbb{Z}[1/p]}(G_m)$ such that: 
$$\theta_{\ZZ[1/p]}^\# = \alpha_{\sigma_m} \circ HC.$$ \end{prop}

\begin{proof} We identify $\sigma_m$ and $C_c^\infty(G_m)$ in an obvious way, so $HC(z)_{\sigma_m} \in \textup{End}_{G_m \times G_m}(C_c^\infty(G_m))$. Let $\rho_l$ and $\rho_r$ be respectively the left and the right action on the regular representation $C_c^\infty(G_m)$. For $z \in \Z_{\mathbb{Z}[1/p]}(G_n)$, we have:
$$\rho_r(\theta_{\ZZ[1/p]}^\#(z))=\rho_l(\theta_{\ZZ[1/p]}^\#(z)^\vee) \in \textup{End}_{G_m \times G_m}(C_c^\infty(G_m)).$$
We can take $\alpha_{\sigma_m}$ to be the composition of $\textup{ev}_{\sigma_m} : \Z_{\mathbb{Z}[1/p]}(M_m^n) \to \textup{End}_{G_m \times G_m}(C_c^\infty(G_m))$, given by the action of $\Z_{\mathbb{Z}[1/p]}(M_m^n)$ on $\sigma_m$, with the isomorphism $\textup{End}_{G_m \times G_m}(C_c^\infty(G_m)) \to \Z_{\mathbb{Z}[1/p]}(G_m)$ that is the inverse of $z \in \Z_{\mathbb{Z}[1/p]}(G_m) \mapsto \rho_l(z) \in \textup{End}_{G_m \times G_m}(C_c^\infty(G_m))$. Note that the morphism $\textup{ev}_{\sigma_m}$ is surjective because $\Z_{\mathbb{Z}[1/p]}(M_m^n) \simeq \Z_{\mathbb{Z}[1/p]}(G_m) \otimes \Z_{\mathbb{Z}[1/p]}(G_{n-m})$. So the map $\alpha_{\sigma_m}$ is surjective and the equality of the proposition holds. \end{proof}

By the surjectivity of the map $\alpha_{\sigma_m}$ and the compatibility with scalar extension of our map $\theta_R^\#$, the finiteness property of Harish-Chandra morphisms \cite[Th 4.3]{dhkm_finiteness} implies:

\begin{corollary}\label{cor:finiteness} Let $R$ be a noetherian $\mathbb{Z}_\ell$-algebra. Then $\theta_R^\# : \Z_R(G_n) \to \Z_R(G_m)$ is finite. \end{corollary}

\subsection{Inductive relations} In Section \ref{banal_theta_corresp_sec} we have defined for fields $R$ of banal characteristic with respect to $G_n$ and $G_m$ some explicit maps between supercuspidal supports:
$$\begin{array}{ccc}
\Omega_R(G_m) & \to & \Omega_R(G_n) \\
 ( M,\rho )_{\textup{scs}} & \mapsto & (M \times T_{n-m},\rho^\vee \chi_M \otimes_R \chi_{T_{n-m}})_{\textup{scs}}
 \end{array}$$
and where the characters are explicit:
$$\chi_M = |\cdot|^{-\frac{n-m}{2}} \textup{ and } \chi_{T_{n-m}} = |\cdot|_1^{(m+1-n)+\frac{(n-1)}{2}} \otimes_R \cdots \otimes_R |\cdot|_1^{\frac{(n-1)}{2}}.$$
To keep track of $m$ and $n$, let us call this map $\theta_{R,n,m} : \Omega_R(G_m) \to  \Omega_R(G_n)$. It is a simple computation to check that these explicit maps give inductive relations such as:
$$\theta_{R,n,m} = \theta_{R,n,k} \circ \theta_{R,k,k} \circ \theta_{R,k,m} \textup{ for } m \leq k \leq n.$$
We prove that these relations also holds in families over the integral Bernstein centers:

\begin{prop}\label{prop:inductive_relations} Let $R$ be an arbitrary $\ZZ[1/p]$-algebra. For $m \leq k \leq n$ we have inductive relations:
$$\theta_{R,n,m}^\# = \theta_{R,k,m}^\# \circ \theta_{R,k,k}^\# \circ \theta_{R,n,k}^\#.$$
\end{prop}

\begin{proof} The relation can be checked directly when $R=\mathbb{C}$ by the explicit description of $\theta_{\C,n,m}$ as a morphism of varieties in the previous section. Moreover, the inclusion $\Z_{\mathbb{Z}[1/p]}(G_r) \hookrightarrow \Z_{\mathbb{C}}(G_r)$ gives the relation over $\mathbb{Z}[1/p]$ and therefore over any $\mathbb{Z}[1/p]$-algebra $R$ by Proposition \ref{defining_theta_n_m_prop}. \end{proof}

\section{Surjectivity of $\theta^\#_{R}$}
\label{sec:surjectivity}

This section is devoted to the proof of the following.
\begin{thm} \label{surjectivity_statement_section6_thm}
For any $\mathbb{Z}[1/p]$-algebra $R$, the morphism $\theta_R^\# : \Z_R(G_n)\to \Z_R(G_m)$ is surjective.
\end{thm}

First, we do several reduction steps. Since $\theta_{R,k,k}^\#$ is the duality isomorphism, Proposition~\ref{prop:inductive_relations} implies it is sufficient to prove the theorem for $n=m+1$, so we will assume for the rest of this section that $n=m+1$. 

Next, since $\theta_R^\#$ is the extension of scalars of $\theta_{\mathbb{Z}[1/p]}^\#$ to $R$, it suffices to prove surjectivity when $R=\mathbb{Z}[1/p]$. A $\mathbb{Z}[1/p]$-module $M$ is zero if and only if $M \otimes W(\overline{\mathbb{F}_{\ell}}) = 0$ for all $\ell \neq p$.  Applying it when $M$ is the cokernel of our map, we get: 
$$\textup{coker}(\theta^\#_{\mathbb{Z}[1/p]})=0 \textup{ if and only if } \textup{coker}(\theta^\#_{W(\overline{\mathbb{F}_{\ell}})})=0 \textup{ for all } \ell \neq p.$$
Thus it suffices to prove surjectivity when $R=W(\overline{\mathbb{F}_{\ell}})$ for $\ell \neq p$. To ease notation in this section, we will abbreviate
\begin{align*}
\Z_n &= \Z_{\wflb}(G_n)\\
\theta_{n,m}^\#&= \theta_{\wflb,n,m}^\#
\end{align*}
if it is clear from the context.

\begin{rem} A similar faithfully flat descent argument can be applied to obtain finiteness of Harish-Chandra morphisms over arbitrary $\mathbb{Z}[1/p]$-algebras. Indeed by \cite[Th 4.3]{dhkm_finiteness} they are finite for noetherian $\mathbb{Z}_\ell$-algebras, so we can deduce finiteness over arbitrary $\mathbb{Z}[1/p]$-algebras from the compatibility of the center of the category with scalar extension proven in the appendix. Therefore Corollary \ref{cor:finiteness} holds over arbitrary $\mathbb{Z}[1/p]$-algebras. We keep it in the current form as this improvement does not simplify later proofs and is a consequence of the surjectivity statement.  \end{rem}

\subsection{The case $n=2$}
In this case, there is a quick proof. For $\lambda$ in $F^{\times}$, the scalar matrix $t_{\lambda} = \left(\begin{smallmatrix} \lambda \\&\lambda\end{smallmatrix}\right)$ in the group center $Z(G_2)$ defines an element, which we will denote $z_{\lambda}$, of the categorical center $\Z_2$ (in fact $\Z_2^{\times}$) in a natural way. Given an object $(\pi,V)\in \Rep_{\wflb}(G_2)$, the action $z_{\lambda}|_V$ of $z_{\lambda}$ on $V$ is by the endomorphism $\pi(t_{\lambda})$. 

Let $\rho_l$ denote the left-translation action of $G_1$ on $C_c^{\infty}(G_1)$. Let $f$ be an element of the induced module $\ind_{P^2_1}^{G_2}(C_c^{\infty}(G_1)\otimes 1)$ \textit{i.e.} $f(tng) = \rho_l(t_1) \cdot f(g)$ where $t = \textup{diag}(t_1,t_2) \in T_2$, $n \in N_2$ and $g \in G_2$. We have, by the definition of parabolic induction, $$(z_{\lambda}\cdot f)(g) = f(g t_\lambda) = \rho_l(\lambda) \cdot f(g)\ ,\ \ \ \text{ for $g\in G_2$}.$$ On the other hand, by Proposition~\ref{prop:theta_HC}, we have 
$$(z_\lambda \cdot f)(g) = \rho_r(\theta^\#_{2,1}(z_{\lambda})) \cdot f(g)\ ,\ \ \ \text{ for $g\in G_2$}.$$
It follows that $\rho_r(\theta^\#_{2,1}(z_\lambda))=\rho_r(\lambda^{-1})$, where the equality takes place in the ring $$\End_{\wflb[G_1\times G_1]}(C_c^{\infty}(G_1)) \cong \Z_1.$$
Fix a depth $k \in \mathbb{N}$. By Section \ref{DEPTH SECTION}, the induction functor preservers depth and $\ind_{P^2_1}^{G_2}(C_c^{\infty}(G_1)\otimes 1)$ preserves depth in the sense that:
$$\theta^\#_{2,1}(e_k^{G_2}) = e_k^{G_1}$$
for the central idempotents defining the respective depth-$k$ subcategories. Therefore $\theta^\#_{2,1}$ can be written coordinate by coordinate, according to the depth, as a direct product of ring morphisms $e_k^{G_2} \Z_2 \to e_k^{G_1} \Z_1$ for each $k$. It is enough to prove that each one of these maps is surjective. Note that the reduction map $F^\times \to F^\times / (1 + \mathcal{P}_F^{k+1})$ induces a surjection of $\wflb$-algebras:
$$\wflb[F^\times] \twoheadrightarrow \wflb[F^\times / (1 + \mathcal{P}_F^{k+1})] = \bigoplus_{i=0}^k e_i^{G_1} \Z_1.$$
As a result of the relation $\rho_r(\theta^\#_{2,1}(z_\lambda))=\rho_r(\lambda^{-1})$ for all $\lambda \in F^\times$, the map $e_k^{G_2} \Z_2 \to e_k^{G_1} \Z_1$ induced by $\theta^\#_{2,1}$ is surjective. So $\theta^\#_{2,1}$ is surjective as well.

\subsection{Definition of gamma factors}
The main tool we use to prove surjectivity of $\theta_{n,n-1}^\#$ for $n\geq 3$ is the theory of gamma factors and converse theorems for $A[G_n]$-modules, where $A$ is a Noetherian $\wflb$-algebra. We briefly recall and consolidate the relevant aspects of the theory developed \cite{Moss16.1,matringe_moss,h_whitt,HMconverse}.

Fix a nontrivial character $\psi:F \to \wflb^{\times}$, and let $\psi_A:F\to A^{\times}$ be its scalar extenion to $A$. We will also use $\psi$ and $\psi_A$, respectively, to denote the corresponding characters on $N_n$ defined in the usual way. For a smooth $A[G_n]$-module $V$ we define the $N_n,\psi$-coinvariants as $V^{(n)} = V/V(N_n,\psi_A)$ where $V(N_n,\psi_A)$ is the $A$-submodule generated by the set $\{nv-\psi_A(n)v:n\in N_n,\ v\in V\}$. 

\begin{definition}
An $A[G_n]$-module $V$ is of Whittaker type if it is admissible, $A[G_n]$-finitely generated, and if $V^{(n)}$ is free of rank one as an $A$-module.
\end{definition}

Any $A[G_n]$ module $V$ of Whittaker type gives rise to a ring homomorphism $f_V:\Z_n\to A$ defined by sending $z\in \Z_n$ to the element of $A$ that gives the endomorphism $z^{(n)}$ of $V^{(n)}$ under the canonical isomorphism
$$\End_A(V^{(n)}) = A.$$ Note that the map $f_V$ factors through the natural action $\Z_n \to \End_{A[G_n]}(V)\to \End_A(V^{(n)})$ and is obtained by composing with the canonical identification. In particular, when Schur's lemma holds \textit{i.e.} $\textup{End}_{A[G_n]}(V) = A$, the map $f_V$ is canonically identified with the natural action of the center $\Z_n \to \textup{End}_{A[G_n]}(V)$.

If $V$ is of Whittaker type, by Frobenius reciprocity $\Hom_A(V^{(n)},A)\cong \Hom_{A[G_n]}(V,\Ind_{N_n}^{G_n}\psi_A)$ is a free module of rank one. We let $\Wh(V,\psi_A)$ denote the Whittaker space of $V$ with respect to $\psi_A$: it is the image of the map
\begin{align*}
\Hom_{A[G_n]}(V,\Ind_{N_n}^{G_n}\psi_A) \otimes V &\to \Ind_{N_n}^{G_n}\psi_A\\
\phi\otimes v &\mapsto \phi(v)
\end{align*}
In fact, $\Wh(V,\psi_A)$ is the image of any morphism $\phi$ that generates $\Hom_{A[G_n]}(V,\Ind_{N_n}^{G_n}\psi_A)$ as an $A$-module. We treat elements $W\in \Wh(V,\psi_A)$ as functions on $G_n$ valued in $A$. If $\lambda$ denotes the composite
$$\lambda: V\to V^{(n)}\to A,$$ then for $v\in V$, the associated Whittaker function $\phi(v) = W_v\in \Wh(V,\psi_A)$ is given by
$$W_v(g) = \lambda(gv),\ \text{ for }g\in G_n$$

The element $u = \textup{diag}(1,-1,\dots,(-1)^{n-1}) \in G_n$ normalizes $N_n$ and conjugates $\psi$ to $\psi^{-1}$. We can then define an $A[G_n]$-module isomorphism $\Ind_{N_n}^{G_n}\psi_A\cong \Ind_{N_n}^{G_n}\psi_A^{-1}$ by $W\mapsto W'$ where
$$W'(g) = W(u g ).$$ In particular, $V^{(n)}$ is isomorphic as an $A$-module to $V/V(N_n,\psi_A^{-1})$.

If $V$ is an $A[G_n]$-module we define ${^\iota V}$ to be the $A[G_n]$-module with underlying $A$-module $V$ but with $G_n$-action twisted by $\iota:g\mapsto  {^tg^{-1}}$. Let $w_n\in G_n$ be the matrix with $1$'s along the antidiagonal, we can define an $A[G_n]$-module isomorphism ${^\iota (\Ind_{N_n}^{G_n}\psi_A)}\cong \Ind_{N_n}^{G_n}\psi_A^{-1}$ by $W\mapsto \widetilde{W}$, where
$$\widetilde{W}(g) = W(w_n {^tg^{-1}}).$$

If $V$ is Whittaker type, it follows that ${^\iota V}$ is of Whittaker type with respect to the character $\psi^{-1}$, and since $({^\iota V})^{(n)}$ is isomorphic to $V^{(n)}$, we conclude ${^\iota V}$ is Whittaker type with respect to $\psi$. In terms of Whittaker models, the map $W\mapsto \widetilde{W}$ gives an $A[G_n]$-module isomorphism $\Wh(V,\psi_A)\cong \Wh({^\iota V},\psi_A^{-1})$.

Let $A$ and $A'$ be commutative finitely generated $\wflb$-algebras. We define the multiplicative subset $S$ of $(A\otimes_{\wflb} A')[X^{\pm 1}]$ consisting of the polynomials in $X$, $X^{-1}$ with first and last coefficient equal to $1$. The functional equation defining gamma factors takes place in the ring $$R = S^{-1}\left((A\otimes_{\wflb} A')[X^{\pm1 }] \right).$$

Let $m<n$ be positive integers, let $V$ be a Whittaker type $A[G_n]$-module and let $V'$ a Whittaker type $A'[G_m]$-module. Let $0\leq j\leq n-m-1$, and let $\M_{a, b}$ denote the set of $a\times b$-matrices with coefficients in $F$. For $W\in \Wh(V,\psi_A)$ and $W'\in \Wh(V',\psi_{A'}^{-1})$ we define
$$I(X,W,W';j) = \sum_{l\in \ZZ}c_l(W,W';j)X^l,$$ where 
$$c_l(W,W';j) = \int_{\M_{j,m}}\int_{N_{n-1}\backslash G_{n-1}^{(l)}}W
\begin{pmatrix}
g\\
x&I_j\\
&&I_{n-m-j}
\end{pmatrix}
\otimes W'(g)\ dg\ ,$$
and $G_{n-1}^{(l)}$ denotes the subset of $G_{n-1}$ consisting of matrices $g$ such that $v_F(\det g) = l$. We write $I(X,W,W') = I(X,W,W';0)$. 

Let $w_{t,r} = \text{diag}(I_t,w_r)$. By \cite[Th 3.2]{Moss16.1}, the formal series $I(X,W,W')$ in fact defines an element of $R$, so the formal series $I(q^{-1}X^{-1}, w_{m,n-m}\widetilde{W},\widetilde{W'};n-m-1)$ also defines an element of $R$.

\begin{thm}[\cite{matringe_moss} Cor 3.10]\label{thm:functional_eqn}
There is a unique element $\gamma(X,V\times V',\psi)$ of $R^{\times}$ such that for all $W\in \Wh(V,\psi_A)$, all $W'\in \Wh(V',\psi_{A'}^{-1})$, $$I(q^{-1}X^{-1}, w_{m,n-m}\widetilde{W}, \widetilde{W'};n-m-1) = \omega_{V'}(-I_{n-1})^{n-2} \gamma(X,V\times V',\psi)I(X,W,W'),$$ where $\omega_{V'}$ is the central character of $\Wh(V',\psi_{A'})$.
\end{thm}

\begin{rem}
\label{rem:n_is_3}
For a single case below, we will need the gamma factor when $V$ is an $A[G_1]$-module and $V'$ is an $A'[G_1]$-module (both Whittaker type), in which case $V$ and $V'$ are simply characters of $F^{\times}$ with values in $A^{\times}$ and $(A')^{\times}$, respectively. Here, the notions of Whittaker type and co-Whittaker are equivalent. In this case we define $\gamma(X,V\times V',\psi)$ to be the unique Godement--Jacquet gamma factor $\gamma(X, V\otimes V',\psi)\in R^{\times}$ satisfying the functional equation in \cite[Th 1.2]{moss_gj}, where $V\otimes V'$ denotes the diagonal tensor product $(A\otimes A')[G_1]$-module.
\end{rem}

\begin{corollary}\label{cor:gamma_inverse}
In the ring $R$, we have
$$\gamma(X,V\times V',\psi)^{-1}=\gamma(q^{-1}X^{-1},{^\iota V}\times {^\iota V'},\psi^{-1}).$$
\end{corollary}
\begin{proof}
The proof given in \cite[Corollary 5.6]{Moss16.1} works in this level of generality.
\end{proof}

The gamma factor is compatible with extension of scalars in the sense of \cite[Cor 3.11]{matringe_moss}, as we now explain. If $f:A\to B$ and $f':A'\to B'$ are ring homomorphisms, and we let  $f\otimes f'$ denote the homomorphism $R\to R'$ obtained by applying $f$ and $f'$ to the coefficients of polynomials, we have
\begin{equation}
\label{eqn:gamma_compatible_scalar_extension}
(f\otimes f')\left(\gamma(X,V\times V',\psi ) \right) = \gamma\left(X,(V\otimes_{A,f}B)\times (V'\otimes_{A',f'}B'),\psi\right).
\end{equation}

Let $a$ be an element of $A^{\times}$ and let $\chi_a:g\mapsto  a^{\text{val}_F(\det g)}$ be the corresponding unramified character on $G_n$. Unramified twisting shifts the variable $X$ in the gamma factor in the classical way:
\begin{lemma}\label{lem:unramified_twisting}
Given $a\in A^{\times}$ and $a'\in (A')^{\times}$, $$\gamma(X, \chi_a V\times \chi_{a'}V',\psi) = \gamma((a\otimes a')X, V\times V',\psi)\ .$$
\end{lemma}
\begin{proof}
Let $\pi:G_n\to \text{Aut}(V)$ denote the homomorphism by which $G_n$ acts on $V$. A Whittaker function $W$ is in $\Wh(\chi_a V,\psi_A)$ if and only if $W(g) = \chi_a(g)W_0(g)$ for $W_0$ in $\Wh(V,\psi_A)$. To see this, note that if $v$ is in $\chi V$ then its Whittaker function $W_v(g)$ is given by $$W_v(g) = \lambda((\chi\pi)(g)v) = \lambda(\chi(g)\pi(g)v) = \chi(g) \lambda(\pi(g)v) = \chi(g) W_{0,v}(g),$$ where $W_{0,v}$ is the Whittaker function associated to $v$ in the space $\Wh(V,\psi_A)$. It follows from the definition that $I(X,W_v,W') = I((a\otimes 1)X,W_{0,v},W')$, and the analogous property is true for $V'$ and $\chi_{a'}V'$. The lemma now follows from the uniqueness in Theorem~\ref{thm:functional_eqn}.
\end{proof}

Finally, we record a basic fact that we will use in the last subsection.

\begin{lemma}[\cite{HMconverse} Corollary 4.2]
\label{lem:rational_functions_finiteness}
Let $B$ and $B'$ be Noetherian $\wflb$-algebras, such that $B'$ is contained in $B$ and $B$ is finitely generated as a $B'$-module. Let $S'$ be the subset of $B'[X,X^{-1}]$ consisting of polynomials with first and last coefficient equal to 1. Then $(S')^{-1}(B'[X,X^{-1}])$ is the intersection of the subrings $B'[[X]][X^{-1}]$ and $S^{-1}(B[X,X^{-1}])$ in $B[[X]][X^{-1}]$.
\end{lemma}

\subsection{Co-Whittaker modules, the universal gamma factor, and a descent theorem}

We define
$$\Gamma_n = \ind_{N_n}^{G_n}\psi$$ and consider it as a $Z_n[G_n]$-module via the natural action of $Z_n$. In this context, we will need the following
\begin{thm}[\cite{h_whitt}]\label{thm:Gamma_cowhitt}
The $\Z_n[G_n]$-module $\Gamma_n$ is admissible over $\Z_n$, $\Gamma_n^{(n)}$ is free of rank one over $\Z_n$, and every nonzero quotient $Q$ of $\Gamma_n$ has $Q^{(n)}\neq 0$.
\end{thm}

This inspired the following definition in \cite{h_whitt}, where $A$ is a Noetherian $\wflb$-algebra:
\begin{definition}
An $A[G_n]$-module $V$ is co-Whittaker if it is admissible over $A$, if $V^{(n)}$ is free of rank one over $A$ and if every nonzero quotient $Q$ of $V$ satisfies $Q^{(n)}\neq 0$.
\end{definition}

It is proved in \cite[Proposition 6.2]{h_whitt} that co-Whittaker modules satisfy Schur's lemma. The property of being co-Whittaker is preserved under the operations of scalar extension and taking quotients. In fact $\Gamma_n$ is the universal co-Whittaker module in this sense:
\begin{thm}[\cite{h_whitt} Theorem 6.3]
If $V$ is a co-Whittaker $A[G_n]$-module, then $V$ is a quotient of $\Gamma_n\otimes_{\Z_n,f_V}A$.
\end{thm}

On one hand, note that a co-Whittaker $A[G_n]$-module $V$ is cyclic: any preimage of a generator of $V^{(n)}$ under the surjection $V\to V^{(n)}$ provides a generator for $V$. In particular it is Whittaker type.

On the other hand, for every Whittaker type $A[G_n]$ module $V$ we can construct a canonical submodule:
$$V_0 := \ker\left(V\to \prod_{\{U\subset V\ :\ U^{(n)} = V^{(n)}\}} V/U\right).$$ 
\begin{lemma}\label{lem:canonical_cowhitt_sub}
The submodule $V_0$ is co-Whittaker and $V_0^{(n)} = V^{(n)}$.
\end{lemma}
\begin{proof}
To see that $V_0^{(n)} = V^{(n)}$, we consider the restriction of $V$ and its submodules to the mirabolic subgroup $P_n$ of $G_n$ consisting of matrices with bottom row having the form $$(0,\dots, 0, 1).$$ In the category of $A[P_n]$-modules there is a natural injection of functors
$$\ind_{N_n}^{P_n}(\psi_{(-)^{(n)}}) \to \text{id},$$ whose image is the so-called Schwartz functions, and denoted $\S$ (to see this, use transitivity of parabolic induction to write, in the traditional Bernstein--Zelevinsky notation, $$\ind_{N_n}^{P_n}(\psi_{V^{(n)}})=(\Phi^+)^{n-1}\Psi^+(V^{(n)})\ ,$$ and apply the exact sequence in, \cite[Prop 5.12 (d)]{bz1}, c.f. \cite[Prop 3.1.3]{emerton_helm_families}). Each submodule $U$ appearing in the definition of $V_0$ satisfies $\S(V) = \S(U)\subset U$, hence $S(V)\subset V_0$. On the other hand, the composition
$$\left(\ind_{N_n}^{G_n}\psi_{V^{(n)}}\right)^{(n)} \overset{\sim}{\to} \S(V)^{(n)}\to V^{(n)}$$ is an isomorphism 
(\cite[Prop 3.1.3 and 3.1.5]{emerton_helm_families}), so the inclusions
$$\S(V) \subset V_0 \subset V$$  induce isomorphisms $\S(V)^{(n)} = (V_0)^{(n)} = V^{(n)}$. 

Now if $Q=V_0/U'$ is a quotient of $V_0$ with $Q^{(n)}=0$ then $(U')^{(n)}=(V_0)^{(n)}=V^{(n)}$ so $V_0$ is in the kernel of the map $V\to V/U'$ and $Q=0$. 
\end{proof}
From the proof of Lemma~\ref{lem:canonical_cowhitt_sub} we find that an equivalent construction of $V_0$ is given by taking the $A[G_n]$ submodule of $V$ generated by the $A[P_n]$-submodule $\S(V)$ of Schwartz functions.

Note also that $f_V = f_{V_0}$.

\begin{lemma}\label{lem:gamma_cowhitt_sub}
We have the following equality: $\gamma(X,V\times V',\psi) = \gamma(X, V_0\times V_0',\psi)$.
\end{lemma}
\begin{proof}
Since $\Wh(V_0,\psi_A)\subset \Wh(V,\psi_A)$, the two gamma factors satisfy the same functional equation for all $W$ in $\Wh(V_0,\psi_A)$, so the equality follows from the uniqueness in Theorem~\ref{thm:functional_eqn}.
\end{proof}

\begin{prop}
\label{prop:scalar_extension_gamma}
Let $V$ be a Whittaker type $A[G_n]$-module, let $V'$ be a Whittaker type $A'[G_m]$-module. Then
$$\gamma(X,V\times V',\psi) = (f_V\otimes f_{V'})\left(\gamma(X,\Gamma_n\times \Gamma_m,\psi)\right)$$ where we have used $f_V\otimes f_{V'}$ to also denote the map on Laurent series defined by applying $f_V\otimes f_{V'}$ to the coefficients.
\end{prop}
\begin{proof}
By \cite{Moss16.1} Corollary 5.5, it is true when $V$ and $V'$ are co-Whittaker. When they are Whittaker type, it is true for their co-Whittaker submodules $V_0$ and $V_0'$, which share the same gamma factor by Lemma~\ref{lem:gamma_cowhitt_sub}.
\end{proof}

We will make use of the following theorem, which states that the coefficients of the gamma factor provide descent data for Whittaker type representations.
\begin{thm}
\label{thm:gamma_descent}
Let $A$ and $A'$ be Noetherian $\wflb$-algebras and asssume $A$ is finitely generated as an $A'$-module. Let $m\geq 2$ be an integer and let $V$ be a Whittaker type $A[G_m]$-module such that for all primitive idempotents $e'$ of $\Z_{m-1}$ the coefficients of $\gamma(V\times e'\Gamma_{m-1},X^{-1},\psi)$ and $\gamma({^\iota V}\times e'\Gamma_{m-1},X,\psi^{-1})$ lie in $A'\otimes e' \Z_{m-1}$. Then the homomorphism $f_V:\Z_m\to A$ factors through $A'$.
\end{thm}
\begin{proof}
This is precisely Theorem 5.1 in \cite{HMconverse} except with the hypothesis ``$V$ is co-Whittaker'' relaxed to ``$V$ is Whittaker type.'' To prove it, use Lemma~\ref{lem:gamma_cowhitt_sub} to replace $V$ with its co-Whittaker submodule $V_0$, then apply Theorem 5.1 in \cite{HMconverse} to $V_0$, and recall that $f_V=f_{V_0}$.
\end{proof}

\begin{rem} We would like to apply this theorem in the case when $A = \Z_m$ is the Bernstein center, but $\Z_m$ is only \emph{locally} Noetherian. However, its connected components are Noetherian and we can always pass to a finite collection of such components by replacing $\Z_m$ by $e \Z_m$ for an appropriate idempotent element. To simplify notation in what follows, we will omit this choice of idempotent from the discussion.
\end{rem}

\subsection{Multiplicativity of the gamma factor}
In this subsection we prove a special case of the multiplicativity property for gamma factors in families, which we will use below. It is more convenient to work with normalized parabolic induction, so we will choose a square root of $q$ in an algebraic closure of $\text{Frac}(\wflb)$ and define $\wflb[\sqrt{q}]$. Since the extension $\wflb\to \wflb[\sqrt{q}]$ is faithfully flat (in fact, trivial unless $\ell=2$), the surjectivity of $\theta^\#_{\wflb[\sqrt{q}]}$ implies the surjectivity of $\theta^\#_{\wflb}$. For the sake of keeping straightforward notation, we will work with $\wflb$ below, with the understanding that when $\ell=2$, the ring $\wflb$ can be replaced with $\wflb[\sqrt{q}]$ without affecting the argument. Let $\delta_P = \delta_{P^n_{n-1}}$ denote the modulus character of the parabolic subgroup $P^n_{n-1}$. The normalized parabolic induction functor $\fri_{P^n_{n-1}}^{G_n}$ on $\Rep_{\mathbb{Z}[\sqrt{q}^{-1}]}(M^n_{n-1})$ is defined by first twisting by $\delta_P^{1/2}$ and then inducing. When $m=1$ we have the following

\begin{prop}\label{prop:multiplicativity_general}
Let $n\geq 3$ and $m\leq n-2$. Let $V_1$ be a Whittaker type $A_1[G_{n-1}]$-module, let $V_2$ be a Whittaker type $A_2[G_1]$-module, and let $V'$ be a Whittaker type $A'[G_m]$-module. Then $\fri_{P^n_{n-1}}^{G_n}(V_1\otimes V_2)$ is a Whittaker type $(A_1\otimes_{\wflb} A_2)[G_n]$-module and 
$$\gamma\left(X,\fri_{P^n_{n-1}}^{G_n}(V_1\otimes V_2)\times V',\psi\right) = \gamma(X,V_1\times V',\psi)\gamma(X, V_2\times V',\psi),$$ where the multiplication takes place in the ring $S^{-1}(D[X,X^{-1}])$ with $D = A_1\otimes A_2\otimes A'$ (c.f. Remark~\ref{rem:n_is_3} when $n=3$, $m=1$). The same equality holds with the lower parabolic $Q^n_{n-1}$ in place of $P^n_{n-1}$.
\end{prop}
\begin{proof}
Since $V_1$ is a Whittaker type $A_1[G_{n-1}]$-module, and parabolic induction preserves admissibility and finite generation, we need only show $\left(\ind_{P^n_{n-1}}^{G_n}(V_1\otimes V_2)\right)^{(n)}$ is free of rank one over $A_1\otimes A_2$.

Over $\Zl$, there is a so-called Leibniz rule for the Bernstein--Zelevinsky functors (\cite[III.1.10]{vig_book}). The proof given in \cite[III.1.10]{vig_book} follows the original proof of Bernstein--Zelevinsky in \cite[p.470-471]{bz2}, which boils down to the geometric lemma for normalized parabolic induction; it works equally well with $\wflb$ (or $\wflb[\sqrt{q}]$ if necessary) in place of $\C$ or $\Zl$. In our setting it simplifies to the following statement: if $\pi_1$ is a $\wflb[G_{n-1}]$-module and $\pi_2$ is a $\wflb[G_1]$-module, $\left(\fri_{P^n_{n-1}}^{G_n}(\pi_1\otimes \pi_2)\right)^{(n)}$ is the $\wflb$-module $\pi_1^{(n-1)}\otimes \pi_2^{(1)}$. Since $
V_1^{(n-1)}\cong A_1$ and $V_2^{(1)}\cong A_2$, we conclude $\fri_{P^n_{n-1}}^{G_n}\left(V_1\otimes V_2\right)^{(n)}\cong A_1\otimes A_2$. 

A priori, the above isomorphisms are only morphisms of $\wflb$-modules, however, for any $\wflb$-algebra $B$, the map $\psi_B$ factors through the structure morphism $N_n\to \wflb^{\times} \to B^{\times}$. Therefore for any $B[N_n]$-module $V$, the submodule module $V(N_n,\psi)$ equals $V(N_n,\psi_B)$, and it follows that the above isomorphisms are morphisms in the categories of $A_1$, $A_2$, and $A_1\otimes A_2$-modules, respectively.

By Proposition~\ref{prop:scalar_extension_gamma} above, it suffices to prove the multiplicativity property with $A_1 = \Z_{n-1}$, $V_1 = \Gamma_{n-1}$, $A_2 = \Z_1$, $V_2 = \Gamma_1$, $A' = \Z_m$ and $V' = \Gamma_m$.

Let $\Kbar$ denote the fraction field of $\wflb$ and choose an isomorphism $\Kbar\cong \mathbb{C}$. All representations are presumed smooth, so this isomorphism identifies $\Rep_{\Kbar}(G_n)$ with $\Rep_{\mathbb{C}}(G_n)$, which allows us to invoke results in \cite{JPSS2}. If $Q^n_{n-1}$ denotes the lower standard parabolic subgroup with Levi subgroup $M^n_{n-1}$, and we let $\pi_1$ be in $\Rep_{\Kbar}(G_{n-1})$ and $\pi_2$ be in $\Rep_{\Kbar}(G_1)$, with both $\pi_1$ and $\pi_2$ of Whittaker type, it is proven in \cite[Th 3.1]{JPSS2} that
$$\gamma\left(X,\fri_{Q^n_{n-1}}^{G_n}(\pi_1\otimes\pi_2) \times \tau, \psi\right) = \gamma(X,\pi_1\times \tau,\psi)\gamma(X,\pi_2\times \tau,\psi),$$ where $\tau\in \Rep_{\Kbar}(G_m)$ is Whittaker type.

Note the following isomorphisms:
\begin{align*}
{^\iota \left(\fri_{P^n_{n-1}}^{G_n}(\pi_1\otimes \pi_2))\right)} &\cong \fri_{Q^n_{n-1}}^{G_n}\left({^\iota (\pi_1\otimes \pi_2)}\right)\\
&\cong \fri_{Q^n_{n-1}}^{G_n}({^\iota \pi_1}\otimes {^\iota \pi_2}).
\end{align*}
Therefore, we can use Corollary~\ref{cor:gamma_inverse}, to deduce the multiplicativity property for the upper parabolic:
\begin{align*}
\gamma\left(X,\fri_{P^n_{n-1}}^{G_n}(\pi_1\otimes\pi_2) \times \tau, \psi\right) &= \gamma(q^{-1}X^{-1}, {^\iota \left(\fri_{P^n_{n-1}}^{G_n}(\pi_1\otimes\pi_2)\right)}\times {^\iota \tau},\psi^{-1})^{-1}\\
&= \gamma(q^{-1}X^{-1}, \fri_{Q^n_{n-1}}^{G_n}({^\iota\pi_1}\otimes{^\iota \pi_2})\times {^\iota \tau},\psi^{-1})^{-1}\\
&= \gamma(q^{-1}X^{-1},{^\iota \pi_1}\times {^\iota \tau},\psi^{-1})^{-1}\gamma(X,{^\iota \pi_2}\times {^\iota \tau},\psi^{-1})^{-1}\\
&= \gamma(X,\pi_1\times \tau,\psi)\gamma(X,\pi_2\times \tau,\psi)
\end{align*}

For a point $x\in \Spec(D)$, let $f_x:D\to \kappa(x)$ denote the corresponding ring homomorphism to the residue field $\kappa(x)$. If $V$ is a $D$-module, let $V_x$ denote the extension of scalars $V\otimes_{D,f_x}\kappa(x)$. The ring $D= \Z_{n-1} \otimes_{\wflb} \Z_1 \otimes_{\wflb} \Z_{n-2}$ is reduced and $\ell$-torsion free, (\cite[Lemma 5.1]{Moss16.1}). We will identify $\Spec(D)(\Kbar) = \Spec(\Z_{n-1})(\Kbar)\times \Spec(\Z_1)(\Kbar)\times \Spec(\Z_m)(\Kbar)$, and for $x = (x_1,x_2,x')\in \Spec(\Z_{n-1})(\Kbar)\times \Spec(\Z_1)(\Kbar)\times \Spec(\Z_m)(\Kbar)$, we decompose
$$f_x = f_{x_1}\otimes f_{x_2}\otimes f'.$$ Equation~(\ref{eqn:gamma_compatible_scalar_extension}) above and \cite[Th 3.1]{JPSS2} give the following equalities:
\begin{align*}
&(f_{x_1}\otimes f_{x_2}\otimes f_{x'})\left(\gamma\left(X,\fri^{G_n}_{P^n_{n-1}}(V_1\otimes V_2)\times V',\psi\right)\right)\\
&=\gamma\left(X,\fri^{G_n}_{P^n_{n-1}}(V_{1,x_1}\otimes V_{2,x_2})\times V'_{x'},\psi\right)\\
&=\gamma\left(X,V_{1,x_1}\times V'_{x'},\psi\right)\gamma\left(X,V_{2,x_2}\times V'_{x'},\psi\right)\\
&=(f_{x_1}\otimes f_{x_2}\otimes f_{x'})\left(\gamma\left(X,V_1\times V',\psi\right)\gamma\left(X,V_2\times V',\psi\right)\right)
\end{align*}
Now consider the difference of gamma factors:
$$\gamma\left(X,\fri^{G_n}_{P^n_{n-1}}(V_1\otimes V_2)\times V',\psi\right)-\gamma\left(X,V_1\times V',\psi\right)\gamma\left(X,V_2\times V',\psi\right).$$ Each of its coefficients is in the kernel of the homomorphism $f_x$ for every $x$ in $\Spec(D)(\Kbar)$. Since $D$ is reduced and $\ell$-torsion free, $\Spec(D)(\Kbar)$ is a Zariski dense subset of $\Spec(D)$, and 
$$\bigcap_{x\in \Spec(D)(\Kbar)}\ker(f_x)=0.$$ Thus the difference of gamma factors is zero.
\end{proof}

\begin{rem} Note we have only proved multiplicativity for parabolic induction from the Levi $M^n_{n-1}$. While a proof of the general case $P^n_k$ instead of $P_{n-1}^n$ might be within reach for $1\leq k<n$, the Leibniz rule for derivatives becomes more complicated. Since Proposition~\ref{prop:multiplicativity_general} is enough for our needs below, we will not pursue this further. \end{rem}

\subsection{Proof of surjectivity}

Recall the notation $f_V:\Z_n\to A$ for the action of $\Z_n$ on $V^{(n)}$, where $V$ is a Whittaker-type $A[G_n]$-module.

The action of $\Z_{n-1}$ on $\Gamma_{n-1}$ gives $\Gamma_{n-1}$ the structure of a Whittaker type $\Z_{n-1}[G_{n-1}]$-module (Theorem~\ref{thm:Gamma_cowhitt}). In this context, the corresponding homomorphism $f_{\Gamma_{n-1}}:\Z_{n-1}\to \Z_{n-1}$ is simply the identity map. Let $1$ denote the trivial character of $G_1$ over $\wflb$. Since $1^{(1)}=\wflb$, we have a natural ring homomorphism $f_1:\Z_1\to \wflb$.

The parabolically induced representation $\ind_{P^n_{n-1}}^{G_n}(\Gamma_{n-1}\otimes 1)$ realizes the theta correspondence in the following sense. On one hand, it has an action of $\Z_{n-1}$ via the composite 
$$\Z_{n-1}\overset{\sim}{\to} \End_{\wflb[G_{n-1}]}(\Gamma_{n-1})\hookrightarrow \End_{\wflb[G_n]}(\ind_{P^n_{n-1}}^{G_n}(\Gamma_{n-1}\otimes 1)).$$ On the other hand, $\ind_{P^n_{n-1}}^{G_n}(\Gamma_{n-1}\otimes 1)$ has an action of $G_n$, which gives it the structure of a smooth $\Z_{n-1}[G_n]$-module. By Proposition~\ref{prop:multiplicativity_general}, it is Whittaker type as a $\Z_{n-1}[G_n]$-module so we can consider the natural homomorphism
$$f_{\ind_{P_{n-1}^n}^{G_n}(\Gamma_{n-1}\otimes 1)}:\Z_n\to \Z_{n-1}.$$ By the following lemma, $f_{\ind_{P_{n-1}^n}^{G_n}(\Gamma_{n-1}\otimes 1)}$ is precisely $\theta_{n,n-1}^\#$.
\begin{lemma}\label{lem:theta_bernstein_action}
The natural $\Z_n$ action $$\Z_n\to \End_{\wflb[G_n]}(\ind_{P^n_{n-1}}^{G_n}(\Gamma_{n-1}\otimes 1))$$ factors through ${\theta}^\#_{n,n-1}:\Z_n\to \Z_{n-1}$.
\end{lemma}
\begin{proof}
Recall that $HC:\Z_n\to \Z_{\wflb}(M^n_{n-1})$ denotes the Harish-Chandra morphism defined by the equation
$$z_{\ind_{P^n_{n-1}}^{G_n}(V)} = \ind_{P^n_{n-1}}^{G_n}(HC(z)_V),\ \ V\in \Rep_{\wflb}(M^n_{n-1}),\ z\in \Z_n,$$ and Proposition~\ref{prop:theta_HC} states that ${\theta}^\#_{n,n-1} = \alpha \circ HC$ where $\alpha:\Z_{\wflb}(M^n_{n-1})\to \Z_{n-1}$ comes from the $\Z_{n-1}[G_{n-1}]$-linear action of $\Z_{\wflb}(M^n_{n-1})$ on $\Gamma_{n-1}$. 
\end{proof}

Our strategy is to use the multiplicativity property to compute the gamma factor $$\gamma\left(X,\ind_{P^n_{n-1}}^{G_n}(\Gamma_{n-1}\otimes 1)\times \Gamma_{n-2},\psi\right)$$ and apply Theorem~\ref{thm:gamma_descent} to relate the image of $f_{\ind_{P_{n-1}^n}^{G_n}(\Gamma_{n-1}\otimes 1)}$ to that of $f_{\Gamma_{n-1}}$.

To translate between normalized and non-normalized parabolic induction, observe that
$$\ind_{P^n_{n-1}}^{G_n}(\Gamma_{n-1}\otimes 1)= \fri_{P^n_{n-1}}^{G_n}\left(\delta_P^{-1/2}(\Gamma_{n-1}\otimes 1)\right)= \fri_{P^n_{n-1}}^{G_n} \left(|\det|^{-1/2}\Gamma_{n-1}\otimes |\cdot |^{\frac{(n-1)}{2}}\right)\ ,$$ so Proposition~\ref{prop:multiplicativity_general} and Lemma~\ref{lem:unramified_twisting} give the following decompositions of gamma factors:
\begin{align*}
&\gamma\left(X,\ind_{P_{n-1}^n}^{G_n}(\Gamma_{n-1}\otimes 1)\times \Gamma_{n-2},\psi\right)=\gamma\left(q^{1/2}X,\Gamma_{n-1}\times \Gamma_{n-2},\psi\right)\gamma\left(q^{-\frac{(n-1)}{2}}X,1\times \Gamma_{n-2},\psi\right)\\
&\gamma\left(X,{^\iota (\ind_{P_{n-1}^n}^{G_n}(\Gamma_{n-1}\otimes 1))}\times \Gamma_{n-2},\psi\right)=\gamma\left(q^{1/2}X,{^\iota \Gamma_{n-1}}\times \Gamma_{n-2},\psi\right)\gamma\left(q^{-\frac{(n-1)}{2}}X,1\times \Gamma_{n-2},\psi\right)
\end{align*}

Let $\swrz$ denote the image of ${\theta}^\#_{n,n-1}$. By Corollary~\ref{cor:finiteness}, $\Z_{n-1}$ is finitely generated as a module over $\swrz$ so we can apply Lemma~\ref{lem:rational_functions_finiteness} in what follows.

By Proposition~\ref{prop:scalar_extension_gamma}, the coefficients of $$\gamma\left(X,\ind_{P_{n-1}^n}^{G_n}(\Gamma_{n-1}\otimes 1)\times \Gamma_{n-2},\psi\right)\text{\ and\ }\gamma\left(X,\ind_{P_{n-1}^n}^{G_n}(\Gamma_{n-1}\otimes 1)\times {^\iota \Gamma_{n-2}},\psi\right)$$ are elements of $\swrz\otimes \Z_{n-2}$. Applying Corollary~\ref{cor:gamma_inverse}, we find that the coefficients of
$$\gamma\left(X,\ind_{P_{n-1}^n}^{G_n}(\Gamma_{n-1}\otimes 1)\times {^\iota \Gamma_{n-2}},\psi\right)^{-1} = \gamma\left(X,{^\iota \left(\ind_{P_{n-1}^n}^{G_n}(\Gamma_{n-1}\otimes 1)\right)}\times \Gamma_{n-2},\psi\right)$$ are also elements of $\swrz\otimes \Z_{n-2}$ (we have used Lemma~\ref{lem:rational_functions_finiteness} here with $B'=\mathcal{S} \otimes \Z_{n-2}$ and $B= \Z_{n-1} \otimes \Z_{n-2}$). Since $\swrz$ is a $\wflb$-algebra, and the coefficients of $\gamma\left(q^{-\frac{(n-1)}{2}}X,1\times \Gamma_{n-2},\psi\right)$ are in $\wflb\otimes \Z_{n-2}$, it follows that the coefficients of
$$\gamma\left(q^{1/2}X,\Gamma_{n-1}\times \Gamma_{n-2},\psi\right)\text{\ and } \gamma\left(q^{1/2}X,{^\iota \Gamma_{n-1}}\times \Gamma_{n-2},\psi\right)$$ must lie in the subring $\swrz\otimes \Z_{n-2}$ of $\Z_{n-1}\otimes \Z_{n-2}$. Shifting the variable $X$ by $q^{-1/2}$, we find the coefficients of $$\gamma\left(X,\Gamma_{n-1}\times \Gamma_{n-2},\psi\right)\text{\ and } \gamma\left(X,{^\iota \Gamma_{n-1}}\times \Gamma_{n-2},\psi\right)$$ also lie in $\swrz\otimes \Z_{n-2}$. Now we invoke Theorem~\ref{thm:gamma_descent} with $m=n-1$ to conclude the identity homomorphism
$f_{\Gamma_{n-1}}:\Z_{n-1}\to \Z_{n-1}$ factors through the inclusion $\swrz\subset \Z_{n-1}$. Therefore $\swrz = \Z_{n-1}$.

\section{Interpretation in terms of Galois parameters: an alternative proof of surjectivity in depth zero}
\label{sec:galois}

In this section we prove the claims made in Subsection~\ref{sec:functoriality_intro} regarding the explicit description of the map ${^L\theta_e}:X_m^e\sslash \widehat{G}_n \to X_n^e\sslash \widehat{G}_m$ and its being a closed immersion in the case $e=0$.

When $e=0$ we have given the description of ${^L\theta_e}$ in Subsection~\ref{sec:functoriality_intro}. When $e>0$, $X_n^e$ is the closed subscheme of $\widehat{G}_n^{k+1}$ representing tuples $(\Fam, \sigma_1,\dots, \sigma_k)$, subject to the relations defining the finitely presented subgroup $\langle \text{Fr},s_1,\dots, s_k\rangle \subset W_F/P_F^e$, where $s_1,\dots, s_k$ are any choice of topological generators of its normal subgroup $I_F/P_F^e$. The map ${^L\theta_e}$ is defined analogously to ${^L\theta_0}$ with $\sigma_1,\dots, \sigma_k$ replacing $\sigma$, as we now show. 

\begin{prop}\label{prop:description_of_Ltheta}
With notation as in Subsection~\ref{sec:functoriality_intro}, the map ${^L\theta_e}$ is the morphism on GIT quotients induced by the morphism $X^e_{m,\ZZ[\sqrt{q}^{-1}]}\to X^e_{n,\ZZ[\sqrt{q}^{-1}]}$ sending $(\Fam, \sigma_1,\dots, \sigma_k)$ to 
\begin{align*}
\left(
\begin{pmatrix}
q^{-\frac{n-m}{2}}I_m\cdot {^t\Fam^{-1}} \\ 
& q^{m+1 - n + \frac{n-1}{2}}\\
&&\ddots\\
&&&q^{\frac{(n-1)}{2}}
\end{pmatrix},
\begin{pmatrix}
{^t\sigma_1^{-1}}\\
&1\\
&&\ddots\\
&&&1
\end{pmatrix},\dots,
\begin{pmatrix}
{^t\sigma_k^{-1}}\\
&1\\
&&\ddots\\
&&&1
\end{pmatrix}\right)
\end{align*}
\end{prop}
\begin{proof}
The stated morphism $X_m^e\to X_n^e$ factors as the composite of:
\begin{enumerate}
\item The morphism
$$\Hom(W_F^0/P_F^e,\widehat{G}_m)_{\ZZ[\sqrt{q}^{-1}]}\to \Hom(W_F^0/P_F^e,\widehat{G}_n)_{\ZZ[\sqrt{q}^{-1}]}$$ induced by 
\begin{align*}
\widehat{G}_m&\to \widehat{G}_n\\
g&\mapsto \begin{pmatrix}{^tg^{-1}} & 0 \\ 0 & I_{n-m} \end{pmatrix}\ , 
\end{align*}
\item The action $$\Hom(W_F^0/P_F^e,\widehat{G}_n)_{\ZZ[\sqrt{q}^{-1}]}\to \Hom(W_F^0/P_F^e,\widehat{G}_n)_{\ZZ[\sqrt{q}^{-1}]}$$ of the unramified twisting given by sending $\Fam'\in \widehat{G}_n$ to
$$\begin{pmatrix}
q^{-\frac{n-m}{2}}
\cdot I_m\\ 
& q^{-m+1 - n + \frac{n-1}{2}}\\
&&\ddots\\
&&&q^{\frac{(n-1)}{2}}
\end{pmatrix}\cdot \Fam'\in \widehat{G}_n\ .$$
\end{enumerate}
Since the schemes $X_n^e\sslash \widehat{G}_n$ are known to be reduced and the $\C$-points are dense (\cite{dhkm_moduli}), it suffices to prove that our proposed morphism coincides with ${^L\theta_e}$ on $\C$-points. Since the theta correspondence induces the map on Weil--Deligne parameters $\phi\mapsto \hat{\phi}$, where $$\hat{\phi} := \phi^{\vee}\cdot\nu^{-\frac{n-m}{2}}\oplus \nu^{-m+1-n+\frac{n-1}{2}}\oplus \cdots \oplus \nu^{\frac{(n-1)}{2}},$$ it suffices to check that $\phi\mapsto \phi^{\vee}$ is induced on Weil--Deligne parameters by pushing forward along the group automorphism $\GL_n\to \GL_n$ sending $g$ to ${^tg^{-1}}$. But the local Langlands correspondence is known to be compatible with automorphisms (c.f. \cite[Prop 5.2.5]{haines}), and for an irreducible complex representation $\pi$ of $G_n$, the representation $g\mapsto \pi({^tg^{-1}})$ is known to be isomorphic to the contragredient $\pi^{\vee}$. It follows that $\phi^{\vee}$ is conjugate to $w\mapsto {^t\phi(w)^{-1}}$.
\end{proof}

The fact that ${^L\theta_e}$ is a closed immersion has already been established as a consequence of the local Langlands correspondence in families together with our result that $\theta_e$ defines a closed immersion. As noted in Subsection~\ref{sec:functoriality_intro}, it is natural to ask for a direct proof using geometry. We now provide such a proof in the depth zero case. 

\begin{prop} For all commutative $\mathbb{Z}[1/p]$-algebra $R$, for all $\lambda \in R^\times$ and all $A \in \widehat{G}_{n-m}(R)$, the closed immersion: 
$$(\mathcal{F},\sigma) \in \widehat{G}_{m,R} \times \widehat{G}_{m,R} \mapsto   \left( \left[ \begin{array}{cc}
\lambda I_m \cdot {}^t \mathcal{F}^{-1} & 0 \\
0 & A
\end{array} \right] , 
 \left[\begin{array}{cc}
{}^t \sigma^{-1} & 0 \\
0 & I_{n-m} \end{array} \right] \right) \in \widehat{G}_{n,R} \times \widehat{G}_{n,R}$$
induces a closed immersion $(\widehat{G}_{m,R} \times \widehat{G}_{m,R}) \sslash \widehat{G}_{m,R} \to (\widehat{G}_{n,R} \times \widehat{G}_{n,R} ) \sslash \widehat{G}_{n,R}$ on GIT quotients. \end{prop}

\begin{proof} First of all, transpose-inverse surely induces an isomorphism of schemes, and so does multiplication by $\lambda \in R^\times$, so it is enough to prove that the closed immersion:
$$f : (M_m,N_m) \in \widehat{G}_{m,R} \times \widehat{G}_{m,R} \mapsto   \left( \left[ \begin{array}{cc}
M_m & 0 \\
0 & A
\end{array} \right] , 
 \left[\begin{array}{cc} N_m & 0 \\
0 & I_{n-m} \end{array} \right] \right) \in \widehat{G}_{n,R} \times \widehat{G}_{n,R}$$
induces a closed immersion between GIT quotients.

Write $R[(m_{i,j})] \otimes_R R[(n_{i,j})]$ for the coordinate ring of $\widehat{G}_{n,R} \times \widehat{G}_{n,R}$ and let $M$ and $N$ be the two universal matrices with coefficients $(m_{i,j})$ and $(n_{i,j})$ respectively. For $\alpha$ a word in two letters, \textit{i.e.} an element of the free monoid with two generators, there is a corresponding $(M N)^\alpha$ in the matrices $M$ and $N$. The coefficients of these matrices belong to $R[(m_{i,j})] \otimes_R R[(n_{i,j})]$. Note that $(\widehat{G}_{n,R} \times \widehat{G}_{n,R}) \sslash \widehat{G}_{n,R} = \textup{Spec}(\mathcal{O}[\widehat{G}_{n,R} \times \widehat{G}_{n,R}]^{\widehat{G}_{n,R}})$ for simultaneous conjugation has coordinate ring generated by the coefficients of the characteristic polynomials $\chi_\alpha(X) = \textup{det}(X I_n - (M N)^\alpha)$ where $\alpha$ runs over all words in two letters. This result is true when $R = \mathbb{Z}$ by the work of Donkin \cite{donkin} (see alternatively \cite[App B, Lem B.9]{jantzen}) and one can deduce it for any $\mathbb{Z}$-algebra using \cite[Part I, Prop 4.18]{jantzen}. In particular this holds over arbitrary $\mathbb{Z}[1/p]$-algebras.

The $k$-th coefficient of $\chi_\alpha$ is denoted $c_k((MN)^\alpha)$. For a word $\alpha$ as above, we define $\alpha_1$ to be the number of occurrences of the first generator in $\alpha$. We obtain an interesting relation by applying $f^\#$ to these coefficients:
\begin{eqnarray*}
f^\#(\textup{det}(X I_n - (M N)^\alpha)) &=& \sum_{k=0}^n f^\#(c_k((M N)^\alpha)) X^k \\
 &=& \textup{det}(X I_m - (M_m N_m)^\alpha) \textup{det}(X I_{n-m} - A^{\alpha_1}) \\
 &=& \left( \sum_{k=0}^{n-m} c_k(A^{\alpha_1}) X^k \right) \left( \sum_{k=0}^m c_k((M_m N_m)^\alpha) X^k \right).
\end{eqnarray*}
Note that the highest coefficient of the characteristic polynomials is always $1$, the following one is minus the trace and the last one is the determinant up to a sign, so we can rewrite this system:
$$\left[ \begin{array}{cccccc}
1 & 0 &  &  \\
c_{n-m-1}(A^{\alpha_1}) & 1 & & \\
c_{n-m-2}(A^{\alpha_1}) & c_{n-m-1}(A^{\alpha_1}) & & \\
\vdots  & \vdots & \ddots & \\
 c_0(A^{\alpha_1}) & c_1(A^{\alpha_1}) &  \ddots & 1 \\
 & c_0(A^{\alpha_1}) & \ddots & \vdots \\
  &  & \ddots & \vdots  \\
  & & \ddots & c_1(A^{\alpha_1}) \\
  & & & c_0(A^{\alpha_1}) \end{array} \right] \left[ \begin{array}{c}
 1 \\
 c_{m-1}((M_m N_m)^\alpha) \\
 c_{m-2}((M_m N_m)^\alpha) \\
 \vdots \\
 c_1((M_m N_m)^\alpha) \\
 c_0((M_m N_m)^\alpha)
\end{array} \right]  = \left[ \begin{array}{c}
 1 \\
 f^\#(c_{n-1}((M N)^\alpha))) \\
 f^\#(c_{n-2}((M N)^\alpha))) \\
 \vdots \\
 \vdots \\
 \vdots \\
 f^\#(c_1((M N)^\alpha))) \\
 f^\#(c_0((M N)^\alpha)))
\end{array} \right].$$
In particular, this system is invertible because the left-hand side matrix has rank $m$. This result in the map $f^\#$ being surjective as all coefficients $c_k((M_m N_m)^\alpha)$ of the characteristic polynomial $\textup{det}(X I_m - (M_m N_m)^\alpha)$ belong to the image of $f^\#$ and these coefficients generate the ring $\mathcal{O}[\widehat{G}_{m,R} \times \widehat{G}_{m,R}]^{\widehat{G}_{m,R}}$. Therefore $f$ induces a closed immersion between GIT quotients. \end{proof}

In order to pullback the closed immersion to the space of parameters, we rely on the following Lemma, which is a consequence of the difficult Theorem VIII.0.2 in \cite{fargues_scholze}.

\begin{lemma} For all $\mathbb{Z}_\ell$-algebra $R$ with $\ell \neq p$, the closed immersion $X_{k,R}^0 \to \widehat{G}_{k,R} \times \widehat{G}_{k,R}$ induces a closed immersion $X_{k,R}^0 \sslash \widehat{G}_{k,R} \to (\widehat{G}_{k,R} \times \widehat{G}_{k,R}) \sslash \widehat{G}_{k,R}$. \end{lemma}

\begin{proof} By \cite[Th VIII.0.2]{fargues_scholze}, the ring of invariants $\mathcal{O}[X_{k,\mathbb{Z}_\ell}^0]^{\widehat{G}_{n,\mathbb{Z}_\ell}}$ is compatible with arbitrary scalar extension. Note that in our situation, their colimit can be taken over maps $F_2 \to W$ because $W = \langle \textup{Frob} , s \rangle$ discretises the tame quotient $W_F / P_F$. In particular, if we consider the map $\phi : F_2 = \langle f_1 , f_2 \rangle \to W = \langle \textup{Frob} , s \rangle$ sending $f_1$ to $\textup{Frob}$ and $f_2$ to $s$, all possible maps $\phi' : F_2 \to W$ factors through the latter in the sense that there exists $\psi : F_2 \to F_2$ such that $\phi \circ \psi = \phi'$. As a result the colimit must be a quotient of the ring of invariants of $\mathcal{O}[\widehat{G}_{k,R} \times \widehat{G}_{k,R}]^{\widehat{G}_{k,R}}$ that is associated to the map $\phi$. \end{proof}

We deduce from the proposition and the lemma:

\begin{thm} The map $\theta : X_m^0 \sslash \widehat{G}_m \to X_n^0 \sslash \widehat{G}_n$ is a closed immersion. \end{thm}

\begin{proof} We have to prove that the corresponding map $\theta^\# : \mathcal{O}[X_n^0]^{\widehat{G}_n} \to \mathcal{O}[X_m^0]^{\widehat{G}_m}$ is surjective. Because $\mathbb{Z}_\ell$ is flat over $\mathbb{Z}[1/p]$, the natural map $\mathcal{O}[X_m^0]^{\widehat{G}_m} \otimes_{\mathbb{Z}[1/p]} \mathbb{Z}_\ell \to \mathcal{O}[X_{n,\mathbb{Z}_\ell}^0]^{\widehat{G}_{m,\mathbb{Z}_\ell}}$ is an injection. Actually the latter map is an isomorphism because in the commutative diagram:
$$\xymatrix{
\mathcal{O}[\widehat{G}_m \times \widehat{G}_m]^{\widehat{G}_m} \otimes_{\mathbb{Z}[1/p]} \mathbb{Z}_\ell \ar[r] \ar[rd] &  \mathcal{O}[X_n^0]^{\widehat{G}_n} \otimes_{\mathbb{Z}[1/p]} \mathbb{Z}_\ell \ar[d] \\
 & \mathcal{O}[X_{n,\mathbb{Z}_\ell }^0]^{\widehat{G}_{n,\mathbb{Z}_\ell}}
}$$
we have $\mathcal{O}[\widehat{G}_m \times \widehat{G}_m]^{\widehat{G}_m} \otimes_{\mathbb{Z}[1/p]} \mathbb{Z}_\ell = \mathcal{O}[\widehat{G}_{m, \mathbb{Z}_\ell} \times \widehat{G}_{m, \mathbb{Z}_\ell}]^{\widehat{G}_{m, \mathbb{Z}_\ell}}$ and the map to $\mathcal{O}[X^0_{m,\mathbb{Z}_\ell }]^{\widehat{G}_{n,\mathbb{Z}_\ell }}$ is the surjective map from the lemma. Hence $\mathcal{O}[X^0_n]^{\widehat{G}_n} \otimes_{\mathbb{Z}[1/p]} \mathbb{Z}_\ell = \mathcal{O}[X^0_{n,\mathbb{Z}_\ell }]^{\widehat{G}_{n,\mathbb{Z}_\ell}}$.

Now the cokernel of $\theta^\#$ is a $\mathbb{Z}[1/p]$-module and the base change of $\theta$ to $\mathbb{Z}_\ell$ fits into a commutative diagram where all maps are known to be closed immersion thanks to the proposition and the lemma:
$$\xymatrix{
X^0_{m,\mathbb{Z}_\ell} \sslash \widehat{G}_{m,\mathbb{Z}_\ell} \ar[r]^{\theta_{\mathbb{Z}_\ell}} \ar[d] &  X^0_{n,\mathbb{Z}_\ell} \sslash \widehat{G}_{n,\mathbb{Z}_\ell}  \ar[d] \\
(\widehat{G}_{m,\mathbb{Z}_\ell} \times \widehat{G}_{m,\mathbb{Z}_\ell} )\sslash \widehat{G}_{m,\mathbb{Z}_\ell} \ar[r] & (\widehat{G}_{n,\mathbb{Z}_\ell} \times \widehat{G}_{n,\mathbb{Z}_\ell} )\sslash \widehat{G}_{n,\mathbb{Z}_\ell}  
}.$$
Therefore $\theta_{\mathbb{Z}_\ell}$ is a closed immersion. This implies that the cokernel of $\theta^\#$ must be trivial, otherwise there would exist an $\ell$ such that, by flatness, this cokernel does not vanish after base change to $\mathbb{Z}_\ell$ \textit{i.e.} $\theta_{\mathbb{Z}_\ell}$ is not a closed immersion. So $\theta$ is a closed immersion. \end{proof}

\section{Applications to the modular theta correspondence}\label{sec:applications}

We draw some conclusions from the previous sections for a modular theta correspondence. In this section $R$ is an algebraically closed field of characteristic $\ell$ (see \cite{trias_theta1} for considerations in the situation of non-algebraically closed fields). To give statements involving the theta correspondence in a more symmetric way, we will use unordered indices $k$ and $k'$ below, as opposed to $n$ and $m$, which always satisfy $m\leq n$ by assumption. We set $\omega = \omega_{k,k'}$.

\subsection{Finiteness of $\Theta$.} Let $\pi \in \textup{Irr}_R(G_k)$ be irreducible. The largest $\pi$-isotypic quotient $\omega_{\pi}$ of the Weil representation has a canonical decomposition $\omega_\pi = \pi \otimes_R \Theta(\pi)$ where $\Theta(\pi) \in \textup{Rep}_R(G_{k'})$. Similarly one can decompose the largest $\pi'$-isotypic quotient for $\pi' \in \textup{Irr}_R(G_{k'})$, thus defining $\Theta(\pi') \in \textup{Rep}_R(G_k)$. To speak of the cosocle $\theta(\pi)$ of $\Theta(\pi)$, we must prove $\Theta(\pi)$ is finite length. This will follow from a more general finiteness result in terms of characters of the Bernstein center. When $\eta : \Z_R(G_k) \to R$ is a character of the Bernstein center, we denote by $I_{\eta} = \langle z - \eta(z) \ | \ z \in \Z_R(G_k) \rangle$ its kernel.

\begin{definition} Set $\omega[\eta] = I_\eta \omega$. The largest $\eta$-isotypic quotient of $\omega$ is: 
$$\omega_\eta = \omega / \omega[\eta] \in \textup{Rep}_R(G_k \times G_{k'}).$$ \end{definition}

\noindent Equivalently, $\omega_{\eta}=\omega \otimes_{\Z_R(G_k),\eta} R$. Note that $\omega_{\eta}$ contains more information than $\omega_{\pi}$. Indeed, for $\pi \in \textup{Irr}_R(G_k)$, if we denote by $\eta_\pi : \Z_R(G_k) \to R$ the character induced by Schur's lemma, we have:
$$\omega_\eta \twoheadrightarrow \bigoplus_{\eta_\pi = \eta} \omega_\pi.$$
This direct sum is finite. Although we will not use this finiteness in what follows we explain it now for the sake of completeness. By the level decomposition \cite[App A.2]{dat}, one can define the depth $r$ associated to $\eta$ in an obvious way \textit{i.e.} $\eta(e_r) \neq 0$ for $e_r$ the central idempotent of the depth $r$ direct factor category. Moreover the depth $r$ subcategory of $\textup{Rep}_R(G_k)$ is equivalent, for a compact open $K$ that is small enough, to a direct factor of the category of modules over $\mathcal{H}_R(G_k,K)$ via the functor of $K$-invariants $V \mapsto V^K$. In particular $e_r \Z_R(G_k)$ becomes a direct factor of the center of $\mathcal{H}_R(G_k,K)$, and because Hecke algebras are finite over their centers \cite{dhkm_finiteness}, the algebra $A = \mathcal{H}_R(G_k,K) \otimes_{e_r \Z_R(G_k),\eta} R$ has finite dimension over $R$. Therefore there is a finite number of simple $A$-modules up to isomorphism, and by equivalence of categories, they are in bijection with the $\pi$ above so the direct sum is finite.

The fact that $\Theta(\pi)$ is finite length is an immediate consequence of the above factorization and the following proposition:

\begin{prop}\label{prop:omega_finite_length} For all $\eta : \Z_R(G_k) \to R$, the $(G_k \times G_{k'})$-representation $\omega_\eta$ is finite length. \end{prop}

Before proving the proposition, we introduce a key lemma, which makes the proposition a simple consequence of the properties of $\theta_R^\#$ and the fact that, over a field, admissible and finite type implies finite length (\cite[II.5.10]{vig_book}). In what follows we use the word \textit{locally} to mean ``after keeping a finite number of terms in the depth decomposition.''

\begin{lemma} \label{weil_rep_adm_fin_type_lem} The Weil representation $\omega \in \textup{Rep}_R(G_k \times G_{k'})$ is locally finitely generated and admissible over $\Z_R(G_k) \otimes_R \Z_R(G_{k'})$. \end{lemma}

\begin{proof} Given the rank-filtration of $\omega$ in Proposition \ref{filtration_weil_rep_proposition} --- it holds whether or not $k \geq k'$ --- it is enough to prove the result for all subquotients: 
$$\omega^{(i)} = \textup{ind}_{P_i^k \times Q_i^{k'}}^{G_k \times G_{k'}}(C_c^\infty(G_i) \otimes_R 1)$$
of the filtration.

First of all, induction preserves finite type \cite[Cor 1.5]{dhkm_finiteness} and the depth \cite[II.5.12]{vig_book} so it is enough to show that $C_c^\infty(G_i) \in \textup{Rep}_R(G_i \times G_i)$ is finite type after bounding the depth. We will consider $C_c^\infty(G_i)$ as the bi-module $\mathcal{H}_R(G_i)$ over the Hecke algebra $\mathcal{H}_R(G_i)$. Let $K_j = 1 + \varpi_F^{1+j} \mathcal{M}_i(\mathcal{O}_F)$ where $\mathcal{O}_F$ is the ring of integers of $F$ and $\varpi_F$ a unformizer. Let $e_{K_j}$ be the idempotent in the Hecke algebra $\mathcal{H}_R(G_i)$ associated to $K_j$. It exists because $K_j$ is a pro-$p$-group, \textit{i.e.}, has invertible pro-order. The depth $r$ subcategory of $\textup{Rep}_R(G_i)$ can be embedded in the category of modules over the relative Hecke algebra $\mathcal{H}_R(G_i,K_j)$ for $j$ big enough, indeed: the progenerator $P(r)$ of this direct factor category as defined in \cite[App A.2]{dat_finitude} is finitely generated, hence it is generated  by its $K_j$-fixed vectors for $K_j$ small enough. If we denote $e_r$ the central idempotent giving the depth $r$ subcategory and if $V \in \textup{Rep}_R(G_i)$ has depth $r$, then $V \mapsto e_{K_j} (e_r V)$ is an equivalence of category on a direct factor category of $\mathcal{H}_R(G_i,K_j)$-modules. 

Note that $e_r \mathcal{H}_R(G_i) =  \mathcal{H}_R(G_i) e_r$ because $e_r$ is a central, so this means that the depth $r$ part of $\mathcal{H}_R(G_i)$ agrees for both the left and right action of $G_i$. As a result, there exists a common $j$ for which the depth $r$ part is identified with $e_{K_j} e_r \mathcal{H}_R(G_i) e_{K_j}$ in the category of bi-modules over $\mathcal{H}_R(G_i,K_j)$. But $\mathcal{H}_R(G_i,K_j)$ is cyclic over $\mathcal{H}_R(G_i,K_j)$, both as a left or right module, so the direct factor $e_{K_j} e_r \mathcal{H}_R(G_i) e_{K_j}$ is \textit{a fortiori} cyclic as a bi-module. Therefore $\omega^{(i)}$ is locally finitely generated as a $(G_k \times G_{k'})$-representation with coefficients in $R$.

Next, as parabolic induction preserves admissibiliy \cite[II.2.1]{vig_book}, the finiteness of Harish-Chandra morphisms \cite[Th 4.3]{dhkm_finiteness} implies that $\omega^{(i)}$ is admissible over $\Z_R(G_k) \otimes_R \Z_R(G_{k'})$ if $C_c^\infty(G_i)$ is admissible over $\Z_R(G_i) \otimes_R \Z_R(G_i)$ -- but this is exactly \cite[Th 1.1]{dhkm_finiteness} \textit{i.e.} the finiteness of relative Hecke algebras over their center. \end{proof}

\begin{proof}[Proof of the proposition] If $k \leq k'$, the morphism $\theta_R^\# : \Z_R(G_{k'}) \to \Z_R(G_k)$ controls the action of $\Z_R(G_{k'})$ on the Weil representation by Propositon \ref{defining_theta_n_m_prop}. Therefore $\Z_R(G_{k'})$ acts as $\eta \circ \theta_R^\#$ on $\omega_\eta$. As $\omega_\eta$ is finite type, and admissible over $R$ by Lemma \ref{weil_rep_adm_fin_type_lem}, it is finite length.

If $k > k'$, the morphism $\theta_R^\# : \Z_R(G_k) \to \Z_R(G_{k'})$ controls the action of $\Z_R(G_k)$ on the Weil representation by Propositon \ref{defining_theta_n_m_prop}. In particular $\omega_\eta = 0$ if $\textup{ker}(\theta_R^\#) \not\subset I_\eta$ by maximality of $I_\eta$. When the inclusion holds, our surjectivity statement from Theorem \ref{surjectivity_statement_section6_thm} ensures that the action of $\Z_R(G_{k'})$ on $\omega_\eta$ must factor through the unique character $\eta'$ such that $\eta = \eta' \circ \theta_R^\#$. Again $\omega_\eta$ is finite type, and admissible over $R$ by Lemma \ref{weil_rep_adm_fin_type_lem}, so it is finite length. \end{proof}

\begin{rem} We could have carried out the proof of the proposition without relying on the surjectivity statement of Theorem \ref{surjectivity_statement_section6_thm}, using only the finiteness statement from Corollary \ref{cor:finiteness}. With this strategy, the unicity of $\eta'$ in the case $k > k'$ is not ensured, but the action of $\Z_R(G_{k'})$ always factors through a quotient $\Z_R(G_{k'})/I'$ that is an Artinian $R$-algebra, and admissibility over $R$ follows. In particular $I'$ is contained in a finite number of maximal ideals corresponding to the finitely many characters $\eta'$ such that $\eta = \eta' \circ \theta_R^\#$. \end{rem}

\subsection{Theta correspondence between characters of the center} Analyzing the proof of the proposition, we can derive a theta correspondence over $R$ between characters of the center. We come back to our original notation with $m\leq n$.

\begin{lemma} \label{defining_bijection_theta_R_scs_lem} We define $\theta_R$ on characters of $\Z_R(G_m)$ by $\eta \mapsto \theta_R (\eta) := \eta \circ \theta_R^\#$.
\begin{enumerate}[label=\textup{\alph*)}]
    \item For all $\eta : \Z_R(G_m) \to R$, we have $(\omega_{n,m})_\eta \neq 0$.
    \item We have $(\omega_{n,m})_\eta = (\omega_{n,m})_{\theta_R(\eta)}$ for all $\eta : \Z_R(G_m) \to R$.
    \item The set $\{ \eta' : \Z_R(G_n) \to R \  | \ (\omega_{n,m})_{\eta'} \neq 0 \}$ is the image of $\theta_R$.
\end{enumerate} \end{lemma}

\begin{proof} a) By the work of M\'inguez \cite{minguez_thesis}, we know that $(\omega_{n,m})_\pi \neq 0$ for all $\pi \in \textup{Irr}_R(G_m)$. As all characters $\eta$ can be realized as some $\eta_\pi$ (\cite[Cor 2.3]{h_whitt}), the result follows.

\noindent b)  By the surjectivity statement of Theorem \ref{surjectivity_statement_section6_thm}, $\theta_R^\#(I_{\theta_R(\eta)}) = I_\eta$. So by definition of $\theta_R^\#$ we obtain $\omega_{n,m}[\theta_R(\eta)] = \omega_{n,m}[\eta]$.

\noindent c) Because of the surjectivity statement, if $(\omega_{n,m})_{\eta'} \neq 0$, then $\Z_R(G_m)$ acts a character as well and this character $\eta$ is the unique one such that $\theta_R(\eta) = \eta'$. \end{proof}

We call the map $\theta_R$ just defined the theta correspondence between supercuspidal supports. It realizes a bijetion between characters of the center according to the lemma.

\begin{definition} Going back to the symmetric notation with $k$, $k'$, we denote $\theta_R$ this bijection in both ways:
$$\{ \eta : \Z_R(G_k) \to R \ | \ \omega_\eta \neq 0 \} \overset{\theta_R}{\simeq} \{ \eta' : \Z_R(G_{k'}) \to R \ | \ \omega_{\eta'} \neq 0 \}.$$ \end{definition}

When indexes $n \geq m$ are ordered, this map is:
$$\begin{array}{ccc}
\Omega_R(G_m) & \to & \Omega_R(G_n) \\
 ( M,\rho )_{\textup{scs}} & \mapsto & (M \times T_{n-m},\rho^\vee \chi_M \otimes_R \chi_{T_{n-m}})_{\textup{scs}}
 \end{array}.$$
Going the other way is simply taking the inverse of this map on its image.

This formulation in terms of characters of the Bernstein center has implications for $\Theta(\pi)$ as we can keep track of the action of the center of the category. Lemma \ref{defining_bijection_theta_R_scs_lem} implies the constituents of $\Theta(\pi)$ are strongly uniform in the following sense:

\begin{corollary} Let $\pi \in \textup{Irr}_R(G_k)$ and assume $\Theta(\pi) \neq 0$. Then $\Z_R(G_{k'})$ acts on $\Theta(\pi)$ via the character $\theta_R(\eta_\pi)$. In particular all constituents of $\Theta(\pi)$ have same supercuspidal support. \end{corollary}

\begin{proof} The first part is an immediate consequence of Lemma \ref{defining_bijection_theta_R_scs_lem}. The second part follows from \cite[Cor 12.12]{helm_bernstein_center}. \end{proof}

\appendix

\section{Geometric lemma} \label{geometric_lemma_appendix}

Let $G$ be a connected reductive group over $F$. We quickly explain why the geometric lemma still holds in the context of smooth representations with coefficients in $R$. Suppose we have fixed a minimal parabolic of $G$, say $P_0$, with Levi decomposition $P_0 = M_0 N_0$. A parabolic subgroup $P$ of $G$ is said to be standard if it contains $P_0$. All such parabolics subgroups $P$ come along with a standard Levi decomposition $M N$ where $M$ is the unique Levi in $P$ containing $M_0$. Let $P' = M' N'$ be an other standard parabolic. For $(\sigma,V) \in \textup{Rep}_R(M)$, we are going to give a filtration of the restriction-induction $\mathfrak{r}_G^{M'} \circ \mathfrak{i}_M^G(\sigma) \in \textup{Rep}_R(M')$. This filtration is famously known as the geometric lemma. In order to define it, we need to introduce the following subset of the Weyl group $W$ of $G$:
$$W^{M,M'} = \{ w \in W \ | \ w(M \cap P_0) \subset P_0 \textup{ and } w^{-1}(M' \cap P_0) \subset P_0 \}.$$
By \cite[II.1.2]{vig_book}, this set $W^{M,M'}$ also is a set of representatives for the double cosets $W_{M'} \backslash W / W_M$.

\subsection{Non-normalised geometric lemma} As we are not using normalised parabolic induction, because we are not assuming the existence of a square root of $q$ in $R$, we recall the version of the geometric lemma we use:
\begin{prop} \label{geometric_lemma_prop} There exists a filtration of $\mathfrak{r}_G^{M'} \circ \mathfrak{i}_P^G(\sigma) \in \textup{Rep}_R(M')$ whose subquotients $(I_w)_w$ are indexed by $W^{M,M'}$ and given by:
$$I_w \simeq \mathfrak{i}_{M' \cap w(M)}^{M'} \left( \delta_w \otimes_R ( w \circ \mathfrak{r}_M^{w^{-1}(M') \cap M}(\sigma)) \right)$$
where $\delta_w = \delta_{N'} / \delta_{N' \cap w(P)}$ is a character of $M' \cap w(M)$. \end{prop}

We will not prove this proposition, but we refer to the many references \cite{bz2,vig_book,ren} for expositions on the geometric lemma. However, the most suitable reference to deal without normalization seems to be the notes \cite[Sec  6]{casselman}. We simply point out the precise results we need and their proofs go along the same way as in the notes. Let $\Omega_w$ be the double coset in $P \backslash G / P'$ associated to $w \in W^{M,M'}$. Choose a total order $<$ on $P \backslash G / P'$, or equivalently on $W^{M,M'}$, such that $U_w = \cup_{w'<w} \Omega_{w'}$ is an open subset of $G$ for all $w \in W^M$. Denote the submodule of functions supported on $U_w$ by $\mathfrak{i}_{U_w} = \{ f \in \mathfrak{i}_M^G(\sigma) \ \vert \ \textup{supp}(f) \subset U_w\}$ and define $\mathfrak{j}_w = \mathfrak{i}_{U_w \cup \Omega_w} / \mathfrak{i}_{U_w}$. Then as in \cite[Prop 6.3.2]{casselman}, we have in $\textup{Rep}_R(P)$:
$$\mathfrak{j}_w \simeq \textup{ind}_{P' \cap w(P)}^{P'}(w \circ \sigma)$$
and the computation of its $N$-coinvariants \cite[Props 6.2.1 \& 6.3.3]{casselman} is still valid so $(J_w)_N$ is the representation $I_w$ we gave above.

\subsection{Maximal parabolics for general linear groups} The general linear group $G_n = \textup{GL}_n(F)$ is a connected reductive group over $F$. We choose as a minimal parabolic subgroup of $G_n$, also called a Borel subgroup in this situation, the subgroup of upper triangular matrices $B_n$ with Levi decomposition $T_n N_n$ where $T_n$ is the subgroup of diagonal matices in $G_n$ and $N_n$ the set of unipotent matrices in $B_k$. For $0 \leq k \leq n$, set:
$$M_k^n = \left\{ \left. \left[ \begin{array}{cc}
a_k & 0 \\
0 & b_{n-k}
\end{array} \right] \in G_n \ \right\vert \ a_k \in G_k \textup{ and } b_{n-k} \in G_{n-k} \right\}.$$
It is a standard Levi of $G_k$ and it is contained in a unique standard parabolic subgroup denoted by $P_k^n = M_k^n N_k^n$. For $k \leq  k' \leq n$, similarly write $P_{k'}^n = M_{k'}^n N_{k'}^n$.

Identify $W$ and the permutation matrices representing $\mathfrak{S}_n$. By setting $r = \textup{max}(0, k - (n - k'))$, the map below induces an isomorphism between $W_{M_k^n} \backslash W / W_{M_{k'}^n} \simeq [ \! [ r , k ] \! ]$:
$$\begin{array}{ccc}
W & \to & [ \! [ r , k ] \! ] \\
\sigma & \mapsto & | \{ \sigma(1) , \dots, \sigma(k) \} \cap \{ 1, \dots, k' \} | \end{array}.$$
As a result of this isomorphism, the set of representatives $W^{M_k^n,M_{k'}^n}$ is in bijection with the set $W(k,k',n) = \{ w_{k,k',i}^n \ | \ 0 \leq i \leq \textup{min}(k,n-k') \}$ where:
$$w_{k,k',i}^n = \left[ \begin{array}{ccc}
\textup{Id}_{k-i} & & \\
 & w_i & \\
 & & \textup{Id}_{n-k'-i}
\end{array} \right] \in G_n \textup{ with } w_i = \left[ \begin{array}{cc}
0 & \textup{Id}_i \\
\textup{Id}_{k'-k+i} & 0 \end{array} \right] \in G_{k'- k + 2i}.$$
These elements all satisfy $w_{k,k',i}^n(M_k^n \cap B_n) \subset B_n$ and ${(w_{k,k',i}^n)}^{-1}(M_{k'}^n \cap B_n) \subset B_n$.
One has:
$$M_{k'}^n \cap w_{k,k',i}^n(M_k^n) = \left[ \begin{array}{cccc}
G_{k-i} & & & \\
 & G_{k'-k+i} & & \\
 & & G_i & \\
 & & & G_{n-k'-i}
\end{array} \right]$$
that we denote $M_{(k-i,k'-k+i,i)}^n$ and:
$$M_k^n \cap {(w_{k,k',i}^n)}^{-1}(M_k^n) = \left[ \begin{array}{cccc}
G_{k-i} & & & \\
 & G_i & & \\
 & & G_{k'-k+i} & \\
 & & & G_{n-k-i}
\end{array} \right]$$
denoted by $M_{(k-i,i,k'-k+i)}^n$.

\begin{rem} When $k'=k$ above, we write $w_{k,i}^n$ and $W(k,n)$ for short. The situation becomes simpler as the element $w_{k,i}^n$ has order at most $2$ and is equal to its inverse. \end{rem}

\subsection{Comparing $H$-induced endomorphisms and $G$-endomorphisms.} Let $G$ be a locally profinite group. Let $H$ be a closed subgroup of $G$. In particular $H$ is a locally profinite group as well. Let $V \in \textup{Rep}_R(H)$. For $f \in \textup{ind}_H^G(V)$ and $\varphi \in \textup{End}_H(V)$, define $\textup{ind}_H^G(\varphi) \cdot f \in \textup{ind}_H^G(V)$ by $(\textup{ind}_H^G(\varphi) \cdot f) (g) = \varphi(f(g))$ for all $g \in G$. Then it easy to see that:

\begin{lemma} The map $\varphi \in \textup{End}_H(V) \mapsto \textup{ind}_H^G(\varphi) \in \textup{End}_G( \textup{ind}_H^G(V))$ is an injective morphism of algebras and the evaluation map $ev_1 : f \in \textup{ind}_H^G(V) \mapsto f(1_G) \in V$ induces a commutative diagram:
$$\xymatrixcolsep{5pc}\xymatrix{
		\textup{ind}_H^G(V) \ar[r]^{\textup{ind}_H^G(\varphi)} \ar[d]^{ev_1} & \textup{ind}_H^G(V) \ar[d]^{ev_1} \\
		V \ar[r]^\varphi & V
		}.$$ \end{lemma}

We are specifically interested in situations when the previous injective map becomes an isomorphism, giving a canonical identification between $\textup{End}_H(V)$ and $\textup{End}_G(\textup{ind}_H^G(V))$.

\begin{corollary} \label{all_morphism_are_coming_from_H_cor1} Suppose that $\textup{Hom}_H(\textup{ker}(ev_1)),V)= 0$. Then the map $\varphi \mapsto \textup{ind}_H^G(\varphi)$ above is an isomorphism and has inverse $\Phi \mapsto \tilde{\Phi}$ where $\tilde{\Phi}$ is, for $\Phi \in \textup{End}_G(\textup{ind}_H^G(V))$, the unique element in $\textup{End}_H(V)$ such that the following diagram commutes:
$$\xymatrixcolsep{5pc}\xymatrix{
		\textup{ind}_H^G(V) \ar[r]^{\Phi} \ar[d]^{ev_1} & \textup{ind}_H^G(V) \ar[d]^{ev_1} \\
		V \ar[r]^{\tilde{\Phi}} & V
		}.$$
In particular $\Phi = \textup{ind}_H^G(\tilde{\Phi})$.
\end{corollary}
		
\begin{proof} First of all we have that $\textup{End}_G(\textup{ind}_H^G(V)) \subset \textup{Hom}_G(\textup{ind}_H^G(V),\textup{Ind}_H^G(V))$ as the inclusion of induced representations $\textup{ind}_H^G(V) \subset \textup{Ind}_H^G(V)$ holds. Using Frobenius reciprocity, we get:
$$\textup{Hom}_G(\textup{ind}_H^G(V),\textup{Ind}_H^G(V)) \simeq \textup{Hom}_H( \textup{ind}_H^G(V) , V).$$
The exact sequence $0 \to \textup{ker}(ev_1) \to \textup{ind}_H^G(V) \to V \to 0$ gives by right exactness of $\textup{Hom}_H(-,V)$:
$$\textup{Hom}_H(\textup{ker}(ev_1),V) \to \textup{Hom}_H(\textup{ind}_H^G(V),V) \overset{p}{\to} \textup{End}_H(V) \to 0.$$
As $\textup{Hom}_H(\textup{ker}(ev_1),V)=0$, the map $p$ must be an isomorphism. So on the one hand we have that $\varphi \in \textup{End}_H(V) \mapsto \varphi \circ ev_1 \in \textup{Hom}_H(\textup{ind}_H^G(V),V)$ is an isomorphism. On the other hand, the isomorphism coming from adjunction is $\psi \in \textup{Hom}_H(\textup{ind}_H^G(V),V) \mapsto A_\psi \in \textup{Hom}_G(\textup{ind}_H^G(V),\textup{Ind}_H^G(V))$ with $A_\psi(f) : g \mapsto \psi(g \cdot f))$ for $f \in \textup{ind}_H^G(V)$. Gathering together the previous two isomorphims yields an isomorphism:
$$\varphi \in \textup{End}_H(V) \mapsto A_{\varphi \circ ev_1} \in \textup{Hom}_H(\textup{ind}_H^G(V),\textup{Ind}_H^G(V)).$$
But the image of $A_{\varphi \circ ev_1}$ is included in $\textup{ind}_H^G(V)$. Indeed, we have $ev_1 (g \cdot f) = f(g)$ for $f \in \textup{ind}_H^G(V)$ and $g \in G$, so $A_{\varphi \circ ev_1}(f) : g \mapsto \varphi (f(g))$ \textit{i.e.} $A_{\varphi \circ ev_1} = i \circ \textup{ind}_H^G(\varphi)$ if $i$ denotes $\textup{ind}_H^G(V) \subset \textup{Ind}_H^G(V)$. As a result $\varphi \in \textup{End}_H(V) \mapsto \textup{ind}_H^G(\varphi) \in \textup{End}_G(\textup{ind}_H^G(V))$ is an isomorphism of $R$-algebras. \end{proof}

When the condition $\textup{Hom}_H(\textup{ker}(ev_1),V) = 0$ holds, we have in particular a canonical isomorphism between $\textup{End}_H(V)$ and $\textup{End}_G(\textup{ind}_H^G(V))$. To refer to this very previse situation, we decorate isomorphims with curved arrows $\curvearrowright$ or $\curvearrowleft$ from $\textup{End}_H(V)$ to $\textup{End}_G(\textup{ind}_H^G(V))$. This means that $\textup{End}_H(V)$ acts on the set of ``images'' -- some would prefer to say the ``fiber'' -- of the representation $\textup{ind}_H^G(V)$ seen as a space of functions.

\section{Around the Bernstein center}

\subsection{Jacobson rings.} We are interested here in Jacobson (commutative) rings. By definition, they are the rings such that every prime ideal is the intersection of maximal ideals. In particular, their Jacobson radical -- which is the intersection of all maximal ideals -- agrees with their nilradical -- which is the set of nilpotent elements, or equivalently, the intersection of all prime ideals. Any finitely generated (commutative) algebra over a Jacobson ring is itself Jacobson. A field is Jacobson, and so is the integers $\mathbb{Z}$, but $\mathbb{Z}_\ell$ is not as its Jacobson radical is $\ell \mathbb{Z}_\ell$.

When $A$ is a Jacobson ring, the topological space $\textup{Spec}(A)$ is Jacobson -- this is even an equivalence \cite[\href{https://stacks.math.columbia.edu/tag/00G3}{Tag 00G3}]{stacks-project}. It ensures that closed points are somehow well-behaved with respect to subsets. For instance, if $X$ is a locally closed subset of $\textup{Spec}(A)$, a closed point $x$ in $X$ will be closed in $\textup{Spec}(A)$ \cite[\href{https://stacks.math.columbia.edu/tag/005X}{Tag 005X}]{stacks-project}. Denoting by $X_{\textup{max}}$ the set of closed points, we will have a natural identification $X_{\textup{max}} = X \cap \textup{Spec}(A)_{\textup{max}}$ and by the proof of \cite[\href{https://stacks.math.columbia.edu/tag/005X}{Tag 005X}]{stacks-project} we have $X_{\textup{max}} \neq \emptyset$ if $X$ is non-empty.

\begin{lemma} \label{reduced_Jacobson_rings_lemma} Let $A$ be a Jacobson reduced ring. Let $U$ be an open dense subset of $\textup{Spec}(A)$. An element of $A$ is determined by its specializations over $U_{\textup{max}}$ \textit{i.e.} we have an injective map: 
$$\begin{array}{ccc}
A & \to & \displaystyle \prod_{\mathfrak{m} \in U_{\textup{max}}} A/\mathfrak{m} \\
a & \mapsto & \displaystyle (a_\mathfrak{m})_{\mathfrak{m}} \end{array}.$$ \end{lemma}

\begin{proof} Consider the ideal $I = \cap_{\mathfrak{m} \in U_{\textup{max}}} \mathfrak{m}$, which is well-defined as $U_\textup{max} \neq \emptyset$. We want to show this ideal is the zero ideal. Let $f \in I$. By definition $D(f) = \textup{Spec}(A[1/f])$ is an open subset of $\textup{Spec}(A)$ and we have $D(f)_{\textup{max}} \cap U_{\textup{max}} = \emptyset$ by \cite[\href{https://stacks.math.columbia.edu/tag/00G6}{Tag 00G6}]{stacks-project}. This implies that $(D(f) \cap U)_{\textup{max}} = \emptyset$ and therefore $D(f) = \emptyset$ by density of $U$ \textit{i.e.} $f$ is nilpotent in $A$ by \cite[ Ex 2.2]{gortz_wedhorn_book}. So $f=0$ because $A$ is reduced and we obtain $I=0$ as claimed. \end{proof}

For all connected reductive groups $G$ over a non-archimedean local field $F$, the block decomposition of the center reads: 
$$\Z_\C(G) = \prod_{\mathfrak{s} \in \mathcal{B}_\C(G)} \Z_\C^\mathfrak{s}(G)$$
where $\mathcal{B}_\C(G)$ is the set of inertial classes and each local component $\Z_\C^\mathfrak{s}(G)$ is an integral domain that is finitely generated as a $\C$-algebra. In particular they are reduced Jacobson rings.

According to Lemma \ref{reduced_Jacobson_rings_lemma} above, if $M$ is a module over a reduced Jacobson ring $A$, it can be sufficient -- when $M$ is ``big'' enough -- to check the action of $A$ on any open dense subset to understand its action on $M$. We make this condition on $M$ more precis by defining a quotient support $\textup{QS}(M) = \{ \mathfrak{p} \in \textup{Spec}(A) \ | \ M \otimes_A A/ \mathfrak{p} \neq 0\}$ for the module $M$. Note that, by Nakayama's lemma, this agrees with the usual definition of support when $M$ is finitely generated. When $\mathfrak{m} \in \textup{Spec}(A)_{\textup{max}}$, an element $a \in A$ acts on $M \otimes_A A/\mathfrak{m}$ through a scalar $a_\mathfrak{m}(M) \in A/\mathfrak{m}$ in the center of $\textup{End}_{A/\mathfrak{m}}(M \otimes_A A/\mathfrak{m})$. Because $M \otimes_A A/\mathfrak{m}$ can be the zero module, we may have $a_\mathfrak{m}(M) = 0$ with this definition even though $a_\mathfrak{m} \neq 0$. The quotient support $\textup{QS}(M)$ is open in $\textup{Spec}(A)$ as its complement is easily seen to be closed.

\begin{corollary} \label{reduced_jacobson_big_modules_cor} Let $M$ be a module over a reduced Jacobson ring $A$. Assume that $\textup{QS}(M)$ is dense in $\textup{Spec}(A)$. Then for all open dense subsets $U$ of $\textup{Spec}(A)$, we have an injection:
$$\begin{array}{ccc}
A & \to & \displaystyle \prod_{\mathfrak{m} \in U_{\textup{max}}} A / \mathfrak{m} \\
a & \mapsto & \displaystyle (a_\mathfrak{m}(M))_{\mathfrak{m}} \end{array}.$$
 \end{corollary}

 \begin{proof} This is a simple application of Lemma \ref{reduced_Jacobson_rings_lemma} to the open dense set $U \cap \textup{QS}(M)$. \end{proof}

\subsection{Generic semi-simplicity.} Let $G$ be a connected reductive group over $F$. Let $R$ be an algebraically closed field of characteristic $\ell$ that is banal with respect to $G$ \textit{i.e.} $\ell$ does not divide the pro-order $|G|$ of the group.

\begin{lemma} \label{bernstein_center_components_lemma} The center of $\textup{Rep}_R(G)$ can be decomposed as a product over inertial classes: 
$$\Z_R(G) = \prod_{\mathfrak{s} \in \mathcal{B}_R(G)} \Z_R^{\mathfrak{s}}(G)$$
where each $\Z_R^{\mathfrak{s}}(G)$ is an integral domain and finite type $R$-algebra. \end{lemma}

\begin{proof} 
By using the results in \cite{vig_book} on the representation theory of $G$ over an algebraically closed field of banal characteristic, the methods of \cite{bernstein_deligne} can be extended to this setting. In particular, the description of these components as ring of invariants also holds.
\end{proof}

\begin{prop} \label{generic_semi_simplicity} There exists an open dense subset $U$ of $\Z_R(G)$ such that for all $\eta \in U_{\textup{max}}$ the category $\textup{Rep}_R^\eta(G)$ is semi-simple and has a unique simple object $\pi_\eta$. \end{prop}

\begin{proof} Considering the block decomposition of $\textup{Rep}_R(G)$, it is enough to prove it for each block. So let $\mathfrak{s} \in \mathcal{B}_R(G)$. Let $P=MN$ and $\sigma_\mathfrak{s} \in \textup{Rep}_R(M)$ be a parabolic and an irreducible supercuspidal associated to this inertial class. The representation $i_P^G(\sigma_\mathfrak{s} \Psi)$ is a pro-generator of the category $\textup{Rep}_R^\mathfrak{s}(G)$ where $\Psi : M \to R[M/M^0]$ is the universal unramified character for $M$. Similarly to \cite[Prop 3.14]{bernstein_deligne}, there exists a compact open subgroup $K$ in $G$ of invertible pro-order in $R$ and a non-zero $f \in \mathfrak{Z}_R^\mathfrak{s}(G)$ such that the $R$-algebra $\mathcal{H}_R(G,K)[1/f]$ is an Azumaya algebra over $\mathfrak{Z}_R(G)[1/f]$ of dimension $N$. Here $\mathfrak{Z}_R(G)[1/f] = \mathfrak{Z}_R^\mathfrak{s}(G)[1/f]$ because $f \in \mathfrak{Z}_R^\mathfrak{s}(G)$. Furthermore $\textup{Rep}_R^\mathfrak{s}(G)$ is naturally equivalent to the category of modules over $\mathcal{H}_R^\mathfrak{s}(G,K)$ where $\mathcal{H}_R^\mathfrak{s}(G,K)$ is a direct factor ring of $\mathcal{H}_R(G,K)$. So $\mathcal{H}_R(G,K)[1/f] = \mathcal{H}_R^\mathfrak{s}(G,K)[1/f]$. Now specializing this algebra to a character $\eta : \Z_R^\mathfrak{s}(G) \to R$ gives an equivalence of categories between $\textup{Rep}_R^\eta(G)$ and the category of modules over $\mathcal{H}_R(G,K)[1/f] \otimes_{\eta} R \simeq \mathcal{M}_N(R)$. The category $\mathcal{M}_N(R)$-mod is Morita equivalent to the category of $R$-vector spaces. So we obtained that $D(f) = \textup{Spec}(\Z_R^\mathfrak{s}(G)[1/f])$ is a non-empty open set in the irreducible variety $\textup{Spec}(\Z_R^\mathfrak{s}(G))$, therefore it is dense and for all $\eta \in U_{\textup{max}}$ the category $\textup{Rep}_R^\eta(G)$ is semi-simple with a single simple object $\pi_\eta = i_P^G(\sigma \Psi) \otimes_{\eta} R$ coming from the generic irreducibility. \end{proof}

\subsection{Regular representation.} We combine the previous two paragraphs to obtain ``generic'' properties about the regular representation. We carry on with the hypotheses with $G$ connected reductive group over $F$ and $R$ an algebraically closed field of banal characteristic with respect to $G$. For $V \in \textup{Rep}_R(G)$ and $\eta : \Z_R(G) \to R$ a character of the center, we recall that the largest $\eta$-quotient of $V$ is defined as $V_\eta = V \otimes_\eta R = V/V[\eta]$.

\begin{prop} \label{regular_representation_generic_quotient_semisimple_prop} There exists an open dense subset $U$ in $\textup{Spec}(\Z_R(G))$ such that for all characters $\eta : \Z_R(G) \to R$ in $U_{\textup{max}}$, we have:
$$C_c^\infty(G)_\eta = \pi_\eta \otimes_R \pi_\eta^\vee.$$ \end{prop}

\begin{proof} This is an easy application of Proposition \ref{generic_semi_simplicity} combined with the classical fact that $C_c^\infty(G)_\pi = \pi \otimes_R \pi^\vee$ as a $(G \times G)$-representation for all irreducible $\pi \in \textup{Rep}_R(G)$. \end{proof}

\subsection{Extension of scalars.} \label{extension_of_scalars_appendix_section} Let $R$ be a $\mathbb{Z}[1/p]$-algebra. We first introduce the Gelfand-Graev representations that will be our cornerstone for the compatibility of the Bernstein center with scalar extension. In this section, all tensor products are over $\ZZ[1/p]$ unless otherwise stated.

Let $N_n$ be the unipotent radical of the standard Borel \textit{i.e.} the group of unipotent upper triangular matrices in $G_n$. We consider the ring $R_0=\mathbb{Z}[1/p,\mu_{p^\infty}]$ that is obtained by adjoining all $p$-power roots of unity. Let $\psi$ be a non-degenerate character of $N_n$ with values in $R_0$. We define the Gelfand-Graev representation with coefficients in $R_0$ by $\textup{ind}_{N_n}^{G_n}(\psi)$. We introduce the term \textit{locally finitely generated} for a representation $V \in \textup{Rep}_R(G_n)$. It means, in terms of the depth decomposition $V = \bigoplus_r V_r$, that each $V_r \in \textup{Rep}_R^r(G_n)$ is finitely generated. By a \textit{local progenerator} we therefore mean a locally finitely generated projective generator of $\textup{Rep}_R(G_n)$. The forthcoming paper \cite{dhkm_banal} proves that:

\begin{prop} \label{dhkm_prop_gelfand_graev_integral} There exists an integral model $W_{N_n,\psi}$ of the representation of Gelfand-Graev over $\mathbb{Z}[1/p]$ such that $W_{N_n,\psi} \otimes R_0$ is isomorphic to $\textup{ind}_{N_n}^{G_n}(\psi)$. Furthermore $W_{N_n,\psi}$ is locally finitely generated and projective. \end{prop}

\begin{definition} The Gelfand-Graev representation over $R$ is defined as $W_{N_n,\psi}^R=W_{N_n,\psi} \otimes R$. \end{definition}

We are going to prove:

\begin{thm} \label{gelfand_graev_endo_isom_center_thm} The map $\Phi_R : z \in \Z_R(G_n) \mapsto z_{\W_{N_n,\psi}^R} \in \textup{End}_{G_n}(W_{N_n,\psi}^R)$ is an isomorphism. \end{thm}

The proof of the theorem breaks down into the following two lemmas, which easily implies on the one hand the surjectivty of $\Phi_R$ and on the other hand its injectivity. As these proofs require several steps, we prove them in a separate section:

\begin{lemma} \label{existence_of_a_section_Phi_A_lemma} There exists a section $\Psi_R : \textup{End}_{G_n}(W_{N_n,\psi}^R) \to \Z_R(G_n)$ of $\Phi_R$. \end{lemma}

\begin{lemma} \label{faithfulness_action_Z_W_n_psi_lemma} The natural action of $\Z_R(G_n)$ on $W_{N_n,\psi}^R$ is faithful. \end{lemma}

We now explain how to deduce the compatibility with scalar extension as a corollary of Theorem \ref{gelfand_graev_endo_isom_center_thm}. First, we observe that restricting to finite depth enables us to work over a module category. Indeed, for a connected reductive group $G$ over $F$ denote by $r_1, r_2, \dots$ the depth sequence as in Section \ref{DEPTH SECTION} and by $e_i^G$ the central idempotent associated to the depth $r_i$ category $\Rep_R(G)_{r_i}$ whose progenerator is $P(r_i)$. If $K$ is a compact open subgroup of $G$ that is small enough, the finitely generated projective representation $\ind_K^G(1_K)$ surjects on each factor of $P(r_i)$ so $e_i^G \ind_K^G(1_k)$ is a progenerator of the depth $r_i$ subcategory and has endomorphism ring $e_i^G \mathcal{H}_R(G,K) e_i^G$. Therefore we have an equivalence of categories:
$$\begin{array}{ccc}
    \Rep_R(G)_{r_i} & \to & (e_i^G \mathcal{H}_R(G,K) e_i^G)\textup{-mod} \\
    V & \mapsto & e_i^G V^K
\end{array}.$$
The following Lemma, whose proof appears at the end of this appendix, gives a fairly explicit interpretation of the compact $K$ for general linear groups if we group the depth pieces according to the ceiling function. 

\begin{lemma}\label{lem:equivalence}
Let $r \in \mathbb{N}$ and $K_r = I_n + \varpi_F^{r+1} \mathcal{M}_n(\mathcal{O}_F)$ be the $r$-th congruence subgroup in $G_n$. The functor of $K_r$-invariants induces an equivalence of categories:
$$\begin{array}{ccc}
    \Rep_R(G_n)_{\leq r} & \to & \mathcal{H}_R(G_n,K_r)\textup{-mod} \\
    V & \mapsto & V^{K_r}
\end{array}.$$
\end{lemma}

Because finitely generated projective modules are finitely presented, the endomorphism ring of the depth-$r$ summand of the Gelfand-Graev representation is compatible with arbitrary scalar extensions according to \cite[Prop I.2.13]{lam_projective}. This means that the functor $- \otimes R$ induces an isomorphism between pieces of depth at most $r$ as above:
$$\textup{End}_{\mathbb{Z}[1/p][G_n]}(W_{N_n,\psi,\leq r}) \otimes R \simeq \textup{End}_{R[G_n]}(W_{N_n,\psi,\leq r}^R).$$
This gives us an isomorphism $\Z_{\mathbb{Z}[1/p]}(G_n)_{\leq r} \otimes R \simeq \Z_R(G_n)_{\leq r}$ thanks to Theorem \ref{gelfand_graev_endo_isom_center_thm}. 

There is a more intrinsic way to describe this isomorphism. Consider the natural map:
$$\Z_{\mathbb{Z}[1/p]}(G_n) \to \Z_R(G_n)$$
induced by the forgetful functor $F : \textup{Rep}_R(G_n) \to \textup{Rep}_{\mathbb{Z}[1/p]}(G_n)$. To describe it explicitly, this ring morphism is $z \mapsto (z_{F(V)})_{V}$ where $V$ runs over all representations in $\textup{Rep}_R(G_n)$. Note that $z_{F(V)}$ is $R$-linear because multiplication by $a \in R$ is $\ZZ[1/p][G]$-linear. One only needs to check that for all $f \in \textup{Hom}_{R[G_n]}(V,V')$ we have $z_{F(V')} \circ f = f \circ z_{F(V)}$. But since $z$ is in $\Z_{\mathbb{Z}[1/p]}(G_n)$, the equality holds because $F(f)$ and $f$ are the same map. This natural map induces a bilinear map $\Z_{\mathbb{Z}[1/p]}(G_n) \times R \to \Z_R(G_n)$ which factors through the tensor product:
$$\eta_R : \Z_{\mathbb{Z}[1/p]}(G_n) \otimes_{\mathbb{Z}[1/p]} R \to \Z_R(G_n).$$
As a consequence of the present discussion, we have:

\begin{corollary} \label{scalar_extension_center_cor} The map $\eta_R$ induces an isomorphism
$$\eta_R : \Z_{\mathbb{Z}[1/p]}(G_n)_{\leq r} \otimes_{\mathbb{Z}[1/p]} R \to \Z_R(G_n)_{\leq r}.$$ \end{corollary}

\subsection{Proofs of Lemmas \ref{existence_of_a_section_Phi_A_lemma} and \ref{faithfulness_action_Z_W_n_psi_lemma}} Central to our approach is the construction of a local progenerator of $\textup{Rep}_R(G_n)$ out of the Gelfand-Graev representation. First of all, we state a general result about compatibility of progenerators with scalar extension and faithfully flat descent. In the lemma below $G$ is a connected reductive group:

\begin{lemma} \label{reduction_progenerator_lemma} Let $R \to S$ be a morphism of $\ZZ[1/p]$-algebras and $P \in \textup{Rep}_R(G)$.
\begin{itemize}
    \item[\textup{(i)}] If $P$ is a local progenerator of $\textup{Rep}_R(G)$, then $P \otimes_R S$ is a local progenerator of $\textup{Rep}_S(G)$.
    \item[\textup{(ii)}] If $P \otimes_R S$ is a local progenerator of $\textup{Rep}_S(G)$ and $S$ is faithfully flat over $R$, then $P$ is a local progenerator of $\textup{Rep}_R(G)$.
\end{itemize}
\end{lemma}

\begin{proof} The depth decomposition $P=\bigoplus P_r$ is compatible to scalar extension in the sense that the depth $r$ factor of $P \otimes_R S$ in $\textup{Rep}_S(G)$ is given by $(P \otimes_R S)_r=P_r \otimes_R S$. By the discussion preceding Lemma \ref{lem:equivalence} the result is equivalent to some central base change statements in module theory. We did not find a reference in terms of descent for finite projective modules over non-commutative rings, so we review the proof now which is very similar to the commutative case. Let $A$ be an $R$-algebra (not necessarily commutative), so that $R$ is central in $A$
and the categories at stake will be $A-\textup{mod}$ and $(A \otimes_R S)-\textup{mod}$ obtained by central base change.

\noindent (i) It is clear that finitely generated modules are preserved by central base change. Projectivity is preserved as well because, by the tensor-hom adjunction, the functors $\Hom_{A \otimes_R S} ( P\otimes_R S , - )$ and $\Hom_A ( P, \Hom_S(S,-) )$ are canonically isomorphic, where the latter is exact as $\Hom_S(S,-)$ is the forgetful functor $(A \otimes_R S)-\textup{mod} \to A-\textup{mod}$. The fact that $\Hom_{A \otimes_R S} ( P\otimes_R S , - )$ is faithful can be easily seen using again the tensor-hom adjunction as it becomes the composition of two faithful functors. So $P\otimes_R S$ is a generator.

\noindent (ii) This claim is similar to the descent for finite projective modules. We first prove that $P$ is finitely presented. Because $P \otimes_R S$ is finitely generated and projective, it is a direct factor of some finite free module \textit{i.e.} there exists $n \in \mathbb{N}$ such that $(A \otimes_R S)^n \simeq (P \otimes_R S) \oplus P'$. Therefore $P'$ is finitely generated, so $P \otimes_R S$ is finitely presented. As $P \otimes_R S$ is finitely generated, there is a finite family $y_1, \dots, y_s$ generating it and we can write each of them as finite sum $y_i = \sum_j x_{i,j} \otimes_R s_{i,j}$. We obtain a finite set made of the $(x_{i,j})_{i,j}$ that we reorder as $x_1, \dots, x_t$. Consider the map $f : (\alpha_i)_i \in A^s \mapsto \sum \alpha_i x_i \in P$. Then its central base change to $S$ is surjective, so $f$ is surjective as well by faithfully flatness. A similar argument proves that its kernel is finitely generated, that is $P$ is finitely presented as an $A$-module.

Next, suppose $V \to W$ is a surjective map of $A$-modules. Let $C$ denote the cokernel of $$\Hom_A(P,V)\to \Hom_A(P,W).$$
Since $P$ is finitely presented and $R \to S$ is flat, we have by \cite[Prop 2.13]{lam_projective} an isomorphism: 
$$\Hom_A(P,V)\otimes_R S \to \Hom_{A \otimes_R S}((P\otimes_R S,V \otimes_R S)$$
and similarly for $W$. Thus $C\otimes_R S$ is the cokernel of the map:
$$\Hom_{A \otimes_R S}(P\otimes_R S , V\otimes_R S)\to \Hom_{A \otimes_R S}(P\otimes_R S,W\otimes_R S)$$
which is zero by projectivity of $P\otimes_R S$. Hence $C$ is zero by faithful flatness of $R \to S$. Therefore we obtain that $P$ is projective.

We now turn to showing that finitely presented and flat implies projective for $A$-modules. Because $P$ is finitely presented, there exists an exact sequence $0 \to P' \to A^s \to P\to 0$ with $s \in \mathbb{N}$ and $P'$ generated by a finite set $c_1,\dots, c_t$ in $A^s$. By the equational criterion of flatness \cite[Th 4.23]{lam_lectures} there exists a homomorphism $\theta: A^s \to P'$ such that $\theta(c_i) = c_i$ for all $i$. Thus $\theta$ is the identity on $P'$ and defines a section of $P' \to A^s$, so the exact sequence splits. We conclude that $P$ is a direct summand of the finite free module $A^s$ \textit{i.e.} $P$ is projective.

For faithfulness, we need to show that the functor $\Hom_A(P,-)$ is faithful. By tensor-hom adjunction, and again by the fact that $P$ is finitely presented and $R \to S$ is flat, the composition $\Hom_A(P,-) \otimes_R S$ is canonically identified with $\Hom_{A \otimes R}(P \otimes_R S,- \otimes_R S)$. Our hypothesis ensures the latter functor is faithful. Now if a composition of functors is faithful, the first functor in the composiiton must be faithful so $\Hom_A(P,-)$ is faithful. This concludes the proof. \end{proof}

The first part of the previous lemma reduces the problem to constructing a progenerator with coefficients in $\ZZ[1/p]$. The second part tells us that it is enough to verify the conditions of being a local progenerator over a faithfully flat extension such as $R_0$. Here is the progenerator we mentioned for general linear groups:

\begin{lemma} \label{progenerator_built_from_GG_lem} Let $\underline{\textup{std}}$ be the set of standard parabolic subgroups of $G_n$ and define:
$$W_{N_n,\psi}^{\textup{gen}}=\bigoplus_{P \in \underline{\textup{std}}} \textup{ind}_{\bar{P}}^{G_n} \circ \textup{res}_{G_n}^{\bar{P}}(W_{N_n,\psi}).$$
It is a local progenerator of $\textup{Rep}_{\mathbb{Z}[1/p]}(G_n)$ and therefore we have:
$$\Z_{\mathbb{Z}[1/p]}(G_n) = \Z(\textup{End}_{G_n}(W_{N_n,\psi}^{\textup{gen}})).$$ \end{lemma}

\begin{proof} According to Lemma \ref{reduction_progenerator_lemma}, it is enough to prove it after a faithfully flat base change as induction and restriction functors commute to scalar extension. Therefore it is enough to prove it over $R_0$ replacing $W_{N_n,\psi}$ by $\textup{ind}_{N_n}^{G_n}(\psi)$. The representation $\textup{ind}_{N_n}^{G_n}(\psi)$ is locally finitely generated and projective by Proposition \ref{dhkm_prop_gelfand_graev_integral}. Note that the induction (and restriction) functors preserve finite generation over arbitrary rings, as opposed to what is written in \cite[Cor 1.5]{dhkm_finiteness}: indeed the noetherianity hypothesis there is superfluous as it relies on the proof of \cite[Lem 4.6]{dat_finitude} and second adjunction\footnote{one of the reasons that could explain it: the proof of the second adjunction in \cite{dhkm_finiteness} initially dealt with noetherian $\mathbb{Z}[1/p]$-algebras and they relaxed the noetherianity assumption at a very late stage of the writing.}, which are both valid over any $\mathbb{Z}[1/p]$-algebra. So $\textup{ind}_{N_n}^{G_n}(\psi)^{\textup{gen}}$ is locally finitely generated and projective. We now prove it is a generator. As all finitely generated objects admits a simple quotient, it is enough to prove that for all simple objects $\pi \in \textup{Rep}_{R_0}(G_n)$: $$\textup{Hom}_{R_0[G_n]}(\textup{ind}_{N_n}^{G_n}(\psi)^{\textup{gen}},\pi) \neq 0.$$
Actually $\pi$ has coefficients in a residue field of $R_0$. Denoting by $\mathcal{P} \in \textup{Spec}(R_0)$ the prime ideal kernel of $R_0 \to \textup{End}_{R_0[G_n]}(\pi)$ given by the action of scalars, we have $k(\mathcal{P}) = \textup{Frac} (R_0 / \mathcal{P})$ and $\pi \in \textup{Rep}_{k(\mathcal{P})}(G_n)$. Therefore by tensor-hom adjunction:
$$\textup{Hom}_{R_0[G_n]}(\textup{ind}_{N_n}^{G_n}(\psi)^{\textup{gen}},\pi) \simeq \textup{Hom}_{k(\mathcal{P})[G_n]}(\textup{ind}_{N_n}^{G_n}(\psi)^{\textup{gen}} \otimes_{R_0} k(\mathcal{P}) , \pi).$$
So we reduce the question to checking that $\textup{ind}_{N_n}^{G_n}(\psi)^{\textup{gen}} \otimes_{R_0} k(\mathcal{P}) = \textup{ind}_{N_n}^{G_n}(\psi_{k(\mathcal{P})})^{\textup{gen}}$ is a local progenerator of the category $\textup{Rep}_{k(\mathcal{P})}(G_n)$. Note that field extensions are faithfully flat, so we can always assume that our base field $k$ is algebraically closed. Now we obtain that $\textup{ind}_{N_n}^{G_n}(\psi_k)^{\textup{gen}}$ is a local progenerator as a consequence of the following three properties: 
\begin{itemize} 
\item the existence of cuspidal support;
\item all cuspidals are generic for general linear groups; \item the restriction $r_{G_n}^{\bar{P}}(\textup{ind}_{N_n}^{G_n}(\psi_k)) \simeq \textup{ind}_{N_M}^M(\psi_k|_{N_M})$ with $P$ of Levi $M$ and $N_M = N_n \cap M$.
\end{itemize} 
So $\textup{Hom}_{R_0[G_n]}(\textup{ind}_{N_n}^{G_n}(\psi)^{\textup{gen}},\pi) \neq 0$ and Lemma \ref{reduction_progenerator_lemma} implies that $W_{N_n,\psi}^{\textup{gen}}$ is a local progenerator by faithfully flat descent. \end{proof}

Let $E_P \in \textup{End}_{G_n}(W_{N_n,\psi}^{\textup{gen}})$ be the projection on the direct factor associated to $P \in \underline{\textup{std}}$. It is an idempotent and the commutator in $\textup{End}_{G_n}(W_{N_n,\psi}^{\textup{gen}})$ of all these idempotents is:
$$\bigoplus_{P \in \underline{\textup{std}}} \textup{End}_{G_n}(\textup{ind}_{\bar{P}}^{G_n} \circ \textup{res}_{G_n}^{\bar{P}}(W_{N_n,\psi})).$$
Therefore we have inclusions:
$$Z(\textup{End}_G(W_{N_n,\psi}^{\textup{gen}})) \subset \bigoplus_{P \in \underline{\textup{std}}} \textup{End}_{G_n}(\textup{ind}_{\bar{P}}^{G_n} \circ \textup{res}_{G_n}^{\bar{P}}(W_{N_n,\psi})) \subset \textup{End}_{G_n}(W_{N_n,\psi}^{\textup{gen}}).$$
We denote by $(M_P)_{P \in \underline{\textup{std}}}$ an element in the second term of the inclusions. Consider the map $\textup{ev}_{G_n} : (M_P)_{P \in \underline{\textup{std}}} \mapsto M_{G_n}$ and denote by $\Phi$ its restriction to $Z(\textup{End}_{G_n}(W_{N_n,\psi}^\textup{gen}))$. The functor $F=\bigoplus_{P \in \underline{\textup{std}}} \textup{ind}_{\bar{P}}^{G_n} \circ \textup{res}_{G_n}^{\bar{P}}$ induces a morphism:
$$\Psi_F : \phi \mapsto (\phi_P)_{P \in \underline{\textup{std}}} = (\textup{ind}_{\bar{P}}^{G_n} \circ \textup{res}_{G_n}^{\bar{P}}(\phi))_{P \in \underline{\textup{std}}}.$$
We want to show:

\begin{lemma} \label{commutative_diagram_Z[1/p]_lemma} The image of $\Psi_F$ is central in $\textup{End}_G(W_{N_n,\psi}^{\textup{gen}})$ \textit{i.e.} there exists a section $\Psi$ of $\Phi$ completing the commutative diagram of solid arrows:
$$\xymatrix{
\textup{End}_{G_n}(W_{N_n,\psi}) \ar[rd]^{\Psi_F} \ar@{-->}[r]^{\Psi}  & Z(\textup{End}_{G_n}(W_{N_n,\psi}^{\textup{gen}})) \ar@{^{(}->}[d]  \ar[r]^{\Phi} & \textup{End}_{G_n}(W_{N_n,\psi}) \\
  & \displaystyle \bigoplus_{P \in \underline{\textup{std}}} \textup{End}_{G_n}(\textup{ind}_{\bar{P}}^{G_n} \circ \textup{res}_{G_n}^{\bar{P}}(W_{N_n,\psi})) \ar@{^{(}->}[d] \ar[ru]^{\textup{ev}_{G_n}} & \\
& \textup{End}_{G_n}(W_{N_n,\psi}^{\textup{gen}}) & 
}.$$
Morevoer $\Phi$ is injective, therefore both $\Psi$ and $\Phi$ are isomorphisms. \end{lemma} 

\begin{proof} We start by the injectivity of $\Phi$. As $W_{N_n,\psi}^\textup{gen}$ is a local progenerator, we know that the map $z \in \Z_{\mathbb{Z}[1/p]}(G_n) \mapsto z_{W_{N_n,\psi}^\textup{gen}} \in Z(\textup{End}_G(W_{N_n,\psi}^\textup{gen}))$ is an isomorphism, and composing by $\Phi$ we obtain $z \in \Z_{\mathbb{Z}[1/p]}(G_n) \mapsto z_{W_{N_n,\psi}} \in \textup{End}_G(W_{N_n,\psi})$. Similarly to \cite{h_whitt}, extending scalars to $\C$ gives the injectivity of $\Phi$.

There remains to prove the existence of the section $\Psi$. The composition $\Psi_F \circ \Phi$ induces:
$$z \in  \Z_{\mathbb{Z}[1/p]}(G_n) \mapsto ( \textup{ind}_{\bar{P}}^{G_n} \circ \textup{res}_{G_n}^{\bar{P}}(z_{W_{N_n,\psi}}))_{P \in \underline{\textup{std}}} \in \bigoplus_{P \in \underline{\textup{std}}} \textup{End}_{G_n}(\textup{ind}_{\bar{P}}^{G_n} \circ \textup{res}_{G_n}^{\bar{P}}(W_{N_n,\psi})).$$
We want to prove that the image of $z$ is central. First the action of $z \in \Z_{\mathbb{Z}[1/p]}(G_n)$ on $W_{N_n,\psi}^\textup{gen}$ is given by $(z_{\textup{ind}_{\bar{P}}^{G_n} \circ \textup{res}_{G_n}^{\bar{P}}(W_{N_n,\psi})})_{P \in \underline{\textup{std}}} \in Z(\textup{End}_{G_n}(W_{N_n,\psi}^\textup{gen}))$. The centrality will be a clear consequence of the following identity:
$$z_{\textup{ind}_{\bar{P}}^{G_n} \circ \textup{res}_{G_n}^{\bar{P}}(W_{N_n,\psi})} = \textup{ind}_{\bar{P}}^{G_n} \circ \textup{res}_{G_n}^{\bar{P}}(z_{W_{N_n,\psi}}).$$
This identity comes from the existence of Harish-Chandra morphisms \cite[Th 4.1]{dhkm_finiteness}. Actually the only property we use, which is weaker than the full Harish-Chandra morphisms, is the fact that $z_{\textup{ind}_{\bar{P}}^{G_n} \circ \textup{res}_{G_n}^{\bar{P}}(W_{N_n,\psi})} = \textup{ind}_{\bar{P}}^{G_n} (f')$ for some $f' \in \textup{End}_M(\textup{res}_{G_n}^{\bar{P}}(W_{N_n,\psi}))$. By Frobenius reciprocity: 
$$f' = \textup{res}_{G_n}^{\bar{P}}(z_{W_{N_n,\psi}})$$
and therefore the identity $z_{W_{N_n,\psi}^\textup{gen}} = \Psi_F \circ \Phi (z_{W_{N_n,\psi}^\textup{gen}})$ holds. So $\Psi_F$ induces the required $\Psi$, which is at the same time injective and a section of $\Phi$. Hence $\Psi$ and $\Phi$ are isomorphisms. \end{proof}

We can extend scalars to $R$ to obtain a new diagam. However we first have to bound the depth in order to avoid complications due to non commutation of infinite direct product and tensor product. As finitely generated projective objects are direct factors of some finite free module, this allows arbitrary scalar extension for their endomorphism rings. In our situation, we obtain a commutative diagram:
$$\xymatrix{
\textup{End}_{G_n}(W_{N_n,\psi,\leq r}^R) \ar[rd]^{\Psi_{F_R}} \ar[r]^{\Psi \otimes R}  & Z(\textup{End}_{G_n}(W_{N_n,\psi,\leq r}^{\textup{gen}})) \otimes R \ar[d]  \ar[r]^{\Phi \otimes R} & \textup{End}_{G_n}(W_{N_n,\psi,\leq r}^R) \\
  & \displaystyle \bigoplus_{P \in \underline{\textup{std}}} \textup{End}_{G_n}(\textup{ind}_{\bar{P}}^{G_n} \circ \textup{res}_{G_n}^{\bar{P}}(W_{N_n,\psi, \leq r}^R)) \ar@{^{(}->}[d] \ar[ru]^{\textup{ev}_{G_n}} & \\
& \textup{End}_{G_n}((W_{N_n,\psi, \leq r}^R)^{\textup{gen}}) & 
}$$
where the endofunctors $F_R = \bigoplus_{P \in \underline{\textup{std}}} \textup{ind}_{\bar{P}}^{G_n} \circ \textup{res}_{G_n}^{\bar{P}}$ and $F_{\ZZ[1/p]} \otimes R$ of $\textup{Rep}_R(G_n)$ are canonically isomorphic because parabolic restriction and induction functors commute with scalar extension. Also $\Psi \otimes R$ and $\Phi \otimes R$ are still inverse isomorphisms and $\Psi_{F_R}$ remains a section of $\textup{ev}_{G_n}$.

Because $\textup{End}_{G_n}(W_{N_n,\psi, \leq r}^\textup{gen}) \otimes R$ identifies with $\textup{End}_{G_n}((W_{N_n,\psi, \leq r}^R)^\textup{gen})$, the image of the first vertical map must lie in $Z(\textup{End}_{G_n}((W_{N_n,\psi, \leq r}^R)^\textup{gen}))$. In other words, we can complete the diagram with a section $\Psi_R$ and a retraction $\Phi_R$, coming respectively from $\Psi \otimes R$ and $\textup{ev}_{G_n}$, into the following:
$$\xymatrix{
\textup{End}_{G_n}(W_{N_n,\psi, \leq r}^R) \ar[rdd]^{\Psi_{F_R}} \ar[rd]^{\Psi_R} \ar[r]^{\Psi \otimes R}  & Z(\textup{End}_{G_n}(W_{N_n,\psi, \leq r}^{\textup{gen}})) \otimes R \ar[d]  \ar[r]^{\Phi \otimes R} & \textup{End}_{G_n}(W_{N_n,\psi, \leq r}^A) \\
 & Z(\textup{End}_{G_n}((W_{N_n,\psi, \leq r}^R)^{\textup{gen}})) \ar[ru]^{\Phi_R} \ar@{^{(}->}[d] & \\
  & \displaystyle \bigoplus_{P \in \underline{\textup{std}}} \textup{End}_{G_n}(\textup{ind}_{\bar{P}}^{G_n} \circ \textup{res}_{G_n}^{\bar{P}}(W_{N_n,\psi, \leq r}^R)) \ar@{^{(}->}[d] \ar[ruu]^{\textup{ev}_{G_n}} & \\
& \textup{End}_{G_n}((W_{N_n,\psi, \leq r}^R)^{\textup{gen}}) & 
}.$$
In particular composing from left to right implies Lemma \ref{existence_of_a_section_Phi_A_lemma}.

There remains to prove Lemma \ref{faithfulness_action_Z_W_n_psi_lemma} to have the compatibility with arbitrary scalar extension:
$$\Z_{\mathbb{Z}[1/p]}(G_n)_{\leq r} \otimes R \simeq \Z_R(G_n)_{\leq r}.$$
Note that the center $\Z_R(G_n)$ acting faithfully on $W^R_{N_n,\psi}$ is equivalent to the injectivity of $\Phi_R$. In the course of the proof of Lemma \ref{commutative_diagram_Z[1/p]_lemma}, we proved the injectivity of $\Phi$ using an identity that was a consequence of the existence of Harish-Chandra morphisms:
$$z_{\textup{ind}_{\bar{P}}^{G_n} \circ \textup{res}_{G_n}^{\bar{P}}(W_{N_n,\psi})} = \textup{ind}_{\bar{P}}^{G_n} \circ \textup{res}_{G_n}^{\bar{P}}(z_{W_{N_n,\psi}}).$$
for $z \in \Z_{\mathbb{Z}[1/p]}(G_n)$. If such an identity held for $z \in \Z_R(G_n)$ and $W_{N_n,\psi}^R$, we would be able to conclude as above that $\Psi_R \circ \Phi_R (z_{(W_{N_n,\psi}^R)^\textup{gen}}) = z_{(W_{N_n,\psi}^R)^\textup{gen}}$. In this case $\Psi_R$ would be an isomorphism and so would $\Phi_R$. Therefore we focus our efforts on proving this identity:
\begin{prop} \label{induced_morphisms_a_la_harish_chandra_prop} Let $z \in \Z_R(G_n)$ and $P \in \underline{\textup{std}}$. Then:
$$z_{\textup{ind}_{\bar{P}}^{G_n} \circ \textup{res}_{G_n}^{\bar{P}}(W_{N_n,\psi})} = \textup{ind}_{\bar{P}}^{G_n} \circ \textup{res}_{G_n}^{\bar{P}}(z_{W_{N_n,\psi}}).$$ \end{prop}

\begin{proof} As in the proof of Lemma \ref{commutative_diagram_Z[1/p]_lemma}, it is sufficient to prove that:
$$z_{\textup{ind}_{\bar{P}}^{G_n} \circ \textup{res}_{G_n}^{\bar{P}}(W_{N_n,\psi})} =\textup{ind}_{\bar{P}}^{G_n}(f') \textup{ for some } f' \in \textup{End}_M(\textup{res}_{G_n}^{\bar{P}}(W_{N_n,\psi}^R)).$$ 
The methods of Bushnell-Henniart \cite{bh_whitt,dhkm_banal} give that $\textup{res}_{G_n}^{\bar{P}}(W_{N_n,\psi}^R)$ is the Gelfand-Graev representation for the Levi $M$ and we denote it by $\sigma$. We want to prove that:
$$Z(\textup{End}_{G_n}(\textup{ind}_{\bar{P}}^{G_n}(\sigma))) \subset \textup{ind}_{\bar{P}}^{G_n}(\textup{End}_M(\sigma)).$$
This will be a consequence of the following proposition:
\begin{prop} We have:
$$\textup{End}_{G_n \times M}(\textup{ind}_{\bar{P}}^{G_n}(C_c^\infty (M))) \overset{\curvearrowleft}{\simeq} \textup{End}_{M \times M}(C_c^\infty(M)) = \Z_R(M).$$ \end{prop}

\begin{proof} By Frobenius reciprocity we have:
$$\textup{End}_{G_n \times M}(\textup{ind}_{\bar{P}}^{G_n}(C_c^\infty (M))) = \textup{Hom}_{M \times M}(\textup{res}_{G_n}^{\bar{P}} \circ \textup{ind}_{\bar{P}}^{G_n}(C_c^\infty (M)),C_c^\infty(M)).$$ 
By the geometric lemma stated in Proposition \ref{geometric_lemma_prop}, the restriction-induction has a filtration by certain $(M \times M)$-modules where the subquotients $I_w$ are indexed by $w \in W^{M,M}$ in the Weyl group. As in Proposition \ref{endomorphisms_of_W_k_nm_prop}, the result of the current proposition will hold as long as:
$$\textup{Hom}_{M \times M}(I_w,C_c^\infty(M))= 0 \textup{ for all } w \neq \textup{Id}.$$
When $w \neq \textup{Id}$ we have isomorphisms of $(M \times M)$-modules:
$$I_w \simeq \bar{\mathfrak{i}}_{M \cap w(M)}^M \left( \delta_w \otimes_R ( w \circ \bar{\mathfrak{r}}_M^{w^{-1}(M) \cap M}(C_c^\infty(M))) \right)$$
where the bar accounts for the fact that are our standard parabolics are opposite to the usual upper triangular ones. By second adjunction $\textup{Hom}_{M \times M}(I_w,C_c^\infty(M))$ is isomorphic to:
$$\textup{Hom}_{(M \cap w(M)) \times M}(\delta_w \otimes_R ( w \circ \bar{\mathfrak{r}}_M^{w^{-1}(M) \cap M}(C_c^\infty(M))),\mathfrak{r}_M^{M \cap w(M)} (C_c^\infty(M))).$$
On the one hand, by integrating over the group $M \cap w(N_n)$ we get that the right-hand side is $C_c^\infty((M \cap w(N_n)) \backslash M)$ as a $((M \cap w(M)) \times M)$-module. An element $(m_w,m)$ will act on a function $f$ in this space by:
$$(m_w,m) \cdot f : (M \cap w(N_n)) m' \mapsto \delta_{M \cap w(N_n)}(m_w) \times f( (M \cap w(N_n)) m_w^{-1} m' m).$$
On the other hand, the left-hand side is $\delta_w \otimes_R (w \circ (C_c^\infty((w^{-1}(\bar{N}_n) \cap M) \backslash M)))$.

Let $z \in Z(M)$. On the right-hand side, the element: 
$$(1,z) - \delta_{M \cap w(N_n)}(z^{-1}) \cdot (z^{-1},1) \in R[(M \cap w(M)) \times M]$$ 
acts as zero. On the left-hand side however, there exists a character $\chi$ such that the element: 
$$(1,z) - \chi(z) \cdot (z^{-w},1)$$
acts as zero. We are not going to make this character $\chi$ explicit as we can carry out our argument for all characters. As a result we see that any morphism in $\textup{Hom}_{M \times M}(I_w,C_c^\infty(M))$ must factor through elements of the form $(\delta_{M \cap w(N_n)}(z^{-1}) \cdot (z,1) - \chi(z^{-1}) \cdot (z^w,1)) \cdot f$. We denote by $I_w'$ the quotient of $I_w$ by the previous elements. In particular the group:
$$H=\{ (z^{-1} z^w,1) | z \in Z(M) \}$$
must act through a character on $I_w'$. 

Let $w \neq \textup{Id}$. Let $z \in Z(M)$ such that $z' = z^{-1} z^w$ is not a compact element. Note that such an element exist because $w \neq \textup{Id}$. Then for compact support reasons we must have $\textup{Hom}_{M \cap w(M)}(I_w',C_c^\infty((w^{-1}(\bar{N}_n) \cap M) \backslash M))) = 0$ as $z'$ acts as a character on the left-hand side but can not on the right-and side because $(w^{-1}(\bar{N}_n) \cap M) H$ is not compact. Therefore $\textup{Hom}_{M \times M}(I_w,C_c^\infty(M))= 0$ and the proposition holds. \end{proof} 

In order to finish the proof, note that $C_c^\infty(M) \twoheadrightarrow \textup{ind}_{N_n \cap M}^M(\psi_{N_n \cap M})=\sigma$ via the largest $\psi^{-1}$-isotypic quotient construction. This map is $(M \times (N_n \cap M))$-equivariant. In particular for $z \in \Z_R(G_n)$ this yields a commutative diagram:
$$\xymatrixcolsep{5pc}\xymatrix{
		\textup{ind}_{\bar{P}}^{G_n}(C_c^\infty(M)) \ar@{->>}[r] \ar[d]^{z_{\textup{ind}_{\bar{P}}^{G_n}(C_c^\infty(M))}} & \textup{ind}_{\bar{P}}^{G_n}(\sigma) \ar[d]^{z_{\textup{ind}_{\bar{P}}^{G_n}(\sigma)}} \\
		\textup{ind}_{\bar{P}}^{G_n}(C_c^\infty(M))  \ar@{->>}[r] & \textup{ind}_{\bar{P}}^{G_n}(\sigma)
		}.$$
By Proposition \ref{induced_morphisms_a_la_harish_chandra_prop} the map $z_{\textup{ind}_{\bar{P}}^{G_n}(C_c^\infty(M))} = \textup{ind}_{\bar{P}}^{G_n}(z'_{C_c^\infty(M)})$ for some $z' \in \Z_R(M)$, so this implies that $z_{\textup{ind}_{\bar{P}}^{G_n}(\sigma)} = \textup{ind}_{\bar{P}}^{G_n}(z'_{\sigma})$ and the proposition holds because $z_{\textup{ind}_{\bar{P}}^{G_n}(\sigma)}$ is induced. \end{proof}

\subsection{Proof of Lemma~\ref{lem:equivalence}}
Consider $Q_{\leq r} = \textup{ind}_{K_r}^{G_n}(1_{K_r})$ where $1_{K_r} = R$ is the trivial module. The lemma will hold as a Morita equivalence statement because $\mathcal{H}_R(G_n,K_r) = \textup{End}_{R[G_n]}(Q_{\leq r})$, so we have to prove that $Q_{\leq r} \in \Rep_R(G_n)_{\leq r}$ is a progenerator of this category. First $Q_{\leq r}$ is:
\begin{itemize}
    \item finitely generated -- it is even cyclic as the characteristic function of $K_r$ generates it;
    \item projective -- because $\textup{ind}_{K_r}^{G_n}$ is left adjoint of the restrcition functor, which is exact.
\end{itemize}
We now prove it is generating the category $\Rep_R(G_n)_{\leq r}$ \textit{i.e.} all irreducible $\pi \in \Rep_R(G_n)_{\leq r}$ is a quotient of $Q_{\leq r}$. Say $\pi$ has depth $r_i \leq r$. Then, as a consequence of the existence unrefined minimal $K$-types \cite{moy_prasad_unrefined,vig_book,dat_finitude}, there exists a point $x$ in the Bruhat-Tits building of $G_n$ such that $\pi$ has a non-trivial fixed vector under $G_{x,r_i^+}$, where $G_{x,-}$ is the Moy-Prasad filtration. We can always assume that the point $x$ belongs to the star of $x_0$, whose stabilizer is $G_{x_0} = \textup{GL}_n(\mathcal{O}_F)$. For all integers $s$, we have $G_{x_0,s^+} \subset G_{x,s^+}$ by \cite[Lemma 4.3]{ABPS_sln}. In particular $G_{x_0,r^+} = K_r \subset G_{x,r^+} \subset G_{x,r_i^+}$ and $\pi$ has a vector fixed under $K_r$. Therefore $Q_{\leq r}$ is a generator.

We also have to prove that $Q_{\leq r}$ has no piece of depth strictly greater than $r$. Take $r_i > r$ and look at the central idempotent $e_{r_i} \in \Z_R(G_n)$ associated to the depth $r_i$ direct factor. Because $Q_{\leq r}$ is cyclic its depth $r_i$ factor is cyclic as well. In particular it is finitely generated. In the abelian category $\Rep_R(G_n)$, a non-zero finitely generated object will admit a non-zero simple quotient \cite[A.VI.3 Prop]{ren}. Let $\pi \in \Rep_R(G_n)_{r_i}$ be irreducible and use adjunction:
$$\textup{Hom}_{G_n}(Q_{\leq r},\pi) = \textup{Hom}_{K_r}(1_{K_r},\pi) = \textup{Hom}_{K_r}(1_{K_r},\pi^{K_r}).$$
We are now interested in $\pi^{K_r}$ and we want to prove it is zero. If it was non-zero, then there would exist an unrefined minimal $K$-type $\psi$ of depth $\leq r$ contained in $\pi$. But all unrefined minimal $K$-types of $\pi$ have same depth $r_i$, so we obtain a contradiction. Therefore $\pi$ has no fixed vectors under $K_r$ and this proves that $Q_{\leq r} \in \Rep_R(G_k)_{\leq r}$.

\bibliographystyle{alpha}
\bibliography{lesrefer}

\end{document}